\documentclass[onefignum,onetabnum]{siamart190516}

\usepackage{amsmath,amssymb,amsxtra}
\usepackage{mathrsfs}
\usepackage{graphicx}
\usepackage{float,verbatim}
\usepackage{threeparttable,booktabs}
\usepackage{underscore}
\usepackage{algorithm}
\usepackage{subfigure}
\usepackage{extarrows}
\usepackage{multirow}
\usepackage{makecell}
\usepackage[noend]{algpseudocode}
 \usepackage{color}
\graphicspath{{./pics/}}
\newtheorem{remark}{Remark}[section]
\newtheorem{example}{Example}[section]
\newtheorem{assumption}{Assumption}[section]

\renewcommand{\d}{\mathrm{d}}

\def\Om {\Omega}

\def\to{\rightarrow}

\def\p {\partial}
\def\tribar{|\!|\!|}
\newcommand{\dx}{\,{\rm d}x}
\newcommand{\dt}{\,{\rm d}t}
\newcommand{\ds}{\,{\rm d}s}
\newcommand{\dd}{\,{\rm d}}
\newcommand{\tV}{\widetilde{V}}
\newcommand{\tA}{\widetilde{A}}
\newcommand{\tB}{\widetilde{B}}

\allowdisplaybreaks

\title{Direct Algorithms for Reconstructing Small Conductivity Inclusions in Subdiffusion\thanks{The work of B. Jin is supported by Hong Kong RGC General Research Fund (Project 14306423), and a start-up fund from The Chinese University of Hong Kong.}}

\author{Jiho Hong\thanks{Department of Mathematics, The Chinese University of Hong Kong, Shatin, N.T., Hong Kong  SAR, P.R. China (\texttt{jihohong@cuhk.edu.hk, b.jin@cuhk.edu.hk, bangti.jin@gmail.com})} \and
Bangti Jin\footnotemark[2]
\and Zhizhang Wu\thanks{Department of Mathematics, The University of Hong Kong, Pokfulam Road, Hong Kong SAR, P.R. China (\texttt{wzz14@tsinghua.org.cn}).}}

\date{\today}

\begin{document}

\maketitle

\begin{abstract}
The subdiffusion model that involves a Caputo fractional derivative in time is widely used to describe anomalously slow diffusion processes.
In this work we aim at recovering the locations of small conductivity inclusions in the model from boundary measurement, and develop novel direct algorithms based on the asymptotic expansion of the boundary measurement with respect to the size of the inclusions and approximate fundamental solutions. These algorithms involve only algebraic manipulations and are computationally cheap. {To the best of our knowledge, they are first direct algorithms for the inverse conductivity problem in the context of the subdiffusion model}. Moreover, we provide relevant theoretical underpinnings for the algorithms. Also we present numerical results to illustrate their performance under various scenarios, e.g., the size of inclusions, noise level of the data, and the number of inclusions, showing that the algorithms are efficient and robust.
\end{abstract}



\pagestyle{myheadings}
\thispagestyle{plain}

\section{Introduction}
In this work, we aim at recovering small conductivity inclusions within a homogeneous background in subdiffusion processes. Let $\Omega\subset\mathbb{R}^d$ ($d=2,3$) be an open  bounded Lipschitz domain with a boundary $\partial\Omega$. Supose that there are $m$ small inclusions $A_{\ell} = \varepsilon B_{\ell} + z_{\ell}$, $\ell = 1, \ldots, m$, where $z_{\ell} \in \mathbb{R}^d$ and $B_{\ell}$ is a bounded Lipschitz domain in $\mathbb{R}^d$ containing the origin, and let $A = \bigcup_{\ell = 1}^m A_{\ell}$. The inclusions $\{ A_{\ell} \}_{\ell = 1}^m$ are assumed to be well separated from each other and from the boundary $\partial \Omega$, i.e., there exist positive constants $c$ and $\varepsilon_0$ such that for all $\varepsilon\in(0,\varepsilon_0)$, there holds
\begin{align}\label{eq:assume:separation}
	\inf_{x \in A_{\ell}} \text{dist}(x, \partial \Omega) > {\color{blue}c}, \quad \forall \ell \in \{ 1, \ldots, m \}, \quad \text{and} \quad \text{dist}(A_j, A_k) > {\color{blue}c}, \quad j \neq k.
\end{align}
The background is homogeneous with a conductivity $\gamma_0$ and the inclusion $A_{\ell}$ has a conductivity $\gamma_{\ell}$, where $\gamma_{\ell}$ and $\gamma_0$ are positive constants with $\gamma_{\ell} \neq \gamma_0$, $\ell = 1, \ldots, m$. The conductivity profile $\gamma_{\varepsilon}$ over the domain $\Omega$ is given by
\begin{equation}
\gamma_{\varepsilon}(x) =\left\{
		\begin{aligned}
			& \gamma_0, && x\in \Omega \backslash \overline{A}, \\
			& \gamma_{\ell}, && x\in A_{\ell}, \quad \ell = 1, \ldots, m.
	\end{aligned}\right.
\end{equation}

Consider the following initial boundary value problem:
\begin{equation}\label{eq:tfd:u}
		\left\{
		\begin{aligned}
			\partial_t^{\alpha} u(x,t) - \nabla \cdot (\gamma_{\varepsilon}(x) \nabla u(x,t)) &= f(x,t), && (x,t)\in \Omega \times (0, T), \\
			u(x,0) &= u_0(x), && x\in \Omega, \\
			\gamma_0 \partial_n u(x,t) &= g(x,t), && (x,t)\in \partial \Omega \times (0, T),
		\end{aligned}
		\right.
	\end{equation}
where $n$ is the unit outward normal vector to the boundary $\partial \Omega$, and the Neumann boundary data $g$ and the initial data $u_0$ satisfy the compatibility condition
$\gamma_0 \partial_nu_0(x) = g(x, 0) \text{ for } x \in \partial \Omega$. Here $\partial_t^{\alpha} u(x,t)$ denotes the left-sided Caputo fractional derivative of $u$ in time of order $\alpha\in(0,1)$, defined by \cite[p. 42]{Jin:2021book}
\begin{align*}
	\partial_t^{\alpha} u(x, t) =\frac{1}{\Gamma(1 - \alpha)}  \int_0^t (t-s)^{-\alpha}\partial_s u(x, s) \ds,
	\end{align*}
where $\Gamma(\cdot)$ denotes Euler's Gamma function.
The model \eqref{eq:tfd:u} is commonly known as the subdiffusion model, and can accurately describe anomalously slow diffusion processes.
It has found many successful applications, e.g., solute transport in heterogeneous media \cite{dentz2004time,berkowitz2006modeling} and protein transport within membranes \cite{kou2004generalized,ritchie2005detection}. The well-posedness of problem \eqref{eq:tfd:u} in Sobolev spaces can be analyzed along the line of \cite{KubicaYamamoto:2020} and \cite[Section 6.1]{Jin:2021book}.
	
The concerned inverse problem is to recover the locations of the  inclusions $\{ A_{\ell} \}_{\ell = 1}^m$ from boundary measurement. It is related to nondestructive testing, i.e., locating anomalies inside a known homogeneous medium. This is a version of the inverse conductivity problem for the model \eqref{eq:tfd:u}. The uniqueness of identifying a general $\gamma(x)$ in the 1D case from one boundary flux was proved in \cite{ChengYamamoto:2009}. The multi-dimensional case was studied in \cite{LiJiaYamamoto:2013} and \cite{CenJinLiuZhou:2023} for a general $\gamma(x)$ using the Dirichlet to Neumann map and  a piecewise constant $\gamma(x)$ from one boundary flux, respectively. Kian et al \cite{KianLiLiuYamamoto:2021} proved the unique determination of multiple coefficients from one boundary flux, excited by one specialized Dirichlet datum. Numerically, the reconstruction is commonly carried out using a regularized least-squares formulation. The numerical analysis of the schemes was given in \cite{JinZhou:2021sicon,JinLuQuanZhou:2024} for space-time or terminal data. Cen et al \cite{CenJinLiuZhou:2023} employed the level set method to recover a piecewise constant $\gamma(x)$. The regularization procedure is very flexible and works well for general conductivities, but solving the resulting optimization problems can be computationally demanding.

In this work, we develop novel direct algorithms for imaging small conductivity inclusions $\{ A_{\ell} \}_{\ell = 1}^m$ from the boundary measurement $I_\Phi$ defined in \eqref{eq:wbm:definition}. First we derive an asymptotic expansion of $I_\Phi$ with respect to the size $\varepsilon$ of the inclusions, cf. Theorem \ref{thm:asymptotic-exp}, which forms the basis for constructing the direct algorithms.
The derivation relies on the coercivity of the Caputo fractional derivative (i.e., Alikhanov's inequality \cite{alikhanov2010priori}). One natural choice of the test function $\Phi$ in $I_\Phi$ is the fundamental solution to the model \eqref{eq:tfd:u}, but its accurate evaluation is highly nontrivial. One crucial step towards developing efficient algorithms is to use approximate fundamental solutions (i.e., truncating the asymptotic expansion of the reduced Green function).
We analyze the perturbation of $I_\Phi$ due to the approximate fundamental solution, and give sufficient conditions on the truncation level and the source location to ensure a similar asymptotic expansion of the resulting boundary measurement $I_\Phi$.
The leading-order terms of the new asymptotic expansions lend themselves to the development of two computationally tractable algorithms, one for locating one inclusion, and the other for locating multiple inclusions, in both 2D and 3D cases.
Moreover, we provide theoretical underpinnings for the algorithms. The numerical results show that the algorithms can accurately predict the locations for both circular and elliptical inclusions.
To the best of our knowledge, they represent first direct algorithms for an inverse conductivity problem associated with the model \eqref{eq:tfd:u}.

The asymptotic expansion of the solutions with respect to the size of small inclusions has been extensively studied for standard elliptic and parabolic problems; see the monographs \cite{ammari2004reconstruction,AmmariKang:2007} for in-depth treatments. The resulting asymptotic expansion can serve as a powerful computational tool for recovering the locations and shapes of small inclusions in elliptic and parabolic problems. Ammari et al \cite{AmmariKangKim:2005} developed direct algorithms to recover the locations of small conductivity inclusions in the heat equation in the 2D case from boundary measurement, and analyzed the properties of the algorithms. Bouraoui et al \cite{BouraouiKhelifi:2016,BouraouiKhelifi:2017} extended the approach in \cite{AmmariKangKim:2005} to identify locations and shape properties of small polygonal conductivity defects in the heat equation from partial boundary measurement. In this work, we utilize asymptotic expansions to recover the locations of small inclusions, which extends the 2D algorithms of the work \cite{AmmariKangKim:2005} to the subdiffusion case, and also provide new algorithms in the 3D case. This development requires several tools from fractional calculus (e.g., integration by parts formula and Alikhanov's inequality) and delicate analysis of the approximate fundamental solution, and the analysis of the algorithms is highly involved. These also represent the main technical novelties of the study.

The rest of the paper is organized as follows. In Section \ref{sec:main}, we present an asymptotic analysis of the boundary measurement. In Section \ref{sec:alg}, we develop two direct algorithms for recovering small inclusions and provide relevant theoretical underpinnings. In Section \ref{sec:exp}, we present numerical experiments to complement the theoretical analysis. Finally, we give concluding remarks in Section \ref{sec: conclusion}. Throughout, $AC[0, T]$ denotes the space of absolutely continuous functions on $(0, T)$. For any $w \in AC[0, T]$, ${}_t\widetilde{\partial}_{T}^{\alpha}w$ denotes the right-sided Caputo derivative defined by
$${}_t\widetilde{\partial}_{T}^{\alpha} w(t) = - \frac{1}{\Gamma(1 - \alpha)} \int_t^{T} (s - t)^{-\alpha} \partial_s w(s) \ds.$$
The notation $C$ denotes a generic constant which may change at each occurrence but it is always independent of $\varepsilon$.

\section{Asymptotic expansion of boundary measurement}\label{sec:main}
Now we establish the asymptotic expansion of the boundary measurement with respect to the size of small inclusions. This result will inspire the construction of direct imaging algorithms in Section \ref{sec:alg}.

\subsection{Asymptotic expansion}

To reconstruct the inclusions, we define the boundary measurement
\begin{align}\label{eq:wbm:definition}
	I_{\Phi} = \int_0^T\!\!\! \int_{\partial \Omega} \gamma_0 (u - U)(x, t) \partial_n \Phi(x, t) \dd \sigma_x \dd t,
\end{align}
where $\dd \sigma_x$ denotes surface integral, $\Phi$ is a test function to be suitably chosen, and $U$ is the background solution satisfying
\begin{equation}\label{eq:bckgrdsol}
		\left\{\begin{aligned}
			\partial_t^{\alpha} U(x,t) - \gamma_0 \Delta U(x,t) &= f(x,t), && (x,t)\in \Omega \times (0, T), \\
			U(x,0) &= u_0(x), && x\in  \Omega, \\
			\gamma_0 \partial_nU(x,t) &= g(x,t), && (x,t)\in  \partial \Omega \times (0, T).
		\end{aligned}\right.
\end{equation}
Throughout, $\Phi$ and $U$ are taken to be smooth functions with continuous derivatives on the cylindrical domain $\overline{\Om}\times[0,T]$.
The existence, uniqueness and regularity of the solutions $u$ and $U$ to problems \eqref{eq:tfd:u} and \eqref{eq:bckgrdsol} can be found in \cite[Theorems 3.5, 3.8 and 3.11]{Mu:2017:ERSTFD}.
		
We also define $\xi_j^{(\ell)}$ for $j = 1, \ldots, d$ and $\ell = 1, \ldots, m$, which is the unique solution to
\begin{equation}
\left\{	\begin{aligned}
	\Delta \xi_j^{(\ell)}(x) &= 0, && x\in \mathbb{R}^d \backslash \partial B_{\ell}, \\
	\xi_j^{(\ell)}|_{+}(x) - \xi_j^{(\ell)}|_{-}(x) &= x_j, && x\in \partial B_{\ell}, \\
	\gamma_0 \partial_n \xi_j^{(\ell)}|_{+}(x) - \gamma_{\ell} \partial_n \xi_j^{(\ell)}|_{-}(x) &= \gamma_0 \partial_n x_j, && x\in   \partial B_{\ell}, \\
	\xi_j^{(\ell)}(x) &= \mathcal{O}(|x|^{1 - d}), && \text{ as } |x| \to \infty,
\end{aligned}\right.
\end{equation}
where $n$ is the outward unit normal vector to $\partial B_{\ell}$, and the subscripts $+$ and $-$ indicate the limit value as $x$ approaches $\partial B_{\ell}$ from outside $B_{\ell}$ and from inside $B_{\ell}$, respectively. Let
\begin{equation} \label{eq: full definition of local coordinates}
\Upsilon_j^{(\ell)}(x) =
\left\{\begin{aligned}
& \xi_j^{(\ell)}(x), && x \in \mathbb{R}^d \backslash \overline{B_{\ell}}, \\
& \xi_j^{(\ell)}(x) + x_j, && x \in B_{\ell}.
\end{aligned}\right.
\end{equation}
Then on the interface $\partial B_\ell$, the following two identities hold:
	\begin{align*}
		\Upsilon_j^{(\ell)}|_{+}(x) - \Upsilon_j^{(\ell)}|_{-}(x) = 0\quad\mbox{and}\quad \gamma_0\partial_n \Upsilon_j^{(\ell)}|_{+}(x) - \gamma_{\ell} \partial_n\Upsilon_j^{(\ell)}|_{-}(x) = (\gamma_0 - \gamma_{\ell})\partial_n x_j.
	\end{align*}
For $\ell = 1, \ldots, m$, let \begin{equation*} M^{(\ell)} = [m_{jk}^{(\ell)}] \in\mathbb{R}^{d\times d},\quad \mbox{with }
m_{jk}^{(\ell)} = \int_{\partial B_{\ell}} x_j \partial_n\xi_{k}^{(\ell)}|_{-}(x) \dd \sigma_x.
\end{equation*}
The matrix $M^{(\ell)}$ is known as the P\'{o}lya-Szeg\"{o} polarization tensor associated with the inclusion $B_{\ell}$ and the conductivity $\gamma_{\ell}$ \cite{polya1951isoperimetric}. The properties of $M^{(\ell)}$ were studied in \cite[Chapter 3]{ammari2004reconstruction}. If the inclusion $B_{\ell}$ is a disk, then with $I_d\in \mathbb{R}^{d \times d}$ being the identity matrix, $M^{(\ell)}\in \mathbb{R}^{d\times d}$ is given by
\begin{equation}\label{eqn:polar-circ}
    M^{(\ell)} = - \frac{d \gamma_0 |B_{\ell}|}{\gamma_{\ell} + (d - 1) \gamma_0} I_d.
\end{equation}

We have the following asymptotic expansion of $I_{\Phi}$ as $\varepsilon \to 0$ under a mild condition on $\Phi$. Below we use the shorthand notation $\ell_{\varepsilon,d} = |\log\varepsilon|$ for $d=2$, and $\ell_{\varepsilon,d} = 1$ for $d=3$.
	\begin{theorem} \label{thm:asymptotic-exp}
		Let the test function $\Phi$ satisfy
  \begin{equation} \label{eqn:PDE-Phi}
	\left\{\begin{aligned}
		{}_t\widetilde{\partial}_{T}^{\alpha} \Phi(x,t) - \gamma_0 \Delta \Phi(x,t) &= 0, && (x,t)\in \Omega \times (0, T), \\
		\Phi(x, T) &= 0, && x\in \Omega.
		\end{aligned}\right.
\end{equation}
  Then, for each fixed $T>0$, with $R_d(\varepsilon) = \mathcal{O}(\ell_{\varepsilon,d}^{1/2}\varepsilon^{d+1})$, there holds as $\varepsilon \to 0$
		\begin{align} \label{eq:asymptotic-exp}
			I_{\Phi} = - \varepsilon^d \sum_{\ell = 1}^m (\gamma_0 - \gamma_{\ell}) \int_0^T \nabla U(z_{\ell}, t) \cdot M^{(\ell)} \nabla \Phi(z_{\ell}, t) \dd t + R_d(\varepsilon).
		\end{align}
\end{theorem}

In view of Theorem \ref{thm:asymptotic-exp}, for small $\varepsilon$, we can approximate $I_{\Phi}$ as
$$ I_\Phi \approx - \varepsilon^d \sum_{\ell = 1}^m (\gamma_0 - \gamma_{\ell}) \int_0^T \nabla U(z_{\ell}, t) \cdot M^{(\ell)} \nabla \Phi(z_{\ell}, t) \dd t.$$
Note that problem \eqref{eqn:PDE-Phi} is not prescribed with boundary conditions, and thus the choice of $\Phi$ has some flexibility. By suitably choosing $\Phi$ and $U$, it motivates the development of direct algorithms for recovering the locations of small inclusions.
	
\subsection{Proof of Theorem \ref{thm:asymptotic-exp} }\label{ssec:anal}
	
First we recall several preliminary results on fractional calculus. See the monographs \cite{Podlubny:1999,kilbas2006theory,Jin:2021book} for further details. We often use the Riemann-Liouville fractional integral of order $\alpha >0$, defined by \cite[p. 22]{Jin:2021book}
$${_0I_t^{\alpha}} w(t) = \frac{1}{\Gamma(\alpha)} \int_0^t (t - s)^{\alpha - 1} w(s) \ds.$$
Then the following version of the fundamental theorem of calculus holds \cite[Theorem 2.13]{Jin:2021book}:
\begin{align} \label{eq: Newton-Lebiniz}
		w(t) = w(0) + {_0I_t^\alpha }(\partial_t^\alpha w) = w(0) + \frac{1}{\Gamma(\alpha)} \int_0^t (t - s)^{\alpha - 1} {\partial_s^{\alpha}} w(s) \ds,\quad \forall w \in AC[0,T].
	\end{align}

The next result gives an \textit{a priori} estimate. We define the norm $\tribar\cdot\tribar$ by
\begin{equation*}
\tribar v \tribar^2 :=
\sup_{0 \le t \le T}{_0I_t^{\alpha}}\big( \| v(\cdot, t) \|_{L^2(\Omega)}^2 +  \| \nabla v(\cdot, t) \|_{L^2(\Omega)}^2 \big) = \sup_{0 \le t \le T}{_0I_t^{\alpha}} \| v(\cdot, t) \|_{H^1(\Omega)}^2.
\end{equation*}
\begin{lemma} \label{lem: regularity estimate of general problem}
For $V_0 \in L^2(\Omega)$, $F\in L^2(0,T;L^2(\Omega))$, $G\in L^2(0,T;L^2(\partial\Omega))$ and $H\in L^2(0,T;H^{2}(A))$, let $v$ solve
		\begin{equation} \label{eq: general problem}
			\left\{\begin{aligned}
				\partial_t^{\alpha} v(x,t)- \gamma_0 \Delta v(x,t) &= F(x,t), && (x,t)\in \Omega \backslash \overline{A} \times (0, T), \\
				\partial_t^{\alpha} v(x,t)- \gamma_{\ell} \Delta v(x,t) &= F(x,t), && (x,t)\in A_{\ell} \times (0, T), \quad \ell = 1, \ldots, m, \\
				v|_{+}(x,t) - v|_{-} (x,t)&= 0, && (x,t)\in  \partial A_{\ell} \times (0, T), \quad \ell = 1, \ldots, m, \\
				\gamma_0 \partial_n v|_{+}(x,t) - \gamma_{\ell} \partial_n v|_{-} (x,t)&= \partial_n H(x,t), && (x,t)\in \partial A_{\ell} \times (0, T), \quad \ell = 1, \ldots, m, \\
				v(x,0) &= V_0(x), && x\in  \Omega, \\
				\gamma_0 \partial_n v(x,t) &= G(x,t), && (x,t)\in \partial \Omega \times (0, T).
			\end{aligned}\right.
		\end{equation}
Then there exists  $C>0$ depending on $\alpha$, $T$, $\{ \gamma_{\ell} \}_{\ell = 0}^m$ and the Lipschitz character of $\Omega$ and $\{ A_{\ell} \}_{\ell = 1}^m$ such that the following estimate holds	\begin{align*}
\tribar v \tribar^2
\le & C\| V_0 \|_{L^2(\Omega)}^2 + C \sup_{0 \le t \le T}{_0I_t^{\alpha}}\Big(  \| F (\cdot,t)\|_{L^2(\Omega\backslash\overline{A})}^2 + \| F(\cdot,t) -\Delta H(\cdot,t)\|_{L^2(A)}^2 \\
 &+ \| \nabla H(\cdot,t) \|_{L^2(A)}^2 +  \| G (\cdot,t)\|_{L^2(\partial \Omega)}^2  \Big).
\end{align*}
\end{lemma}
\begin{proof}
By integration by parts, divergence theorem and the governing equation \eqref{eq: general problem} for $v(x,t)$, we obtain
		\begin{align*}
			&\int_{\Omega} \gamma_{\varepsilon}(x) |\nabla v(x,t)|^2{\rm d}x  = \int_{\Omega \backslash A} \gamma_0 | \nabla v (x,t)|^2{\rm d}x + \sum_{\ell = 1}^m \int_{A_{\ell}} \gamma_{\ell} |\nabla v(x,t)|^2{\rm \d}x\\
			 =& \int_{\partial \Omega} \gamma_0 v(x,t) \partial_n v(x,t){\rm d}\sigma_x+ \sum_{\ell = 1}^m \int_{\partial A_{\ell}} \left( \gamma_{\ell} \partial_n v|_{-}(x,t) - \gamma_0\partial_n v|_{+}(x,t) \right){\rm d}\sigma_x \\
             &- \int_{\Omega} (\partial_t^{\alpha} v(x,t) - F(x,t))v(x,t){\rm d}x\\
			 = &\int_{\partial \Omega} G(x,t) v(x,t) {\rm d}\sigma_x - \int_{\partial A}\partial_n H(x,t)v(x,t){\rm d}\sigma_x - \int_{\Omega} (\partial_t^{\alpha} v(x,t) - F(x,t))v(x,t)\dx.
		\end{align*}
		By Green's identity on the subdomain $A$, we obtain
		\begin{align*}
			\int_{\Omega} \gamma_{\varepsilon}(x) |\nabla v(x,t)|^2{\rm d}x = & \int_{\partial \Omega} G(x,t) v(x,t) {\rm d}\sigma_x - \int_{\Omega} (\partial_t^{\alpha} v(x,t))v(x,t){\rm d}x \\
            &+  \int_{\Om\backslash\overline{A}}F (x,t)v(x,t){\rm d}x
             + \int_A(F(x,t)-\Delta H(x,t))v(x,t){\rm d}x\\
             &- \int_{A} \nabla H(x,t) \cdot\nabla v(x,t){\rm d}x.
		\end{align*}
By Alikhanov's inequality \cite[Lemma 1]{alikhanov2010priori}
$$w(t) \partial_t^\alpha w(t) \ge \tfrac{1}{2} \partial_t^{\alpha} w^2(t),\quad \forall w\in AC[0,T],$$
we have for all $t\in[0,T]$,
\begin{align*}
		 {\rm I}:=&{_0I_t^{\alpha}} \int_{\Omega} \frac{1}{2} \partial_t^{\alpha} v(x,t)^2\dx + {_0I_t^{\alpha}} \int_{\Omega}\gamma_{\varepsilon}(x) |\nabla v(x,t)|^2\dx\\
			\le & {_0I_t^{\alpha}} \Big(\int_{\partial \Omega} G(x,t) v(x,t){\rm d}\sigma_x + \int_{\Om\backslash\overline{A}}F(x,t) v(x,t)\dx \\
            &+ \int_A(F(x,t)-\Delta H(x,t))v(x,t)\dx - \int_{A} \nabla H(x,t) \cdot\nabla v(x,t)\dx\Big) := \sum_{j = 1}^4 {\rm I}_j.
   \end{align*}
Then the Cauchy-Schwarz inequality, and the trace inequality $\| w \|_{L^2(\partial \Omega)} \le C(\Omega) \| w \|_{H^1(\Omega)}$ (for $w\in H^1(\Omega)$) \cite{evans2018measure} imply
$${\rm I}_1\le C({_0I_t^{\alpha}}\| G(x,t) \|_{L^2(\partial \Omega)}^2)^{1/2} ({_0I_t^{\alpha}}\| v(x,t) \|_{H^1(\Omega)}^2)^{1/2},
$$
and similar bounds on the terms ${\rm I}_j$, $j=2,3,4$. Consequently,
\begin{align}\label{eq:cauchyschwarz:lem}
		{\rm  I}	\le &  \Big( C({_0I_t^{\alpha}}\| G (x,t)\|_{L^2(\partial \Omega)}^2)^{1/2} + ({_0I_t^{\alpha}}\| F (x,t)\|_{L^2(\Omega\backslash\overline{A})}^2)^{1/2} \nonumber\\
        & + ({_0I_t^{\alpha}}\| F(x,t)-\Delta_x H(x,t) \|_{L^2(A)}^2)^{1/2}+ ({_0I_t^{\alpha}}\| \nabla H(x,t) \|_{L^2(A)}^2)^{1/2}\Big) \tribar v\tribar,
		\end{align}
for some $C=C(\alpha,T,\Omega,\{A_{\ell} \}_{\ell = 1}^m)$.
Clearly, we have for any $w\in L^\infty(0,T)$,
\begin{align} \label{eq: monotonicity of integral}
\sup_{t\in[0,T]}{_0I_t^{\alpha}} w(t)
 &\le \|w\|_{L^\infty(0,T)} \sup_{t\in[0,T]}\frac{1}{\Gamma(\alpha)}\int_0^t (t-s)^{\alpha-1}\,{\rm d}s=\frac{T^\alpha}{\Gamma(1+\alpha)}\|w\|_{L^\infty(0,T)}.
\end{align}
Then, the choice $\mu = \min \{ \Gamma(1+\alpha)/2, \min_{0 \le \ell \le m} \gamma_{\ell} \}$, the identity \eqref{eq: monotonicity of integral} with $w(t) = \| v(x, t) \|_{L^2(\Omega)}$ and \eqref{eq: Newton-Lebiniz} give
		\begin{align*}
			\mu\tribar v\tribar^2  \le &\frac{T^\alpha}{2}\sup_{0 \le t \le T}\| v(x,t) \|_{L^2(\Omega)}^2 + \sup_{0 \le t \le T} {_0I_t^{\alpha}}  \int_{\Omega}\gamma_{\varepsilon}(x) |\nabla v(x,t)|^2\dx\\
			=&\frac{T^\alpha}{2}\| v(x,0) \|_{L^2(\Omega)}^2 +  T^\alpha \sup_{0 \le t \le T} {_0I_t^{\alpha}} \int_{\Omega} \tfrac{1}{2} \partial_t^{\alpha} v(x,t)^2\dx\\
             &+  \sup_{0 \le t \le T}{_0I_t^{\alpha}}\int_{\Omega}\gamma_{\varepsilon} |\nabla v(x,t)|^2 \dx.
		\end{align*}
		Using the estimate \eqref{eq:cauchyschwarz:lem} and Young's inequality, we arrive at
\begin{align*}
   \mu \tribar v \tribar^2  \le & \frac{T^\alpha}{2}\| V_0 \|_{L^2(\Omega)}^2 + 2 \kappa \tribar v \tribar^2  +  \frac{C \max(1, T^{\alpha})}{2 \kappa} \sup_{0 \le t \le T}{_0I_t^{\alpha}}\Big( \| F(x, t) \|_{L^2(\Omega\backslash\overline{A})}^2\\
   &+ \| F(x,t)-\Delta H(x,t) \|_{L^2(A)}^2  +  \| \nabla H(x,t) \|_{L^2(A)}^2 + \| G(x,t) \|_{L^2(\partial \Omega)}^2\Big) ,
		\end{align*}
which completes the proof of the lemma by choosing $\kappa$ sufficiently small.
\end{proof}	
	
The next result gives an integral representation of the boundary measurement $I_\Phi$.
\begin{lemma} \label{lem:equiv-wbm}Under the assumption of Theorem \ref{thm:asymptotic-exp}, we have
\begin{align}
I_{\Phi} = \sum_{\ell = 1}^m (\gamma_0 - \gamma_{\ell}) \int_0^T \!\!\!\int_{A_{\ell}} \nabla u (x,t)\cdot \nabla \Phi(x,t)\dd x \dd t.
\end{align}
\end{lemma}
\begin{proof}
The representation follows from the divergence theorem and the following integration by parts formula  \cite[Lemma 2.6]{Jin:2021book}
		\begin{align*}
			\int_0^T (\partial_t^{\alpha} w(t)) (v(t) - v(T^-)) \dd t = \int_0^T (w(t) - w(0^+)) ({}_t\widetilde{\partial}_{T}^{\alpha} v(t)) \dd t,\quad \forall v,w\in AC[0,T],
		\end{align*}
and the rest of the derivation is similar to that of \cite[Lemma 3.3]{AmmariKangKim:2005}. Thus the detail is omitted.
\end{proof}
	
The next result gives a crucial asymptotic expansion.
\begin{lemma} \label{lem:est-trunc}
Let $\Upsilon_j^{(\ell)}$ be defined in \eqref{eq: full definition of local coordinates}, and
		\begin{align} \label{eq: auxiliary asymptotic expansion}
			V(x, t) = u(x, t) - U(x, t) + \varepsilon \sum_{\ell = 1}^m \sum_{j = 1}^d \partial_{x_j} U(z_{\ell}, t) \Upsilon_{j}^{(\ell)}\left( \frac{x - z_{\ell}}{\varepsilon} \right).
		\end{align}
		Then there exists a constant $C$ independent of $\varepsilon$ such that
\begin{align} \label{eq:est-trunc}
{_0I_T^{\alpha}} \| V(\cdot, t) \|_{H^1(\Omega)}^2 \le C\ell_{\varepsilon,d}\varepsilon^{d+2}.
\end{align}
\end{lemma}
\begin{proof}
The proof is based on mathematical induction on the number $m$ of inclusions. Actually we shall prove a stronger assertion \eqref{eq:est-trunc:general}. In the induction step, Lemma \ref{lem: regularity estimate of general problem} plays a crucial role. Specifically, for any fixed $m$, $\gamma_0$, $\{ \gamma_{\ell} \}_{\ell = 1}^m$ and $\{ A_{\ell} \}_{\ell = 1}^m$, let $u^0 = U$, $u^{k}$ be the solution to problem \eqref{eq:tfd:u} with the presence of the first $k$ inclusions $\{ A_{\ell} \}_{\ell = 1}^k$, for $k=1,\ldots, m$, and  let
\begin{equation*}
V^k(x, t) = u^k(x, t) - U(x, t) + \varepsilon \sum_{\ell = 1}^k \nabla U(z_{\ell}, t) \cdot \Upsilon^{(\ell)}(\varepsilon^{-1}(x - z_{\ell})),
\end{equation*}
with $\Upsilon^{(\ell)} = (\Upsilon^{(\ell)}_1, \ldots, \Upsilon^{(\ell)}_d)$. Also let $\tV_k = V^k - V^{k - 1}$ for $k = 1, \ldots, m$, with $V^0 = 0$.
For each fixed $k$, with $\tA_k = \bigcup_{\ell = 1}^k A_{\ell}$, the function $\tV_k$ satisfies
\begin{equation*}
\left\{\begin{aligned}
\partial_t^{\alpha} \tV_k(x,t)- \gamma_0 \Delta \tV_k(x,t) &= F_k(x,t), && (x,t)\in \Omega \backslash \overline{\tA_k} \times (0, T), \\
\partial_t^{\alpha} \tV_k(x,t)- \gamma_{\ell} \Delta \tV_k(x,t) &= F_k(x,t), && (x,t)\in A_{\ell} \times (0, T), \quad \ell = 1, \ldots, k, \\
\tV_k|_{+}(x,t) - \tV_k|_{-} (x,t)&= 0, && (x,t)\in \partial A_{\ell} \times (0, T), \quad \ell = 1, \ldots, k, \\
\gamma_0 \partial_n \tV_k|_{+}(x,t) - \gamma_{\ell} \partial_n \tV_k|_{-} (x,t)&= \partial_n H_k(x,t), && (x,t)\in \partial A_{\ell} \times (0, T), \quad \ell = 1, \ldots, k, \\
\tV_k(x,0) &= V_{0, k}(x), && x\in \Omega, \\
\gamma_0 \partial_n \tV_k(x,t) &= G_k(x,t), && (x,t)\in \partial \Omega \times (0, T),
\end{aligned}\right.
\end{equation*}
where the functions $F_k(x,t):\Omega\times (0,T)\to \mathbb{R}$, $H_k(x,t):\widetilde{A}_k\times(0,T)\to \mathbb{R}$, $V_{0, k}(x):\Omega \to\mathbb{R}$ and $G_k(x,t):\partial\Omega\to \mathbb{R}$ are given by
\begin{align*}
F_k(x,t) &= (\gamma_k - \gamma_0) \Delta u^{k - 1}(x,t) \chi_{A_k}(x) + \varepsilon \partial_t^{\alpha} \nabla U(z_k, t) \cdot \Upsilon^{(k)} (\varepsilon^{-1}(x - z_{k})),\\
H_k(x,t) & =
\begin{cases}
\varepsilon (\gamma_0-\gamma_\ell) \nabla U(z_{k}, t) \cdot \Upsilon^{(k)}(\varepsilon^{-1}(x - z_{k})),&\,\, (x,t)\in A_\ell \times (0, T), 1\le\ell< k,\\
(\gamma_k - \gamma_0) (u^{k - 1}(x, t) - \nabla U(z_{k}, t) \cdot x),&\,\, (x,r)\in A_k \times (0, T),
\end{cases} \\
V_{0, k}(x) & = \varepsilon \nabla u_0(z_k) \cdot \Upsilon^{(k)} (\varepsilon^{-1}(x - z_{k})), \\
G(x,t) & = \gamma_0 \nabla U(z_k, t) \cdot \partial_n \Upsilon^{(k)}( \varepsilon^{-1} (x - z_{k})).
\end{align*}
It follows directly from Lemma \ref{lem: regularity estimate of general problem} that
\begin{align} \label{eq: regularity estimate of difference}
\sup_{0\le t\le T} {_0I_t^{\alpha}} &\| \tV_k(\cdot,t)\|_{H^1(\Omega)}^2 \le  C \| V_{0, k} \|_{L^2(\Omega)}^2 + C \sup_{0 \le t \le T}{_0I_t^{\alpha}} \Big( \| F_k(\cdot,t) \|_{L^2(\Omega\backslash\overline{\tA_k})}^2 \\
& + \| F_k(\cdot,t) -\Delta H_k(\cdot,t)\|_{L^2(\tA_k)}^2 +  \| G_k(\cdot,t)\|_{L^2(\partial \Omega)}^2 +  \| \nabla H_k(\cdot,t) \|_{L^2(\tA_k)}^2 \Big).\nonumber
\end{align}
Next we bound the first four terms of the inequality \eqref{eq: regularity estimate of difference} for each fixed $k = 1, \ldots, m$. Let $\widetilde{\Omega}_k = \varepsilon^{-1} (\Omega - z_k)$. Note that $\text{diam}(\widetilde{\Omega}_k) = \mathcal{O}(\varepsilon^{-1})$ due to the assumption \eqref{eq:assume:separation}. Moreover, since $\Upsilon^{(k)}_j$ is bounded and $\Upsilon^{(k)}_j(x) = \mathcal{O}(|x|^{1 - d})$ as $|x| \to \infty$ for each $j$ and $k$, there exist $C,C'>0$ such that
\begin{equation*}
\|\Upsilon^{(k)}_j\|_{L^2(\widetilde{\Om}_k)}^2 \le \begin{cases}
C'+\int_{1<|x|<C\varepsilon^{-1}}\, (C|x|^{-1})^2{\rm d}x=C'+2\pi C^2\log(C\varepsilon^{-1}),&\mbox{if }d=2,\\
C'+\int_{|x|>1}(C|x|^{-2})^{2}\,{\rm d}x<\infty,&\mbox{if }d=3,
\end{cases}
\end{equation*}
i.e., $\|\Upsilon^{(k)}_j\|_{L^2(\widetilde{\Omega}_k)}^2 =\mathcal{O}(\ell_{\varepsilon,d}).$
Then, changing variables $\widetilde{x} = \varepsilon^{-1} (x - z_k)$ leads to
\begin{align} \label{eq: estimate of initial value}
&\| V_{0, k}(x) \|_{L^2(\Omega)}^2 = \varepsilon^d \|V_{0, k}(\varepsilon\widetilde{x} + z_k)\|_{L^2(\widetilde{\Om}_k)}^2 \\
=& \varepsilon^d \|\varepsilon \nabla u_0(z_k) \cdot \Upsilon^{(k)}(\widetilde{x})\|_{L^2(\widetilde{\Om}_k)}^2 \le C\ell_{\varepsilon,d}\varepsilon^{d+2}.\nonumber
\end{align}
Let $\tB_k = \varepsilon^{-1}(\tA_k - z_k)$. Then similarly, due to the assumption \eqref{eq:assume:separation}, we have $$\|\Upsilon^{(k)}_j\|_{L^2(\widetilde{\Omega}_k \backslash \tB_k)}^2 =\mathcal{O}(\ell_{\varepsilon,d}).$$
Consequently,
\begin{align} \label{eq: estimate of source}
&{_0I_t^{\alpha}}\| F_k(x, t) \|_{L^2(\Omega \backslash \overline{\tA_k})}^2  = {_0I_t^{\alpha}}\big(\varepsilon^d \|F_k(\varepsilon\widetilde{x} + z_k, t)\|_{L^2(\widetilde{\Om}_k \backslash \overline{\tB_k})}^2\big) \\
=& {_0I_t^{\alpha}}\big(\varepsilon^d \|\varepsilon \partial_t^{\alpha} \nabla U(z_k,t) \cdot \Upsilon^{(k)}(\widetilde{x}) \|_{L^2(\widetilde{\Om}_k \backslash \overline{\tB_k})}^2 \big)\le C\ell_{\varepsilon,d}\varepsilon^{d+2}. \nonumber
\end{align}
Next, note that the following identity holds:
$$F_k(x,t) -\Delta H_k(x,t) = \varepsilon \partial_t^{\alpha} \nabla U(z_k, t) \cdot \Upsilon^{(k)} (\varepsilon^{-1}(x - z_{k})),\quad  x \in \tA_k.$$
Since $\varepsilon^{-1} | x - z_k | = \mathcal{O}(\varepsilon^{-1})$ for $x \in A_{\ell}$ with $\ell \neq k$ (due to the assumption \eqref{eq:assume:separation}), $\Upsilon^{(k)}_j$ is bounded and $\Upsilon^{(k)}_j = \mathcal{O}(|x|^{1 - d})$ as $|x| \rightarrow \infty$. Thus we have
\begin{equation*}
\Upsilon^{(k)}_j(\varepsilon^{-1}(x - z_{k})) =\left\{\begin{aligned}
\mathcal{O}(\varepsilon^{d - 1}), &\quad x \in A_{\ell}, \ell \neq k,\\ \mathcal{O}(1), &\quad  x \in A_k.
\end{aligned}\right.
\end{equation*}
Consequently, we deduce
\begin{align} \label{eq: estimate of source and interface}
& {_0I_t^{\alpha}}\| F_k(x, t) - \Delta H_k(x,t) \|_{L^2({\tA_k})}^2 \\
=& \sum_{\ell = 1}^k {_0I_t^{\alpha}}\big(\| \varepsilon \p_t^\alpha\nabla U(z_k,t)\cdot \Upsilon^{(k)}(\varepsilon^{-1}(x - z_k))\|_{L^2({A_{\ell}})}^2 \big) \le  C \varepsilon^{d+2}. \nonumber
\end{align}
Furthermore, for each $j$, since $| \nabla \Upsilon^{(k)}_j(x) | = \mathcal{O}(|x|^{-d})$ as $|x| \to \infty$, we have $| \nabla \Upsilon^{(k)}_j(x) | = \mathcal{O}(\varepsilon^d)$ for $x \in \partial \widetilde{\Omega}_k$ due to \eqref{eq:assume:separation}. Hence,
\begin{align} \label{eq: estimate of boundary value}
{_0I_t^{\alpha}}\| G_k(x, t) \|_{L^2(\p\Omega)}^2 & = {_0I_t^{\alpha}}\big(\varepsilon^{d-1} \|G_k(\varepsilon\widetilde{x} + z_k, t) \|_{L^2(\p\widetilde{\Om}_k)}^2\big) \\
& = {_0I_t^{\alpha}}\big(\varepsilon^{d-1} \left\|\gamma_0\nabla U(z_k, t) \cdot \partial_n \Upsilon^{(k)}(\widetilde{x}) \right\|_{L^2(\p\widetilde{\Om}_k)}^2 \big) \nonumber \\
& \le C \varepsilon^{d-1}(\varepsilon^{d})^2|\p\widetilde{\Om}_k|=\mathcal{O}( \varepsilon^{d-1}(\varepsilon^{d})^2(\varepsilon^{-1})^{d-1})=\mathcal{O}(\varepsilon^{2d}). \nonumber
\end{align}
Note that the desired assertion \eqref{eq:est-trunc} is one special case (with $t=T$ and $k=m$) of the following relation:
\begin{equation}\label{eq:est-trunc:general}
    \sup_{0 \le t \le T}{_0I_t^{\alpha}} \| V^{k}(\cdot, t) \|_{H^1(\Omega)}^2\le C\ell_{\varepsilon,d}\varepsilon^{d+2},\quad k=1,\ldots,m.
\end{equation}
Below, we prove the assertion \eqref{eq:est-trunc:general} by mathematical induction on $k$. It suffices to estimate the last term ${_0I_{t}^{\alpha}} \| \nabla H_k(\cdot,t) \|_{L^2(\tA_k)}^2$ in the fundamental estimate \eqref{eq: regularity estimate of difference}. For $k = 1$, we have $V^1 = \tV_1$. Note that
\begin{align*}
\nabla H_1(x,t) = (\gamma_1 - \gamma_0) \left( \nabla U(x,t) - \nabla U(z_1,t) \right) = \mathcal{O}(|x|) = \mathcal{O}(\varepsilon),\quad (x, t) \in  A_1 \times (0, T).
\end{align*}
Then we deduce
\begin{align} \label{eq: estimate of interface for 1}
{_0I_{t}^{\alpha}} \| \nabla H_1(\cdot, t) \|_{L^2(A_1)}^2 \le C\varepsilon^2{_0I_{t}^{\alpha}} \| 1 \|_{L^2(A_1)}^2\le C\varepsilon^{d+2}.
\end{align}
Thus, the assertion \eqref{eq:est-trunc:general} is proved for the case $k = 1$ by combining the estimates \eqref{eq: regularity estimate of difference}--\eqref{eq: estimate of interface for 1}. Now fix $2\le k \le m$ and assume that the assertion \eqref{eq:est-trunc:general} holds for the case $k-1$ (instead of $k$). Then we have
\begin{align} \label{eq: recurrence of V}
{_0I_t^{\alpha}} \| V^{k}(\cdot,t) \|_{H^1(\Omega)}^2 \le 2 \big({_0I_t^{\alpha}} \| \tV_{k}(\cdot,t) \|_{H^1(\Omega)}^2 + {_0I_t^{\alpha}} \| V^{k - 1}(\cdot,t) \|_{H^1(\Omega)}^2\big).
\end{align}
Also, we have the identity
\begin{equation*}
\nabla H_k(x, t) = \left\{\begin{aligned}
\varepsilon (\gamma_0 - \gamma_{\ell}) \sum_{j = 1}^d \partial_{x_j} U(z_k, t) \nabla \Upsilon^{(k)}_j(\varepsilon^{-1}(x - z_k)), & \,\, (x, t) \in  A_{\ell} \times (0, T), 1 \le \ell < k.\\
(\gamma_k - \gamma_0) ( \nabla u^{k - 1}(x, t) - \nabla U(z_k, t) ), &\,\, (x, t) \in  A_k \times (0, T).
\end{aligned}\right.
\end{equation*}
Since $| \nabla \Upsilon^{(k)}_j(x) | = \mathcal{O}(|x|^{-d})$ as $|x| \to \infty$, we have $| \nabla \Upsilon^{(k)}_j(\varepsilon^{-1}(x - z_k)) | = \mathcal{O}(\varepsilon^d)$ for $x \in A_{\ell}$ with $\ell \neq k$ due to the assumption \eqref{eq:assume:separation}. Consequently,
\begin{align} \label{eq: estimate of interface for ell}
& {_0I_{t}^{\alpha}} \| \nabla H_k(x, t) \|_{L^2(\tA_{k}\backslash A_{k})}^2 \\
= &\sum_{\ell = 1}^{k - 1} {_0I_{t}^{\alpha}} \big( \| \varepsilon (\gamma_0 - \gamma_{\ell}) \sum_{j = 1}^d \partial_{x_j} U(z_k, t)  \nabla \Upsilon^{(k)}_j(\varepsilon^{-1}(x - z_k)) \|_{L^2(A_{\ell})}^2 \big) \le C \varepsilon^{3d + 2}. \nonumber
\end{align}
By repeating the argument for \eqref{eq: estimate of interface for 1}, we obtain
\begin{align*}
{_0I_{t}^{\alpha}} \| (\gamma_k - \gamma_0) ( \nabla U(\cdot, t) - \nabla U(z_k, t) ) \|_{L^2(A_k)}^2 \le C \varepsilon^{d + 2}.
\end{align*}
Also, the argument for the estimate \eqref{eq: estimate of interface for ell} leads to
\begin{align*}
{_0I_{t}^{\alpha}} \| \varepsilon \sum_{\ell = 1}^{k - 1} \sum_{j = 1}^d \partial_{x_j} U(z_{\ell}, t) \nabla \Upsilon^{(\ell)}_j(\varepsilon^{-1}(\cdot - z_{\ell})) \|_{L^2(A_k)}^2 \le C \varepsilon^{3d + 2}.
\end{align*}
In addition, the induction hypothesis for the case $k - 1$ gives
\begin{align*}
& {_0I_{t}^{\alpha}} \bigg\| \nabla u^{k - 1}(\cdot, t) - \nabla U(\cdot, t) + \varepsilon \sum_{\ell = 1}^{k - 1} \sum_{j = 1}^d \partial_{x_j} U(z_{\ell}, t) \nabla \Upsilon^{(\ell)}_j(\varepsilon^{-1}(\cdot - z_{\ell})) \bigg\|_{L^2(A_k)}^2 \\
\le& {_0I_{t}^{\alpha}} \| V^{k - 1} \|_{H^1(\Omega)}^2  \le C \ell_{\varepsilon, d} \varepsilon^{d + 2}.
\end{align*}
Combining the preceding three bounds yields
\begin{align} \label{eq: estimate of interface for m}
{_0I_{t}^{\alpha}} \| \nabla H_k(\cdot, t) \|_{L^2(A_{k})}^2\le C \ell_{\varepsilon, d} \varepsilon^{d + 2}.
\end{align}
Finally, the assertion \eqref{eq:est-trunc:general} is proved for all $2\le k \le m$ by the induction hypothesis for $k - 1$ with \eqref{eq: regularity estimate of difference}--\eqref{eq: estimate of boundary value} and \eqref{eq: recurrence of V}--\eqref{eq: estimate of interface for m}. This completes the proof of the lemma.
\end{proof}

Now we can give the proof of Theorem \ref{thm:asymptotic-exp}.
	\begin{proof}[Proof of Theorem \ref{thm:asymptotic-exp}] Let
		\begin{align*}
			\widehat{I}_{\Phi} &= \sum_{\ell = 1}^m (\gamma_0 - \gamma_{\ell}) \int_0^T \!\!\!\!\int_{A_{\ell}} \Big( \nabla U(x, t) - \sum_{j = 1}^d \partial_{x_j} U(z_s, t) \nabla \Upsilon_j^{(\ell)} \left( \varepsilon^{-1}(x - z_{\ell}) \right) \Big) \cdot \nabla \Phi(x, t) \dd x \dd t,\\
			\widetilde{I}_{\Phi} &= \sum_{\ell = 1}^m (\gamma_0\! -\! \gamma_{\ell}) \int_0^T \sum_{j = 1}^d \sum_{k = 1}^d \partial_{x_j} U(z_{\ell}, t) \partial_{x_k} \Phi(z_{\ell}, t) \dd t \int_{A_{\ell}} \nabla ( x_j \! - \! \Upsilon_j^{(\ell)}(\varepsilon^{-1}(x\! -\! z_{\ell})) ) \cdot \nabla x_k \dd x.
		\end{align*}
		Also let $V$ be defined in \eqref{eq: auxiliary asymptotic expansion}. By the representation in Lemma \ref{lem:equiv-wbm}, we have
		\begin{align*}
			&| I_{\Phi} - \widehat{I}_{\Phi}|
			=  \left| \sum_{\ell = 1}^m (\gamma_0 - \gamma_{\ell}) \int_0^T \!\!\!\!\int_{A_{\ell}} \nabla V(x, t) \cdot \nabla \Phi(x, t)  \dd x \dd t \right| \\
			\le & C \sum_{\ell = 1}^m \left( \int_0^T \!\!\!\!\Gamma(\alpha) (T - t)^{1 - \alpha} \|\nabla \Phi(\cdot, t)\|_{L^2(A_\ell)}^2 {\rm d}t \right)^{\frac{1}{2}}  \left( {_ 0I_T^\alpha}\|\nabla V(\cdot, t)\|_{L^2(\Omega)}^2\right)^{\frac{1}{2}} \\
			\le & C \varepsilon^{d / 2} ({_0I_T^{\alpha}} \| \nabla V(\cdot, t) \|_{L^2(\Omega)}^2)^{1/2}.
		\end{align*}
Since $\nabla U(x, t) = \sum_{j = 1}^d \partial_{x_j} U(z_{\ell}, t) \nabla x_j + \mathcal{O}(\varepsilon)$ and $\nabla \Phi(x, t) = \sum_{j = 1}^d \partial_{x_j} \Phi(z_{\ell}, t) \nabla x_j + \mathcal{O}(\varepsilon)$ for $x \in A_{\ell}$, we have
		$| \widehat{I}_{\Phi} - \widetilde{I}_{\Phi} | \le C \varepsilon^{d + 1}$.
		Furthermore, by the definition of the matrix $M^{(\ell)}$,
		\begin{align*}
			\widetilde{I}_{\Phi} & = \varepsilon^d \sum_{\ell = 1}^m (\gamma_0 - \gamma_{\ell}) \int_0^T \sum_{j = 1}^d \sum_{k = 1}^d \partial_{x_j} U(z_{\ell}, t) \partial_{x_k} \Phi(z_{\ell}, t) \dd t \int_{B_{\ell}} \nabla \left( x_j - \Upsilon_j^{(\ell)}\left( x \right) \right) \cdot \nabla x_k \dd x\\
			& = - \varepsilon^d \sum_{\ell = 1}^m (\gamma_0 - \gamma_{\ell}) \int_0^T \sum_{j = 1}^d \sum_{k = 1}^d \partial_{x_j} U(z_{\ell}, t) \partial_{x_k} \Phi(z_{\ell}, t) \dd t \int_{\partial B_{\ell}} x_k \partial_n \xi_j^{(\ell)} \dd \sigma_x\\
			& = - \varepsilon^d \sum_{\ell = 1}^m (\gamma_0 - \gamma_{\ell}) \int_0^T \nabla U(z_{\ell}, t) \cdot M^{(\ell)} \nabla \Phi(z_{\ell}, t) \dd t.
		\end{align*}
		Combining the preceding estimates with Lemma \ref{lem:est-trunc} completes the proof of the theorem.
	\end{proof}

\section{Direct reconstruction algorithms}\label{sec:alg}

Now we develop two direct algorithms for recovering small conductivity inclusions, and discuss relevant theoretical underpinnings.

\subsection{Asymptotic expansions of fundamental solutions}
The construction of the proposed algorithms relies on the fundamental solution to the model \eqref{eq:tfd:u}. For any $(x_0,t_0)\in\mathbb{R}^d \times[0,T]$, let $\Psi_{(x_0,t_0)}$ be the fundamental solution satisfying
\begin{equation*}
\left\{\begin{aligned}
	\p_t^\alpha \Psi_{(x_0,t_0)}(x,t) - \gamma_0 \Delta \Psi_{(x_0,t_0)}(x,t) &= \delta_{(x_0,t_0)}(x,t),\quad(x,t)\in \mathbb{R}^d\times[0,T],\\
	\Psi_{(x_0,t_0)}(x,0) &=0,\quad x\in \mathbb{R}^d,\\
\Psi_{(x_0,t_0)}(x,t) &\quad\mbox{decays to zero as }|x|\to\infty,
\end{aligned}\right.
\end{equation*}
where $\delta_{(x_0, t_0)}$ denotes the Dirac delta function at $(x_0, t_0)$.
Let $\Phi(x,t)=\Psi_{(x_0,0)}(x,T-t)$, with $x_0\in\mathbb{R}^d \backslash\overline{\Om}$. Then by changing variables $s'=T-s$, we have
\begin{align*}	
{}_t\widetilde{\partial}_{T}^{\alpha} \Phi (x, t)
&= \frac{1}{\Gamma(1-\alpha)}\int_0^{T - t}\frac{\p_{s'} \Psi_{(x_0,0)}(x,s')}{(T-t-s')^\alpha}\,\ds'=\p_t^\alpha \Psi_{(x_0,0)}(x,T-t).
\end{align*}
Thus, $\Phi$ satisfies \eqref{eqn:PDE-Phi}, and $I_{\Phi}$ admits the asymptotic expansion \eqref{eq:asymptotic-exp} in Theorem \ref{thm:asymptotic-exp}. However, the accurate evaluation of $\Psi_{(x_0,t_0)}$ is very challenging; see \cite{aceto2022efficient} for the 1D case (using Wright function). Thus, one has to use approximate fundamental solutions in the algorithms so as to ensure their  computational tractability, and this greatly complicates the analysis of the algorithms.
Following \cite{Qiu:2024:DSM}, we define the reduced Green function $\Psi_{d,\alpha}$ by
\begin{align} \label{eq: fundamental solution and reduced Green function}
  \Psi_{(x_0,t_0)}(x,t)=\Psi_{d,\alpha}\left(\frac{x-x_0}{\sqrt{\gamma_0 (t-t_0)^\alpha}}\right)\cdot(\gamma_0 (t-t_0)^\alpha)^{-\frac{d}{2}},
\end{align}
and we employ the following result from \cite[Theorems 3.1 and 3.6]{Qiu:2024:DSM}. Also  the explicit expressions for the coefficients $a_0$ and $a_{j,k}$ for $j=1,2$ and $k=0,1,2$ can be found in \cite[Remarks 3.2, 3.7 and Tables 1--2]{Qiu:2024:DSM}.
\begin{lemma}
\label{lemma:QS:funsolasympt}
For any $N\in\mathbb{N}^+$, there exist constants $a_0>0$, $a_{j,k}\in\mathbb{R}$ with $j=1,2$ and $0 \le k < N$ such that the reduced Green functions $\Psi_{d,\alpha}(x)$ for $d = 1, 2, 3$ admit the following asymptotic expansions as $|x| \to \infty$:
\begin{align*}
	\Psi_{1,\alpha}(x)&=\exp\left(-a_0|x|^{\frac{2}{2-\alpha}}\right)|x|^{-\frac{1-\alpha}{2-\alpha}}\Big(\sum_{k=0}^{N-1} a_{1,k} |x|^{\frac{-2k}{2-\alpha}} + \mathcal{O}(|x|^{\frac{-2N}{2-\alpha}})\Big),\\
	\Psi_{2,\alpha}(x)&=\exp\left(-a_0|x|^{\frac{2}{2-\alpha}}\right)|x|^{-\frac{2-2\alpha}{2-\alpha}}\Big(\sum_{k=0}^{N-1} a_{2,k} |x|^{\frac{-2k}{2-\alpha}} + \mathcal{O}(|x|^{\frac{-2N}{2-\alpha}})\Big),\\
	\Psi_{3,\alpha}(x)&=-(2\pi|x|)^{-1}\Psi_{1,\alpha}'(|x|).
\end{align*}
Moreover, the asymptotics of $\nabla\Psi_{d,\alpha}$ can be obtained by termwise differentiation of $\Psi_{d,\alpha}$.
\end{lemma}

We denote by $\Psi_{d,\alpha,N}(x)$ the truncated asymptotic expansion of $\Psi_{d,\alpha}(x)$ with $N$ terms (i.e., dropping the remainder $\mathcal{O}(|x|^{\frac{-2N}{2-\alpha}})$ in Lemma \ref{lemma:QS:funsolasympt}),
and define the approximate fundamental solution $\Psi_{(x_0,t_0),N}$ with $N$ terms by
\begin{equation}\label{eq:Psi:nd}		
 \Psi_{(x_0,t_0),N}(x,t):=\Psi_{d,\alpha,N}\left(\frac{x-x_0}{\sqrt{\gamma_0 (t-t_0)^\alpha}}\right)\cdot(\gamma_0 (t-t_0)^\alpha)^{-\frac{d}{2}}.
\end{equation}
Below we take the test function $\Phi(x,t)=\Psi_{(x_0,0), N}(x,T-t)$. Note that $\Phi=\Psi_{(x_0,0), N}(x,T-t)$ no longer satisfies \eqref{eqn:PDE-Phi} for any $N$ and $x_0$, and thus Theorem \ref{thm:asymptotic-exp} does not directly apply to $I_{\Phi}$. We give suitable conditions on $N$ and $x_0$ to guarantee that the perturbation of $I_{\Phi}$ is comparable with the remainder in \eqref{eq:asymptotic-exp} in Theorem \ref{thm:asymptotic-exp} and that $I_{\Phi}$ admits a similar asymptotic expansion.
We shall see that this approximation requires the source point $x_0$ to be far away from $\Om$, but $N$ is allowed to be small.
The next result gives the expression of $\nabla\Psi_{(x_0,0),N}(x,t)$.
\begin{lemma}\label{lem:Phi-deriv}
For each fixed $N\in\mathbb{N}^+$ and $d = 2, 3$,  $\nabla\Psi_{(x_0,0),N}(x,t)$ is given by
\begin{align*}
\nabla\Psi_{(x_0,0),N}(x,t)&=(\gamma_0 t^\alpha)^{-\frac{d+1}{2}}\nabla\Psi_{d,\alpha,N}\left(\frac{x-x_0}{\sqrt{\gamma_0 t^\alpha}}\right),\quad \mbox{with } \nabla\Psi_{d,\alpha,N}(x)=x S_{d,N}(|x|^2),	
\end{align*}
where the auxiliary function $S_{d,N}$ is defined by for $y>0$,
\begin{align*}
  S_{2,N}(y)&=2\exp(-a_0 y^{\frac{1}{2-\alpha}})\sum_{k=0}^{N-1}a_{2,k}\left(-\frac{a_0}{2-\alpha}y^{\frac{\alpha-1}{2-\alpha}}+\frac{\alpha-1-k}{2-\alpha}y^{-1}\right)y^{\frac{\alpha-1-k}{2-\alpha}},\\
   {S}_{3,N}(y) &=\! -(2\pi)^{-1}\exp\left(\!-a_0 y^{\frac{1}{2-\alpha}}\right)\sum_{k=0}^{N-1} a_{1,k}\left[\left(\frac{2a_0}{2-\alpha}\right)^2 y^{\frac{5\alpha-5-2k}{2(2-\alpha)}}\! +\! \frac{8a_0(k+1-\alpha)}{(2-\alpha)^2} y^{\frac{5\alpha-7-2k}{2(2-\alpha)}}\right]\\
&\quad-(2\pi)^{-1}\exp\left(-a_0 y^{\frac{1}{2-\alpha}}\right)\sum_{k=0}^{N-1} a_{1,k}\frac{(2k+1-\alpha)(2k+5-3\alpha)}{(2-\alpha)^2} y^{\frac{5\alpha-9-2k}{2(2-\alpha)}}.	
\end{align*}
\end{lemma}
\begin{proof}
 The desired identities follow directly from \eqref{eq: fundamental solution and reduced Green function} and Lemma \ref{lemma:QS:funsolasympt}.
	\end{proof}
	
Below we employ the boundary measurements $I_\Phi$ for two choices of the background solution $U$: harmonic functions and approximate fundamental solutions, and develop direct reconstruction algorithms based on these measurement data. We also provide relevant theoretical underpinnings. To indicate the dependence on $U$, we use the notation $I_\Phi(U)$.

\subsection{Reconstructing one inclusion using  harmonic functions}\label{ssec:single}
	
First we consider one inclusion for which we choose the background solution $U$ to be a harmonic function (i.e.,  $f\equiv 0$), with harmonic $u_0$ and time-independent $g$ such that $U(x,t)  = U(x) = u_0(x)$. This choice has its origin in the constant current method \cite{KwonSeoYoon:2002}. We only discuss one circular inclusion, but the algorithm also works for inclusions of other shapes, cf. Section \ref{sec:exp}. To recover one circular inclusion of radius $\varepsilon$ centered at $z_1$, we take  $\Phi(x,t)=\Psi_{(x_0,0), N}(x,T-t)$ as the test function, with the truncation level $N \in \mathbb{N}^+$ and the source point $x_0 \in \mathbb{R}^d\backslash\overline{\Om}$ properly chosen, and use $d$ measurements corresponding to $d$ background solutions $U(x, t)  = U(x) = u_0^j(x) = a_j \cdot x$, $j = 1, \ldots, d$, where the vectors $\{a_j\}_{j=1}^d\subset \mathbb{R}^d$ are linearly independent.

First we specify the assumption on $N$ and $x_0$ in $\Phi(x,t)=\Psi_{(x_0,0),N}(x,T-t)$.
\begin{assumption} \label{ass:N-T}
$N\in\mathbb{N}^+$ and $x_0\in\mathbb{R}^d\backslash\overline{\Om}$ satisfy
\begin{align}
	& N \ge \lceil\delta |\log\varepsilon|\rceil, \label{eq: N in assumption} \\
 C_0^{-1}\exp\left(C_0'\delta^{-1}\right) \le & \operatorname{dist}(x_0,\Om)\le{C_0}\exp\left(C_0'\delta^{-1}\right), \label{eq:source:howfar:ti}
\end{align}
where the constants $C_0>1$, $C_0'>\frac{(2-\alpha)d}{4}$ and $0< \delta \ll 1$ are independent of $\varepsilon$.
\end{assumption}

Below let $b_2= 2, b_3 =\frac{5}{2}$, and for $d=2,3$,
\begin{equation*}
 g_d(x,t):= \tau(x,t)^{b_d}\exp(-a_0\tau(x,t)),
 \end{equation*}
 with
 \begin{equation*}
 \tau(x,t) := \left( \frac{|x-x_{0}|^2}{\gamma_0 (T-t)^\alpha} \right)^{\frac{1}{2-\alpha}}, \quad \forall (x,t)\in\Om\times (0,T).
\end{equation*}
Note that $a_0\in(0,1)$ by \cite[Remarks 3.2, 3.7]{Qiu:2024:DSM}, and thus over $\tau \in (0, \infty)$, the function $\tau \mapsto \tau^{b_d}\exp(-a_0\tau)$ is positive and attains the maximum at $\tau=\tfrac{b_d}{a_0}$ for fixed $b_d$. Then for some constant $C>0$, we have
\begin{align} \label{eq:bdd-Int}
0 < g_d(x,t) < C, \quad \forall (x,t)\in\Om\times (0,T).
\end{align}
Now we give an estimate of the perturbation due to the truncation of the approximate fundamental solution.
\begin{lemma}\label{lemma:epsilon:dover2:estimate}
    If Assumption \ref{ass:N-T} holds for $x_0$ and $N$ with sufficiently small $\delta$,
    then
    \begin{equation}\label{eq:epsilon:dover2:estimate}
    \sup_{(x,t)\in\Om\times(0,T)}\left|(\nabla \Psi_{(x_0,0)}-\nabla \Psi_{(x_0,0),N})(x,T-t)\right| = \mathcal{O}(\varepsilon^{\frac{d}{2}}).
    \end{equation}
\end{lemma}
\begin{proof}
    From \eqref{eq: fundamental solution and reduced Green function} and \eqref{eq:Psi:nd}, we have, for all $(x,t)\in\Om\times(0,T)$,
    \begin{align}\label{eq:fundsol:difference:gradients}
        &\left|(\nabla \Psi_{(x_0,0)}-\nabla \Psi_{(x_0,0),N})(x,T-t)\right|\\
        =&|x-x_0|^{-d-1}\tau(x,t)^{\frac{(2-\alpha)(d+1)}{2}}\left|(\nabla \Psi_{d,\alpha}-\nabla \Psi_{d,\alpha,N})(\tau(x,t)^{\frac{2-\alpha}{2}})\right|.\nonumber
    \end{align}
    To estimate the right-hand side, by Lemmas \ref{lemma:QS:funsolasympt} and \ref{lem:Phi-deriv}, we derive, as $\tau\to\infty$,
    \begin{equation}\label{eq:reducedfundsol:difference}
    |(\nabla \Psi_{d,\alpha}-\nabla \Psi_{d,\alpha,N})(\tau^{\frac{2-\alpha}{2}})|=\begin{cases}
        \mathcal{O}\left(\exp(-a_0\tau)\tau^{\frac{2-\alpha}{2}}\tau^{2\alpha-2-N}\right),&\mbox{ in 2D},\\
        \mathcal{O}\left(\exp(-a_0\tau)\tau^{\frac{2-\alpha}{2}}\tau^{\frac{5\alpha-5-2N}{2}}\right),&\mbox{ in 3D}.
    \end{cases}
    \end{equation}
Note that $\frac{(2-\alpha)(d+2)}{2}+b_d(\alpha-1)=b_d$ and $\sup_{\tau>0}\exp(-a_0\tau)\tau^{b_d}<\infty$ for $d=2,3$.
Also the estimate \eqref{eq:source:howfar:ti} in Assumption \ref{ass:N-T} gives
$$\tau(x,t)\ge\tau(x,0)\ge \tau_0 := \widetilde{C}\exp \Bigg( \frac{2}{2-\alpha}C_0'\delta^{-1} \Bigg),$$
for $\widetilde{C} > 0$ independent of $\delta$. By combining these relations with \eqref{eq:fundsol:difference:gradients} and \eqref{eq:reducedfundsol:difference}, we obtain, for sufficiently small $\delta$,
    $$
    \sup_{(x,t)\in\Om\times(0,T)}\left|(\nabla \Psi_{(x_0,0)}-\nabla \Psi_{(x_0,0),N})(x,T-t)\right|=\mathcal{O}(\tau_0^{-N}).
    $$
Finally, the condition \eqref{eq: N in assumption} and the condition $\frac{2}{2-\alpha}C_0'>\frac{d}{2}$ in Assumption \ref{ass:N-T} imply that, for sufficiently small $\delta$,
$$\tau_0^{-N} \le \tau_0^{-\delta|\log\varepsilon|}  = \exp( - ( \tfrac{2}{2-\alpha}C_0' + \delta \log \widetilde{C}) |\log \varepsilon|)  = \mathcal{O}(\varepsilon^{\frac{d}{2}}),$$ which gives the desired estimate \eqref{eq:epsilon:dover2:estimate}.
\end{proof}

The following theorem gives an asymptotic expansion of $I_{\Phi}$ under Assumption \ref{ass:N-T}.
\begin{theorem}
\label{thm:2d:1inc:ti}
Let Assumption \ref{ass:N-T} hold for $x_0$ and $N$ with sufficiently small $\delta$.
Then for any harmonic function $U$, $I_{\Phi}$ with $\Phi(x,t)=\Psi_{(x_0,0), N}(x,T-t)$ satisfies
\begin{align} \label{eq:I:asymp:ti}
I_{\Phi}= - \varepsilon^d  (\gamma_0 - \gamma_{1}) \int_0^T \nabla U(z_1) \cdot M^{(1)} \nabla {\Phi(z_1, t)}
\,{\rm d}t + \widetilde{R}_d(\varepsilon),
\end{align}
with $\widetilde{R}_d(\varepsilon) = \mathcal{O}(\varepsilon^{d+1}\ell_{\varepsilon,d}^{\frac12})$.
Furthermore, with $r:=-(\gamma_0-\gamma_1)\nabla U(z_1)\cdot\frac{z_1 -x_0}{ |z_1 - x_0|}$, there exist $C_1$, $C_2$, and $C_3>0$ independent of $\varepsilon$ such that	
\begin{align}\label{eq:I:est:ti}
&C_1 r \varepsilon^d - C_3 \varepsilon^{d+1} \ell_{\varepsilon,d}^{1/2} \le I_{\Phi} \le C_2 r\varepsilon^d + C_3 \varepsilon^{d+1} \ell_{\varepsilon,d}^{\frac12}.
\end{align}
\end{theorem}

\begin{proof}
{Let $\widetilde{\Phi}(x,t) = \Psi_{(x_{0},0)}(x,T-t)$.}
First, under Assumption \ref{ass:N-T}, by \eqref{eq:bdd-Int} and the known asymptotics of $\Psi_{(x_0,t)}$ \cite[Section 2]{EidelmanKochubei:2004}, we have
$\nabla\Psi_{(x_{0},0)}(x,T-t)\simeq\mathcal{O}(1)$ uniformly in $(x,t)\in\Om\times(0,T)$,
and thus the assertion of Theorem \ref{thm:asymptotic-exp} holds for $I_{\widetilde{\Phi}}$.
Next, using Lemma \ref{lem:est-trunc} and the trace lemma, we obtain the following estimate on $V$ defined in \eqref{eq: auxiliary asymptotic expansion}:
$$\|V\|_{L^2(\p\Om\times(0,T))}^2=\mathcal{O}({_0I_T^{\alpha}} \|V(\cdot, t) \|_{H^1(\Omega)}^2)=\mathcal{O}(\ell_{\varepsilon,d}\varepsilon^{d+2}).$$ Moreover, the asymptotic relation $\Upsilon(x)=\mathcal{O}(|x|^{1-d})$ as $|x|\to\infty$ gives $$\|\varepsilon\nabla_x U(z_1,t)\cdot\Upsilon(\varepsilon^{-1}(x-z_1))\|_{L^2(\p\Om\times(0,T))}^2=\mathcal{O}(\varepsilon^2\varepsilon^{d-1}\varepsilon^{2(d-1)})=\mathcal{O}(\varepsilon^{3d-1}).$$
Combining these relations gives
$$\|u-U\|_{L^2(\p\Om\times(0,T))}^2=\mathcal{O}(\ell_{\varepsilon,d}\varepsilon^{d+2}).$$
Then under Assumption \ref{ass:N-T} with sufficiently small $\delta$, Lemma \ref{lemma:epsilon:dover2:estimate} yields
\begin{align*}
&\Big|\int_0^T\!\!\! \int_{\partial \Omega} \gamma_0 (u - U)(x, t) \partial_n ({\widetilde{\Phi}(x,t) - \Phi(x,t)}) \dd \sigma_x \dd t \Big|
\le C (\ell_{\varepsilon,d} \varepsilon^{d+2})^{\frac12}\varepsilon^{\frac{d}{2}}= C \ell_{\varepsilon,d}^{\frac12} \varepsilon^{d+1},\\
&\Big|\varepsilon^d (\gamma_0 - \gamma_{1}) \int_0^T \nabla U(z_1, t) \cdot M^{(1)} \nabla ({\widetilde{\Phi}(z_1, t) - \Phi(z_1,t)})\, {\rm d}t \Big|
\le C \varepsilon^{d}\varepsilon^{\frac{d}{2}}=C\varepsilon^{\frac{3d}{2}}.
\end{align*}
These and Theorem \ref{thm:asymptotic-exp} give \eqref{eq:I:asymp:ti}.

Next, Lemmas \ref{lemma:QS:funsolasympt} and \ref{lem:Phi-deriv} with $N=1$ and \eqref{eq:source:howfar:ti} imply
\begin{align}\label{eq:dphi:ti}
\nabla \Psi_{(x_0,0)}(x,T-t)
=&\begin{cases}
\displaystyle-\frac{2a_0 a_{2,0}}{2-\alpha}  \frac{x-x_{0}}{|x-x_{0}|^4}g_2(x,t)\left(1+\delta_{2}(\tau(x,t))\right),&d=2,\\[3mm]
\displaystyle-\frac{a_{1,0}}{2\pi}\left(\frac{2a_0}{2-\alpha}\right)^2\frac{x-x_{0}}{|x-x_{0}|^3}g_3(x,t)\left(1+\delta_{3}(\tau(x,t))\right),&d=3,
\end{cases}
\end{align}
with 
{$\max\{ |\delta_2(\tau)|, |\delta_3(\tau)| \} \le C\exp\left(-C'\delta^{-1}\right)$}
for some $C$ and $C'$ independent of $\tau$ and $\delta$.
Moreover, by \eqref{eq:dphi:ti} with $x=z_1$ and \eqref{eq:I:asymp:ti} with $M^{(1)}$ the polarization tensor of one circular inclusion, cf. \eqref{eqn:polar-circ}, we have
\begin{equation*}
\begin{aligned}
	&I_\Phi+ \mathcal{O}(\varepsilon^{d+1}\ell_{\varepsilon,d}^{\frac{1}{2}})\\
   =&\begin{cases}
\displaystyle\varepsilon^2\frac{2a_0a_{2,0}}{2-\alpha}r|m_{11}^{(1)}| |z_1 - x_0|^{-3}\int_0^T g_2(z_1,t)\left(1+\delta_2(\tau(z_1,t))\right)\,{\rm d} t, &d=2,\\
\displaystyle\varepsilon^3\frac{a_{1,0}}{2\pi}\left(\frac{2a_0}{2-\alpha}\right)^2 r|m_{11}^{(1)}||z_1-x_{0}|^{-2}\int_{0}^Tg_3(z_1,t)\left(1+\delta_{3}(\tau(z_1,t))\right)\,{\rm d} t, &d=3.
\end{cases}
\end{aligned}	
\end{equation*}
The relation \eqref{eq:source:howfar:ti} in Assumption \ref{ass:N-T} implies $|z_1 - x_0|\simeq\mathcal{O}(1)$.
Also note that for $\delta>0$ small enough,
{we have $\max\{ |\delta_2(\tau)|, |\delta_3(\tau)| \} < 1$.}
These relations and \eqref{eq:bdd-Int}
yield the inequality \eqref{eq:I:est:ti}.
\end{proof}

\begin{figure}[hbt!]
    \centering
    \setlength{\tabcolsep}{0pt}
    \begin{tabular}{cc}
	    \includegraphics[width=0.45\linewidth]{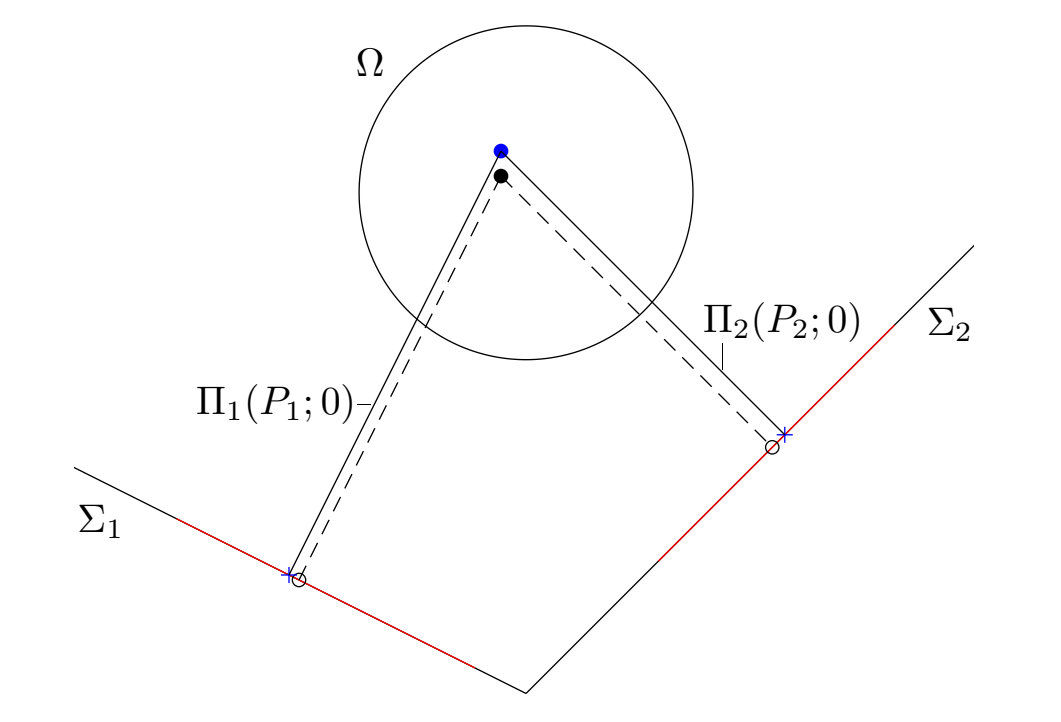}&
    \includegraphics[width=0.45\linewidth]{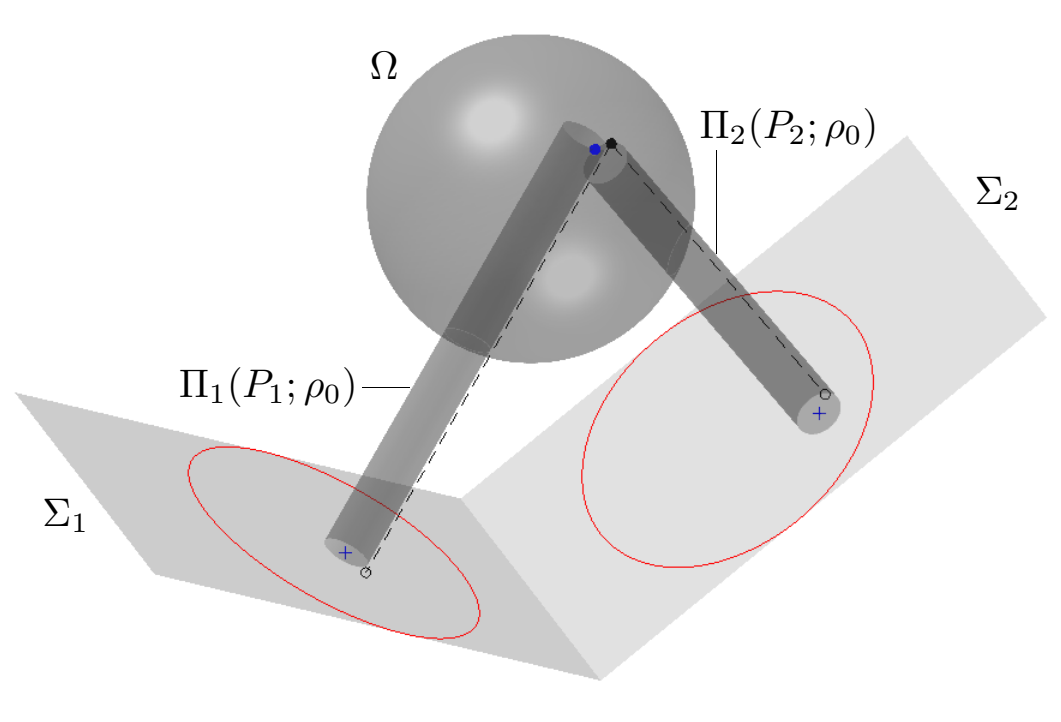}\\
    (a) 2D & (b) 3D
    \end{tabular}
    \caption{Illustration of the algorithm for locating one circular inclusion. The black solid dot denotes the target inclusion. The solutions $P_1$ and $P_2$ to \eqref{eq:reconst:bytwolines} are located at the blue $+$ marks, and the reconstruction $P$ is indicated by the blue solid dot. The red line segments in (a) and the interiors of red circles in (b) indicate the orthogonal projections of $\Omega$ onto $\Sigma_1$ and $\Sigma_2$.}
    \label{fig:algorithm for 1 inclusion}
\end{figure}

Now we describe a direct algorithm for locating one circular inclusion.
We have $d$ measurements, corresponding to $U(x, t) = U(x) = u_0^j(x) = a_j \cdot x$ for $j = 1, \ldots, d$.
For $j = 1, 2$, let ${\Sigma}_j$ be a line parallel to $a_j$ in 2D, and be a hyperplane orthogonal to $a_j\times a_3$ in 3D. By the linear independence of $\{a_j\}_{j=1}^d$, we have $\Sigma_1 \nparallel \Sigma_2$. The first step of the algorithm is to find the source point $P_j\in \Sigma_j$, $j=1,2$, by solving the following equation for some $N \in \mathbb{N}^+$:
\begin{equation}\label{eq:reconst:bytwolines}
I_{\Phi}(U)=0, \mbox{ with }\begin{cases}
	    (\Phi(x,t), U(x))=(\Psi_{(P_j,0), N}(x,T-t),u_0^j(x)),&\mbox{in }2D,\\
     (\Phi(x,t), U(x))=(\Psi_{(P_j,0), N}(x,T-t),u_0^j(x)) \\
       \quad\mbox{ and }(\Psi_{(P_j,0), N}(x,T-t),u_0^3(x)),&\mbox{in }3D.
	\end{cases}
	\end{equation}
The existence of a solution $P_j\in\Sigma_j$ to \eqref{eq:reconst:bytwolines} inside the orthogonal projection of $\Om$ onto $\Sigma_j$ for $j=1,2$ is given in Theorem \ref{thm:rec1inc}(a) below.
Let $\Pi_j(P_j;\rho)$ be the line if $\rho=0$ and the infinite cylinder of radius $\rho$ if $\rho>0$, endowed with the axis of symmetry $\{x\in\mathbb{R}^d\,:\,(x-P_j)\perp \Sigma_j\}$, $j=1,2$.
Let
$$\rho_0:=\min\left\{\rho\ge0\,:\,\Pi_1(P_1;\rho)\bigcap\Pi_2(P_2;\rho)\ne\emptyset\right\}.$$
Note that there always holds $\rho_0=0$ in 2D (i.e., $\Pi_j(P_j;\rho_0)$ is a line for $j = 1, 2$). The recovered  location  $P$ is given by
\begin{equation}\label{eq:1inc:recon:def}
 \{P\}:=\Pi_1(P_1;\rho_0) \bigcap \Pi_2(P_2;\rho_0).
\end{equation}
See Fig. \ref{fig:algorithm for 1 inclusion} for schematic illustrations.
We shall give an error estimate for this algorithm in Theorem \ref{thm:rec1inc}(b). This algorithm in 2D is inspired by that in \cite[Section 5.1]{ammari2004reconstruction}. We have extended it to the model \eqref{eq:tfd:u} using $\Phi_{(x_0,t_0),N}$ and presented a new algorithm in 3D.

Next we discuss theoretical properties of the algorithm. We first give a result that motivates the construction of the algorithm. Let
\begin{equation}\label{eq:Itild:art}
\widetilde{I}_{\Phi}:=- \varepsilon^d  (\gamma_0 - \gamma_{1}) \int_0^T \nabla U(z_1) \cdot M^{(1)} \nabla \Phi(z_{1}, t) \,{\rm d}t
\end{equation}
be the leading-order term in \eqref{eq:I:asymp:ti}, which approximates $I_\Phi$ for small $\varepsilon$, cf. Theorem \ref{thm:2d:1inc:ti}.
The following result shows that the algorithm can exactly recover the location $z_1$ of the inclusion, if we use $\widetilde{I}_{\Phi}$ instead of $I_{\Phi}$ when solving \eqref{eq:reconst:bytwolines}. Note that for $\widetilde{I}_\Phi$, the remainder $\widetilde{R}_d(\varepsilon)$ is absent, and the impact of $\widetilde{R}_d(\varepsilon)$ on the reconstruction $P$ will be studied in Theorem \ref{thm:rec1inc}.

\begin{proposition}\label{prop:asympexp:rec1inc}
Let $\Sigma_j \subset \mathbb{R}^d\backslash\overline{\Om}$. Then for $N=1$, there exists a unique pair of solutions $(P_1,P_2)\in\Sigma_1\times\Sigma_2$ to \eqref{eq:reconst:bytwolines} with $I_\Phi$ replaced by $\widetilde{I}_\Phi$. Furthermore, we have $$\Pi_1(P_{1};0)\bigcap\Pi_2(P_{2};0)=\{z_1\}.$$
The results also hold for any $N \in \mathbb{N}^+$ under the additional assumption $\operatorname{dist}(\p\Om,\Sigma_j)\gg1$.
\end{proposition}
\begin{proof}
By  \cite[Remarks 3.2, 3.7 and Tables 1--2]{Qiu:2024:DSM}, the coefficients $a_0$, $a_{1,0}$ and $a_{2,0}$
are all positive.
Then there holds $S_{d,N}(y)<0$ if $N=1$ and $y>0$, or $N \in \mathbb{N}^+$ and $y\gg1$.
Then by the expression of $M^{(1)}$ in \eqref{eqn:polar-circ} (for a circular inclusion), the condition $\widetilde{I}_{\Phi}(u_0^j)=0$ (and additionally $\widetilde{I}_{\Phi}(u_0^3)=0$ in 3D) for $\Phi(x,t)=\Psi_{(P_j,0), N}(x,T-t)$ is equivalent to $a_j\cdot(z_1 - P_{j})=0$ (and also $a_3\cdot(z_1 - P_{j})=0$ in 3D). This implies that $P_{j}$ is the orthogonal projection of $z_1$ onto $\Sigma_j$, and hence $\Pi_1(P_{1};0)\bigcap\Pi_2(P_{2};0)=\{z_1\}$.
\end{proof}

Let $\widetilde{\Sigma}_j$ be the orthogonal projection of $\Om$ onto $\Sigma_j$ for $j=1,2$.
In practice, the algorithm solves \eqref{eq:reconst:bytwolines} using the data $I_\Phi(u_0^j)$, $j=1,\ldots,d$, with $\Phi(x,t)=\Psi_{(x_0,0),N}(x,T-t)$ for $x_0\in\widetilde{\Sigma}_j$ to obtain $P_1$ and $P_2$, and to locate $P$ as \eqref{eq:1inc:recon:def}.
The following error bound on the recovered inclusion location holds.
\begin{theorem}\label{thm:rec1inc}
Let $\Sigma_j \subset \mathbb{R}^d\backslash\overline{\Om}$, and Assumption \ref{ass:N-T} hold for $N$ and all $x_0\in\widetilde{\Sigma}_j, j = 1, 2$, with sufficiently small $\delta$. There exist $C,\varepsilon_0>0$ such that the following assertions hold for all $\varepsilon\in(0,\varepsilon_0)$:
\begin{itemize}
    \item[\rm(a)] There exists a pair of solutions $(P_1,P_2)\in\widetilde{\Sigma}_1\times\widetilde{\Sigma}_2$ to \eqref{eq:reconst:bytwolines}.
    \item[\rm(b)] The reconstruction $P$ defined in \eqref{eq:1inc:recon:def} for every pair $(P_1,P_2)$ in part (a) satisfies
$$| P - z_1 |\le C\varepsilon\ell_{\varepsilon,d}^{\frac12}.$$
\end{itemize}
\end{theorem}
\begin{proof}
Let $r=r(x_0,U)$ be as in Theorem \ref{thm:2d:1inc:ti}. It follows from \eqref{eq:I:est:ti} that there exists $C>1$ such that for all $P_j\in\widetilde{\Sigma}_j$, with $\Phi(x,t)=\Psi_{(P_j,0),N}(x,T-t)$, $j=1,2$, there hold
\begin{align}
\label{eq:conditions:itmvalthm}
&I_\Phi(u_0^j)\begin{cases}
	>0,&\mbox{if }r(P_j,u_0^j)>C\varepsilon\ell_{\varepsilon,d}^{\frac12},\\
	<0,&\mbox{if }r(P_j,u_0^j)<-C\varepsilon\ell_{\varepsilon,d}^{\frac12},
\end{cases} \\
& \mbox{ and also in 3D, } I_\Phi(u_0^3)\begin{cases}
	>0,&\mbox{if }r(P_j,u_0^3)>C\varepsilon\ell_{\varepsilon,d}^{\frac12},\\
	<0,&\mbox{if }r(P_j,u_0^3)<-C\varepsilon\ell_{\varepsilon,d}^{\frac12}.
\end{cases}\nonumber
\end{align}
Next we discuss the 2D and 3D cases separately.\\
(i) 2D case. Let $F_I^j(P_j) = I_{\Phi}(u_0^j)$, with $\Phi(x,t)=\Psi_{(P_j,0),N}(x,T-t)$ for $j = 1, 2$, and $\Lambda_j$ be the line segment
\begin{equation}\label{eq:root:bounds}
	\Lambda_j = \big\{P_j\in\widetilde{\Sigma}_j\,:\,|r(P_j,u_0^j)|\le C\varepsilon\ell_{\varepsilon,2}^{\frac12}\big\}.
\end{equation}
Then $F_I^j(P_j)$ is continuous on $\Lambda_j$.
By \eqref{eq:conditions:itmvalthm} and the intermediate value theorem, $F_I^j(P_j) = 0$ has at least one solution in $\Lambda_j$.
By the definition of $r(P_j,u_0^j)$, as $\varepsilon \rightarrow 0$, the line segment $\Lambda_j$ shrinks to $\operatorname{proj}(z_1\to\Sigma_j)$,  the projection of $z_1$ onto the line $\Sigma_j$, and it is an interior point of $\widetilde{\Sigma}_j$, so part (a) holds in 2D. Let $(P_1,P_2)$ be any such pair of solutions.
Since $P_j$ is in the region \eqref{eq:root:bounds} for $j=1,2$, the reconstruction $P$ defined by \eqref{eq:1inc:recon:def} satisfies
\begin{equation*}
	 \left|\operatorname{proj}((P-z_1)\to\Sigma_j)\right|=\left|\frac{a_j}{|a_j|}\cdot(P-z_1)\right|=\frac{|z_1-P_j|}{|\gamma_0-\gamma_1||a_j|}\left|r(P_j,u_0^j)\right|=\mathcal{O}(\varepsilon\ell_{\varepsilon,2}^{\frac12}),\quad j=1,2.
\end{equation*}
This and the condition $\Sigma_1 \nparallel \Sigma_2$ implies the desired assertion in (b).\\
(ii) 3D case. For each $j=1,2$, we consider the solid parallelogram
\begin{equation}\label{eq:parallelogram}
    \widetilde{\Lambda}_j = \left\{ P_j\in\widetilde{\Sigma}_j\,:|a_k\cdot(P_j-z_1)|\le \operatorname{dist}(\Om,\Sigma_j)|\gamma_0-\gamma_1|^{-1}C\varepsilon, \quad k = j,3 \right\}.
\end{equation}
For each $j=1,2$, the function $\widetilde{F}_I^j(P_j) = ( I_\Phi(u_0^j), I_\Phi(u_0^3) )$, with $\Phi(x,t)=\Psi_{(P_j,0),N}(x,T-t)$, is continuous on $\widetilde{\Lambda}_j$.
Note that for $j=1,2$, every $P_j \in \partial \widetilde{\Lambda}_j = \{P_j \in \widetilde{\Sigma}_j: |a_k\cdot(P_j-z_1)|= \operatorname{dist}(\Om,\Sigma_j)|\gamma_0-\gamma_1|^{-1}C\varepsilon, \quad \forall k = j, 3\}$ satisfies $|r(P_j,u_0^j)|\ge C\varepsilon$ and $|r(P_j,u_0^3)|\ge C\varepsilon$.
Then by \eqref{eq:conditions:itmvalthm} and Poincar\'{e}-Miranda theorem \cite{Kulpa:1997}, $\widetilde{F}_I^j(P_j) = (0, 0)$ has at least one solution in $\widetilde{\Lambda}_j$.
As $\varepsilon\rightarrow 0$, the parallelogram $\widetilde{\Lambda}_j$ shrinks to $\operatorname{proj}(z_1\to\Sigma_j)$, which is again an interior point of $\widetilde{\Sigma}_j$, so part (a) also holds in 3D.
Let $(P_1,P_2)$ be any such pair of solutions.
Since $P_j \in \widetilde{\Lambda}_j$ for $j=1,2$, the reconstruction $P$ satisfies
\begin{equation*}
	 \left|\operatorname{proj}((P-z_1)\to\Sigma_j)\right|=\mathcal{O}\left(\left|\frac{a_j}{|a_j|}\cdot(P-z_1)\right|+\left|\frac{a_3}{|a_3|}\cdot(P-z_1)\right|\right)
=\mathcal{O}(\varepsilon).
\end{equation*}
This and the condition $\Sigma_1 \nparallel \Sigma_2$ imply the desired assertion in part (b).
\end{proof}

\begin{remark}
Theorem \ref{thm:rec1inc} requires linearly independent directions $a_j$. The analysis shows that mutually orthogonal vectors $a_j$ lead to best conditioning for the reconstruction from the projection, and nearly parallel $a_j$ may impair the reconstruction accuracy.
\end{remark}

\subsection{Reconstructing  multiple inclusions using approximate fundamental solutions}\label{ssec:multiple}

Now we consider the case of recovering multiple small circular inclusions. To this end, we choose the test function $\Phi(x, t) = \Psi_{(x_1,0), N}(x, T - t)$ and the background solution $U(x,t) = \Psi_{(x_2,0), N}(x, t)$ for some $N \in \mathbb{N}^+$ with finitely many and properly chosen $x_1, x_2 \in \mathbb{R}^d \backslash \overline{\Om}$. The basic idea of the approach is related to the factorization method and MUSIC algorithm \cite{Kirsch:2008} in the sense that we shall construct an approximate indicator function to indicate the presence of small inclusions. We choose a truncation level $N \in \mathbb{N^+}$ and source points $\{ x_{j, n} \}_{j = 1}^n \subset \mathbb{R}^d \backslash \overline{\Om}$ suitably, and obtain the data
\begin{equation}\label{eq:datamatrix:realistic}
B_I = \bigg[ I_{\Phi}(U)\Big{|}_{U(x, t)=\Psi_{(x_{j',n},0), N}(x, t), \Phi(x, t) = \Psi_{(x_{j'',n},0), N}(x, T - t)} \bigg]_{j',j''=1}^n \in \mathbb{R}^{n\times n},
\end{equation}
with which we aim to reconstruct the locations of multiple inclusions. In the idealized case, i.e. choosing $\Phi(x, t) = \Psi_{(x_1,0)}(x, T - t)$ and $U(x,t) = \Psi_{(x_2,0)}(x, t)$ for infinitely many $x_1, x_2$, we can prove a sufficient criterion for locating inclusions. This important observation motivates the construction of an indicator function to indicate the inclusion locations in practice, using the finite amount of data \eqref{eq:datamatrix:realistic}.
{We focus the discussion on circular inclusions, but the algorithm is also applicable to inclusions of other shapes; see the numerical experiments in Section \ref{sec:exp}.}

First we discuss the idealized case. Let $\Gamma$ be a nonempty open subset of $\{x\in\mathbb{R}^d\,:\,|x|=R\}$, with $R>\max\{|x|\,:\,x\in\overline{\Om}\}$.
Let $\mathcal{L}_z$ and $\mathcal{M}$ be the Hilbert-Schmidt integral operators defined on $L^2(\Gamma)$, respectively, by
\begin{align*}
	\mathcal{L}_{z}\left[f\right](y) &:=   \int_\Gamma f(x)\!\!\int_0^T\!\!\nabla\Psi_{(x,0)}(z,t)\cdot\nabla\Psi_{(y,0)}(z,T-t)\,{\rm d}t{\rm d}\sigma_x,\quad\forall f\in L^2(\Gamma),\, y\in\Gamma,\, z\in\Om,\\
		\mathcal{M}[f](y)&:=\int_\Gamma f(x) I_{\Psi_{(y,0)}(\cdot,T-\cdot)}(U=\Psi_{(x,0)})\,{\rm d}\sigma_x,\quad\forall f\in L^2(\Gamma),\, y\in\Gamma.
\end{align*}
The operator $\mathcal{L}_z$ serves as an idealized probing function to locate small inclusions from the measurement data. In view of the integral representation of $I_\Phi$ in Lemma \ref{lem:equiv-wbm}, the ranges of the operators $\mathcal{L}_z$ and $\mathcal{M}$ consist of analytic functions that analytically extend to $\mathbb{R}^d\backslash \{z\}$ and $\mathbb{R}^d\backslash \overline{A}$, respectively.
We denote their extensions (to $\mathbb{R}^d\backslash \{z\}$ and $\mathbb{R}^d\backslash \overline{A}$, respectively) also by $\mathcal{L}_z$ and $\mathcal{M}$. The next result gives a sufficient criterion for locating small conductivity inclusions. In the statement, we let
$$
Z:=\{z_{\ell}\}_{\ell=1}^m=\bigcap_{\varepsilon>0}A\quad\mbox{and}\quad B_{\rho }({z}):=\{y: |y-{z}|<\rho\},
$$
a ball centered at $z$ with a radius $\rho$. Clearly $Z$ is the set of inclusion centers. The proof of the theorem is technical and hence it is deferred to Appendix \ref{app:thm:multi-char}.

\begin{theorem}\label{thm:multi-char}
For any $z\in\Om$, we have $z\in Z$ if, for every $f\in L^2(\Gamma)$, the following identity  holds
\begin{align}\label{eq:error:oneside}
\lim_{\varepsilon \to 0} \inf_{g\in L^2(\Gamma)} \lim_{\rho\to0} \left|\int_{B_{\rho}({z})\backslash\overline{A}}\gamma_0 {\Delta_y}\left(\mathcal{L}_z[f](y)-\mathcal{M}[g](y)\right)F_{z,f}(y)\,{\rm d}y\right| = 0,
\end{align}
with
$$F_{z,f}(y):=\exp\left((y-z)\cdot\int_\Gamma  \nabla\Psi_{(x,0)}(z,T)f(x)\,\dd \sigma_x\right).$$
Further, for all $z\in\Om$, there exists an $f\in L^2(\Gamma)$ satisfying $\|f\|_{L^\infty(\Gamma)}=1$ and some $C>0$ independent of $\varepsilon$ such that
\begin{equation}\label{eq:compare}
\lim_{\rho\to0}\int_{B_{\rho}({z})}\gamma_0 \Delta_y\mathcal{L}_z[f](y)F_{z,f}(y)\,{\rm d}y\ge C>0.
\end{equation}
\end{theorem}
	
The second identity \eqref{eq:compare} in Theorem \ref{thm:multi-char} shows that for all $z\in\Omega$, the function $\mathcal{L}_z[f]$ has a singularity (only) at $z$. The first identity \eqref{eq:error:oneside} in Theorem \ref{thm:multi-char} is obviously true for any $z\in Z$ in view of the inclusion $Z\subset A$ and shows that if the singularity of $\mathcal{L}_z[f]$ is canceled out by $\mathcal{M}[g]$ with an optimal choice of $g$ and the degenerate size of the conductivity inclusions (i.e. $\varepsilon\to0^+$), then we have $z\in Z$.
This characterization of the set $Z$ relies on the values of the functional $\mathcal{L}_z[f]$ inside $\Om$, which are determined by the analytic extension from the open set $\Gamma$. Instead of the singularity cancellation in the domain $\Omega$, we exploit the degeneracy of $\mathcal{L}_z[f]-\mathcal{M}[g]$ on the surface $\Gamma$ as a criterion for the numerical reconstruction of $Z$, which can be  approximated using the discrete data $B_I\in \mathbb{R}^{n\times n}$, cf. \eqref{eq:datamatrix:realistic}.
More precisely, we seek for points $z\in\Om$ such that for each $f\in L^2(\Gamma)$, the infimum of the functional $\|\mathcal{L}_z[f]-\mathcal{M}[g]\|_{L^2(\Gamma)}$ over $g\in L^2(\Gamma)$ tends to zero as the inclusion size $\varepsilon$ tends to zero. We restate the optimality of $g\in L^2(\Gamma)$ in terms of an orthogonal projection.
Specifically, let $P_{\mathcal{M}}$ be the orthogonal projection from $L^2(\Gamma)$ onto the closure $\overline{\mathcal{R}(\mathcal{M})}$ of the range $\mathcal{R}(\mathcal{M})$ of the operator $\mathcal{M}$.
Theorem \ref{thm:multi-char} suggests that an inclusion is located at $z$
if the mapping
$$\varphi(z) = \frac{\|\mathcal{L}_z \|_{\rm HS}}{\|({\rm id}-P_{\mathcal{M}}) \mathcal{L}_z\|_{\rm HS}}$$ blows up at $z$, where $\|\cdot\|_{\rm HS}$ denotes the Hilbert-Schmidt norm on $L^2(\Gamma)$, and ${\rm id}$ is the identity operator. Thus, it provides an indicator function for locating the inclusion centers $z$ in the idealized case.

In practice, to use the finite data  \eqref{eq:datamatrix:realistic}, we choose $N$ and $\{ x_{j, n} \}_{j = 1}^n$ suitably and construct an indicator function that approximates the mapping $\varphi(z)$.
The next result gives an asymptotic expansion of $I_\Phi$ in \eqref{eq:datamatrix:realistic}, and its proof is similar to that of Theorem \ref{thm:2d:1inc:ti}. Note that the leading-order term in \eqref{eq:I:asymp:pointsources} is a linear combination of the kernels of $\{ \mathcal{L}_{z_{\ell}} \}_{\ell = 1}^m$ {(see \eqref{eqn:polar-circ} for the expressions of $\{ M^{(\ell)} \}_{\ell = 1}^m$ for circular inclusions)} but uses $\Psi_{(x,0),N}$ instead of $\Psi_{(x,0)}$.
\begin{theorem}\label{thm:2d:1inc:pointsources}
Let $N \in \mathbb{N}^+$ and each point of $\{ x_{j, n} \}_{j = 1}^n \subset \Gamma$ satisfy Assumption \ref{ass:N-T} with sufficiently small $\delta$.
Then, for $U(x,t) = \Psi_{(x_{j'}, 0), N}(x,t)$ and $\Phi(x,t)=\Psi_{(x_{j''},0), N}(x, T - t)$ with fixed $j', j'' \in \{ 1, \ldots, n \}$, we have
\begin{align} \label{eq:I:asymp:pointsources}
&I_{\Phi}= - \varepsilon^d\sum_{\ell=1}^m (\gamma_0 - \gamma_{\ell}) \int_0^T \nabla U(z_\ell, t) \cdot M^{(\ell)} \nabla {\Phi(z_{\ell},t)}
\, {\rm d}t + \widetilde{R}_d(\varepsilon),
\end{align}
with $\widetilde{R}_d(\varepsilon)= \mathcal{O}(\varepsilon^{d+1}\ell_{\varepsilon,d}^{1/2}),$
and moreover, for some $C_1, C_2, C_3 > 0$ independent of $\varepsilon$
\begin{align*}
&C_1 \min_{1\le\ell\le m} r_\ell \varepsilon^d - C_3 \varepsilon^{d+1} \ell_{\varepsilon,d}^{1/2} \le I_{\Phi} \le  C_2 \max_{1\le\ell\leq m} r_\ell \varepsilon^d + C_3 \varepsilon^{d+1} \ell_{\varepsilon,d}^{1/2}.
\end{align*}
with
$$r_{\ell}:=(\gamma_0-\gamma_\ell)\frac{(z_{\ell}-x_{j'})\cdot(z_{\ell}-x_{j''})}{|z_{\ell}-x_{j'}||z_{\ell}-x_{j''}|}.$$
\end{theorem}

To approximate the idealized indicator function $\varphi(z)$,
we need one additional assumption on the source points $\{x_{j,n}\}_{j=1}^n \subset \Gamma$.
\begin{assumption} \label{ass:partition}
$x_{j,n} \in \Gamma_{j,n}$ for all $n \in \mathbb{N}^+$ and $j = 1, \dots, n$, where $\{\Gamma_{j,n}\}_{j=1}^n$ is a family of pairwise disjoint open subsets of $\Gamma$ satisfying $|\Gamma_{j,n}|=|\Gamma|/n$ for all $j = 1, \dots, n$ and $\bigcup_{j=1}^n\overline{\Gamma}_{j,n}=\overline{\Gamma}$ for each fixed $n$, and
$$\lim_{n \to \infty} \max_{j = 1, \ldots, n} \operatorname{diam}(\Gamma_{j,n}) = 0.$$
\end{assumption}
Since a Hilbert-Schmidt integral operator with a square-integrable kernel is compact, $\mathcal{M}$ is compact and admits Schmidt-representation \cite[Theorem 1.1]{Gohberg1990}
$$\mathcal{M}[f] = \sum_{j=1}^{\nu(\mathcal{M})} s_j(\mathcal{M})(f,v_j)_{L^2(\Gamma)} u_j, \quad\forall f\in L^2(\Gamma),$$
where $\{s_j(\mathcal{M})\}_{j=1}^{ \nu(\mathcal{M})}$ is the sequence of positive singular values of $\mathcal{M}$ (in descending order), and $\{u_j\}_{j=1}^{\nu(\mathcal{M})}$, $\{v_j\}_{j=1}^{\nu(\mathcal{M})} \in L^2(\Gamma)$ are orthonormal singular vectors.
{Since $\{u_j\}_{j=1}^{\nu(\mathcal{M})}$ forms an orthonormal basis of $\overline{\mathcal{R}(\mathcal{M})}$,}
we can define a sequence of projections onto $\overline{\mathcal{R}(\mathcal{M})}$ by
$$
P_{\mathcal{M},k}[f]:=\sum_{j=1}^k(f,u_j)_{L^2(\Gamma)}u_j,\quad\forall f\in L^2(\Gamma),\, k \in \{1,\dots,\nu(\mathcal{M})\}\cap\mathbb{N} .
$$
The construction involves two steps:
first approximate $\varphi(z)$ by
$$\varphi_k(z) = \frac{\|\mathcal{L}_z\|_{\rm HS}}{\|({\rm id}-P_{\mathcal{M},k})\mathcal{L}_z\|_{\rm HS}},$$
and then approximate $\varphi_k(z)$ by a fully discrete one. The next lemma is useful in analyzing the discrete approximation of $\varphi_k(z)$.
\begin{lemma}
\label{lemma:svd:perturb}
Let $\{\mathcal{M}_n\}_{n=1}^\infty$ be a sequence of compact linear operators on $L^2(\Gamma)$ represented by
\begin{equation}\label{eq:svd:In}
	\mathcal{M}_n[f] = \sum_{j=1}^{\nu(\mathcal{M}_n)} s_j(\mathcal{M}_n)(f,v_{j,n})_{L^2(\Gamma)} u_{j,n}, \quad\forall f\in L^2(\Gamma),
\end{equation}
with positive singular values $s_j(\mathcal{M}_n)$ {\rm(}ordered nonincreasingly{\rm)} and orthonormal singular vectors $\{u_{j,n}\}_{j=1}^{\nu(\mathcal{M}_n)}$, $\{v_{j,n}\}_{j=1}^{\nu(\mathcal{M}_n)} \in L^2(\Gamma)$.
If
$$\lim_{n \to \infty} \|\mathcal{M}-\mathcal{M}_n\|_{\mathcal{L}(L^2(\Gamma))} = 0,$$ then for each $j \in \{1,\dots,\nu(\mathcal{M})\}$, we have
\begin{align*}
\max\{ \|u_j - \mathcal{T}_{1,n} u_{j,n}\|_{L^2(\Gamma)}, \|v_j - \mathcal{T}_{2,n} v_{j,n}\|_{L^2(\Gamma)} \} \le C(\operatorname{gap}_j(\mathcal{M}))^{-1} \|\mathcal{M}-\mathcal{M}_n\|_{\mathcal{L}(L^2(\Gamma))},
\end{align*}
for some unitary transformations $\mathcal{T}_{1,n}, \mathcal{T}_{2,n} \in \mathcal{L}(L^2(\Gamma))$, where the constant $C$ depends only on the multiplicity of $s_j(\mathcal{M})$, and $\operatorname{gap}_j(\mathcal{M})$ is the distance between $s_j(\mathcal{M})$ and its nearest distinct singular value of $\mathcal{M}$.
\end{lemma}
\begin{proof}
By \cite[Corollary 1.6 in Chapter VI]{Gohberg1990}, we have
\begin{equation}\label{I:singval:converge}
\sup_{1\le j\le\nu(\mathcal{M})}|s_j(\mathcal{M})-s_j(\mathcal{M}_n)|\le\|\mathcal{M}-\mathcal{M}_n\|_{\mathcal{L}(L^2(\Gamma))}\to0,\quad \mbox{as } n\to\infty.
\end{equation}
Let $m_\ell$ be the multiplicity of the $\ell$-th singular value of $\mathcal{M}$, $J_\ell:=\{j\in\mathbb{N}\,:\,\sum_{i=1}^{\ell-1} m_i<j\le \sum_{i=1}^{\ell} m_i\}$ and $s_\ell:=s_j(\mathcal{M})$ for $j\in J_\ell$. From \eqref{I:singval:converge} and \cite[Corollary 9]{Crane:2020:ASV}, we deduce that for every $\ell$, there exists an $n_0\in\mathbb{N}$ such that for all $n\ge n_0$, $k\in J_\ell$,
		\begin{align*}
			\bigg\| u_k - \sum_{j \in J_\ell} (u_k, u_{j, n})_{L^2(\Gamma)} u_{j, n} \bigg\|_{L^2(\Gamma)} & \le 2\sqrt{2}(\operatorname{gap}_k(\mathcal{M}))^{-1} \|\mathcal{M}-\mathcal{M}_n\|_{\mathcal{L}(L^2(\Gamma))}, \\
		\bigg\| v_k - \sum_{j \in J_\ell} (v_k, v_{j, n})_{L^2(\Gamma)} v_{j, n} \bigg\|_{L^2(\Gamma)} & \le 2\sqrt{2}(\operatorname{gap}_k(\mathcal{M}))^{-1} \|\mathcal{M}-\mathcal{M}_n\|_{\mathcal{L}(L^2(\Gamma))}.
		\end{align*}
	Thus, upon choosing some unitary transformations $\mathcal{T}_{1,n}, \mathcal{T}_{2,n} \in \mathcal{L}(L^2(\Gamma))$ suitably, one can guarantee that for all $n\geq n_0$, $k\in J_\ell$,
	\begin{align*}
		 \|u_k - \mathcal{T}_{1,n} u_{k,n}\|_{L^2(\Gamma)} &\le  C(\operatorname{gap}_k(\mathcal{M}))^{-1} \|\mathcal{M}-\mathcal{M}_n\|_{\mathcal{L}(L^2(\Gamma))}, \\
		 \|v_k - \mathcal{T}_{2,n} v_{k,n}\|_{L^2(\Gamma)} &\le  C(\operatorname{gap}_k(\mathcal{M}))^{-1} \|\mathcal{M}-\mathcal{M}_n\|_{\mathcal{L}(L^2(\Gamma))},
	\end{align*}
	with the constant $C=2\sqrt{2} m_\ell$.
\end{proof}
		
To construct a discrete indicator function that approximates $\varphi_k(z)$, for $j',j''=1,\ldots,n$, we define for any $z\in \Omega$
	\begin{align*}
		C_{j',j'',n}(z)&:=\int_0^T  S_{d,N}\left(\frac{|z-x_{j',n}|^2}{\gamma_0 t^\alpha}\right) S_{d,N}\left(\frac{|z-x_{j'',n}|^2}{\gamma_0 (T-t)^\alpha}\right) (\gamma_0 t^\alpha)^{\frac{-d-2}{2}}(\gamma_0 (T-t)^\alpha)^{\frac{-d-2}{2}}\, {\rm d}t,\\
		G_n(z)&:=\left[(z-x_{j',n})\cdot  (z-x_{j'',n})C_{j',j'',n}(z)\right]_{j',j''=1}^n.
	\end{align*}
Also let $V_k\in \mathbb{R}^{n\times k}$ be the matrix consisting of the leading $k$ left singular vectors of $B_I\in\mathbb{R}^{n\times n}$ in \eqref{eq:datamatrix:realistic}, and $Q_{n,k}=I_n-V_{k} V_{k}^*$. We construct a discrete indicator function $W_{n,k}$ by
\begin{equation}\label{eq:W:definition}
	W_{n,k}(z)=\frac{\|G_n(z)\|_{\rm F}}{\|Q_{n,k} G_n(z)\|_{\rm F}},
\end{equation}
where $\|\cdot\|_{\rm F}$ denotes the Frobenius norm for matrices. In practice, the truncation level $k$ in the ratio  $W_{n,k}(z)$ (or $\varphi_k$) is determined by the decay of the singular values of $B_I$.
The next result indicates that $W_{n,k}$ well approximates $\varphi_k(z)$ for $n\gg 1$ and $\varepsilon \ll 1$.
\begin{theorem}\label{thm:convergence:indicatorftn}
Let $k\in\{1, \dots, {\nu(\mathcal{M})}\}\cap\mathbb{N}$, $N \in \mathbb{N}^+$, and $\{ x_{\ell,n} \}_{\ell = 1}^n \subset \Gamma$. Suppose that Assumption \ref{ass:N-T} holds for $N$ and each point of $\{ x_{\ell,n} \}_{\ell = 1}^n$ with sufficiently small $\delta$, and Assumption \ref{ass:partition} holds for $\{ x_{\ell,n} \}_{\ell = 1}^n$.
Then, the following identities hold:
\begin{align}
\label{eq:indicator:approximation}
\lim_{n\to\infty} \frac{|\Gamma|}{n}\|G_n(z)\|_{\rm F}&=\|\mathcal{L}_z\|_{\rm HS}+\mathcal{O}(\varepsilon),\\
\label{eq:indicator:approximation:denominator}
\lim_{n\to\infty}\frac{|\Gamma|}{n}\|Q_{n,k}G_n(z)\|_{\rm F}&=\|({\rm id}-P_{\mathcal{ M},k})\mathcal{L}_z\|_{\rm HS}+\mathcal{O}\Big(\varepsilon + \ell_{\varepsilon,d}^{1/2}\varepsilon^{d+1}\sum_{j=1}^k (\operatorname{gap}_j(\mathcal{M}))^{-1}\Big),
\end{align}
where $\operatorname{gap}_j(\mathcal{M})$ is the distance between $s_j(\mathcal{M})$ and its nearest distinct singular value of $\mathcal{M}$, and the constants in $\mathcal{O}$ in \eqref{eq:indicator:approximation}--\eqref{eq:indicator:approximation:denominator} are independent of $z$.
\end{theorem}
\begin{proof}
The integral kernels of the operators $\mathcal{L}_z$ and $\mathcal{M}$ are given, respectively, by
\begin{align*}
K(x,y;\mathcal{L}_z)&=\int_0^T\nabla\Psi_{(y,0)}(z,t)\cdot\nabla\Psi_{(x,0)}(z,T-t)\,{\rm d}t,\quad\forall x,y\in\Gamma,\, z\in\Om,\\
 K(x,y;\mathcal{M})&=I_{\Psi_{(x,0)}(\cdot,T-\cdot)}(U=\Psi_{(y,0)}),\quad\forall x,y\in\Gamma.
\end{align*}
Let $K_n(x,y;\mathcal{L}_z)$ be the piecewise constant approximation of the kernel function $K(x,y;\mathcal{L}_z)$ defined almost everywhere on $\overline{\Gamma}\times\overline{\Gamma}$ by
$$
K_n(x,y;\mathcal{L}_z) =K(x_{j',n},x_{j'',n};\mathcal{L}_z),\quad \mbox{for }
(x,y)\in\Gamma_{j',n}\times \Gamma_{j'',n},
$$
where $\Gamma_{j,n}$ is given in Assumption \ref{ass:partition}. Then, Assumption \ref{ass:partition} implies
\begin{equation}\label{eq:Kn:Lz:conv}	
\lim_{n \to \infty} \|K_n(x,y;\mathcal{L}_z)\|_{L^2(\Gamma\times\Gamma)} = \|K(x,y;\mathcal{L}_z)\|_{L^2(\Gamma\times\Gamma)} = \|\mathcal{L}_z\|_{\rm HS},
\end{equation}
with $$\|K_n(x,y;\mathcal{L}_z)\|_{L^2(\Gamma\times\Gamma)}^2=\sum_{j',j''=1}^n |K(x_{j',n},x_{j'',n};\mathcal{L}_z)|^2|\Gamma|^2/n^2.$$
By Assumptions \ref{ass:N-T} and \ref{ass:partition} and Lemmas \ref{lemma:QS:funsolasympt} and \ref{lem:Phi-deriv}, we obtain
\begin{align}\label{eq:Lz:vs:C}	
K(x_{j',n},x_{j'',n};\mathcal{L}_z)&=[G_n(z)]_{j',j''} + \mathcal{O}(\varepsilon), \quad \mbox{as } \varepsilon \to 0,
\end{align}
where the constant in $\mathcal{O}$ is independent of $j'$, $j''$, $n$ and $z\in\Om$. Combining the identity \eqref{eq:Kn:Lz:conv} with the estimate \eqref{eq:Lz:vs:C} gives the first desired relation \eqref{eq:indicator:approximation}.

Next, we define the following two operators $\mathcal{L}_{n,z}$ and $\mathcal{M}_n$ on $L^2(\Gamma)$: for any $f \in L^2(\Gamma)$,
\begin{align*}
\mathcal{M}_n[f] = & \sum_{j,k=1}^n [B_I]_{k,j} f_{j,n}\chi_{\Gamma_{k,n}},\quad \mbox{with } f_{j,n}=\frac{1}{|\Gamma_{j,n}|}\int_{\Gamma_{j,n}}f{\rm d}\sigma_y,\\
\mathcal{L}_{n,z}[f](y) = & \int_{\Gamma} f(x) \widehat{K}_{n, z}(x, y) {\rm d} \sigma_x, \quad \text{with } \widehat{K}_{n, z}(x, y) = [G_n(z)]_{ij}, \text{ for } (x,y) \in \Gamma_{i,n}\times\Gamma_{j,n}.
\end{align*}
For each $n$, $\mathcal{M}_n$ has finite rank and thus is compact. By \cite[Theorem 1.1]{Gohberg1990}, $\mathcal{M}_n$ admits the representation \eqref{eq:svd:In} for some positive sequence $s_j(\mathcal{M}_n)$ of descending order and orthonormal sequences $\{u_{j,n}\}_{j=1}^{\nu(\mathcal{M}_n)}, \{v_{j,n}\}_{j=1}^{\nu(\mathcal{M}_n)} \subset L^2(\Gamma)$. Then, the following identity holds:
\begin{align*}
\frac{|\Gamma|}{n}\|Q_{n,k} G_n(z)\|_{\rm F}&=\bigg\| \Big({\rm id}-\sum_{j=1}^k(\cdot,u_{j,n})_{L^2(\Gamma)}u_{j,n}\Big) \mathcal{L}_{n,z} \bigg\|_{\rm HS}.
\end{align*}
By the triangle inequality, we get
\begin{align} \label{eq: triangle inequality}
&\left| \frac{|\Gamma|}{n}\|Q_{n,k} G_n(z)\|_{\rm F} - \bigg\| \Big({\rm id}-P_{\mathcal{M}, k}\Big) \mathcal{L}_{z} \bigg\|_{\rm HS} \right| \\
\le & \bigg\| \Big({\rm id}-P_{\mathcal{M}, k}\Big)( \mathcal{L}_{n, z} - \mathcal{L}_{z} ) \bigg\|_{\rm HS} + \bigg\| \Big( P_{\mathcal{M}, k} - \sum_{j=1}^k(\cdot,u_{j,n})_{L^2(\Gamma)}u_{j,n} \Big) \mathcal{L}_{n,z} \bigg\|_{\rm HS}. \nonumber
\end{align}
Now we estimate the two terms separately. First, the estimate \eqref{eq:Lz:vs:C} implies
\begin{align} \label{eq: limit of integral operator}
&\limsup_{n\to\infty}\|\mathcal{L}_{z}-\mathcal{L}_{n,z}\|_{\rm HS}\\
=&{\limsup_{n\to\infty}\bigg\|K(\cdot,\cdot;\mathcal{L}_z)-\sum_{i,j=1}^n[G_n(z)]_{ij}\chi_{\Gamma_{i,n}\times\Gamma_{j,n}}(\cdot,\cdot)\bigg\|_{L^2(\Gamma\times\Gamma)}}\le C\varepsilon,\nonumber
\end{align}
with the constant $C$ independent of $z$. Thus $\| \mathcal{L}_{n,z} \|_{\rm HS}$ is uniformly bounded (independent of $n$). Meanwhile, since the kernel function $K(\cdot,\cdot;\mathcal{M})$ is smooth, Assumptions \ref{ass:N-T}--\ref{ass:partition} and Lemmas \ref{lemma:QS:funsolasympt}--\ref{lem:Phi-deriv} give
$$\lim_{n \to \infty, \varepsilon \to 0} \|\mathcal{M}-\mathcal{M}_n\|_{L^2(\Gamma)} = 0.$$
Thus, by Lemma \ref{lemma:svd:perturb}, we have, for $n\gg1$ and $\varepsilon\ll1$,
\begin{align*} 
\bigg\|P_{\mathcal{M},k}- \sum_{j=1}^k(\cdot,u_{j,n})_{L^2(\Gamma)}u_{j,n}\bigg\|_{\mathcal{L}(L^2(\Gamma))}&\le \sum_{j=1}^k 2\|u_j - u_{j,n}\|_{L^2(\Gamma)}\\
&\le \sum_{j=1}^k \frac{C}{\operatorname{gap}_j(\mathcal{M})}\|\mathcal{M}-\mathcal{M}_n\|_{\mathcal{L}(L^2(\Gamma))}.\nonumber
\end{align*}
To bound $\|\mathcal{M}-\mathcal{M}_n\|_{\mathcal{L}(L^2(\Gamma))}$, we define two kernel functions related to the operator $\mathcal{M}$:
\begin{align*}
\widetilde{K}(x,y)&:=I_\Phi(U),\quad\mbox{with }U=\Psi_{(x,0),N}\mbox{ and }\Phi(\cdot,t)=\Psi_{(y,0),N}(\cdot,T-t),\quad\forall (x,y)\in\Gamma\times\Gamma,\\
\widetilde{K}_n(x,y)&:=\widetilde{K}(x_{j,n},x_{k,n})=[B_I]_{k,j},\quad\forall(x,y)\in\Gamma_{j,n}\times\Gamma_{k,n}\quad\forall j,k\in\{1,\dots,n\}.
\end{align*}
Then, we have
\begin{align*}
\lim_{n \to \infty} \|\mathcal{M}-\mathcal{M}_n\|_{\mathcal{L}(L^2(\Gamma))} \le & \lim_{n \to \infty} \|K(\cdot,\cdot;\mathcal{M})-\widetilde{K}_n(\cdot,\cdot)\|_{L^2(\Gamma\times\Gamma)} \\
= & \|K(\cdot,\cdot;\mathcal{M})-\widetilde{K}(\cdot,\cdot)\|_{L^2(\Gamma\times\Gamma)}.
\end{align*}
In addition, by Theorem \ref{thm:2d:1inc:pointsources}, we have
\begin{align*}
\| K(\cdot,\cdot;\mathcal{M}) - \widetilde{K}(\cdot,\cdot) \|_{L^2(\Gamma\times\Gamma)} \le C \ell_{\varepsilon,d}^{1/2} \varepsilon^{d+1}.
\end{align*}
Combining the preceding three estimates gives
\begin{align} \label{eq: limit of projection}
\limsup_{n\to\infty}\bigg\|P_{\mathcal{M},k}- \sum_{j=1}^k(\cdot,u_{j,n})_{L^2(\Gamma)}u_{j,n}\bigg\|_{\mathcal{L}(L^2(\Gamma))}\le C\ell_{\varepsilon,d}^{1/2}\varepsilon^{d+1}\sum_{j=1}^k (\operatorname{gap}_j(\mathcal{M}))^{-1}.
\end{align}
Finally, combining \eqref{eq: limit of integral operator} and \eqref{eq: limit of projection}, we obtain \eqref{eq:indicator:approximation:denominator} by letting $n \rightarrow \infty$ in \eqref{eq: triangle inequality}.
\end{proof}

\section{Numerical experiments and discussions} \label{sec:exp}
Now we present numerical results to illustrate the algorithms in Section \ref{sec:alg}. Throughout we take $\Omega$ to be the unit disk $\{x\in\mathbb{R}^2\,:\,|x|<1\}$, $\alpha=0.5$, $T=1$ and $\gamma_0=1$. We generate the data $u$ using a fully discrete scheme for problem \eqref{eq:tfd:u} which employs the Galerkin finite element method with conforming linear elements in space and the L1 scheme in time \cite[Chapters 2 and 4]{JinZhou:2023book}. To accurately resolve the discontinuous conductivity, we employ finite element meshes graded around small inclusions, and divide the time interval $[0,T]$ into $2^7$ sub-intervals. The noisy data $\widetilde{u}$ is generated by $\widetilde{u} := u + \zeta$ on $\p\Om\times[0,T]$, where the noise $\zeta$ is randomly sampled in the finite element nodal values with a noise level $\sigma\geq0$, i.e.,
${\|\zeta\|_{L^1(\p\Om\times[0,T])}}/{\|u\|_{L^1(\p\Om\times[0,T])}}\sim \mathcal{N}(0,\sigma^2)$. Below we evaluate the algorithms in Sections \ref{ssec:single} and \ref{ssec:multiple} separately.

\subsection{One inclusion}
\begin{figure}[hbt!]
\centering\setlength{\tabcolsep}{0pt}
\begin{tabular}{cccc}
	\includegraphics[width=0.24\linewidth]{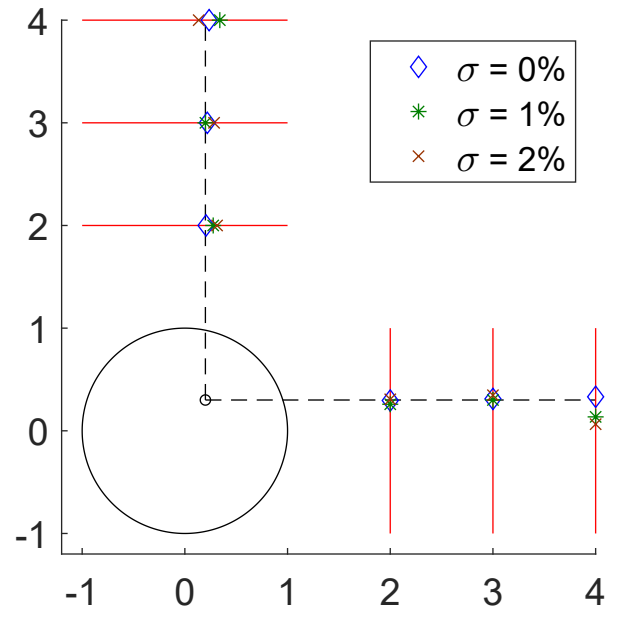}&
    \includegraphics[width=0.24\linewidth]{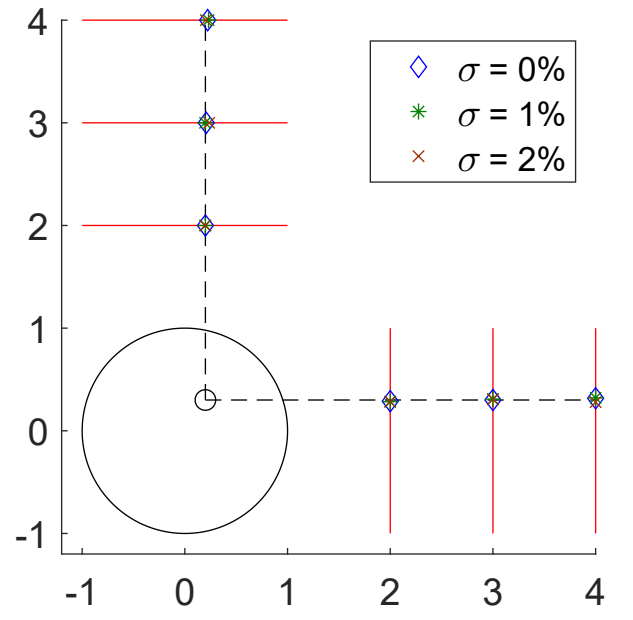}&
    \includegraphics[width=0.24\linewidth]{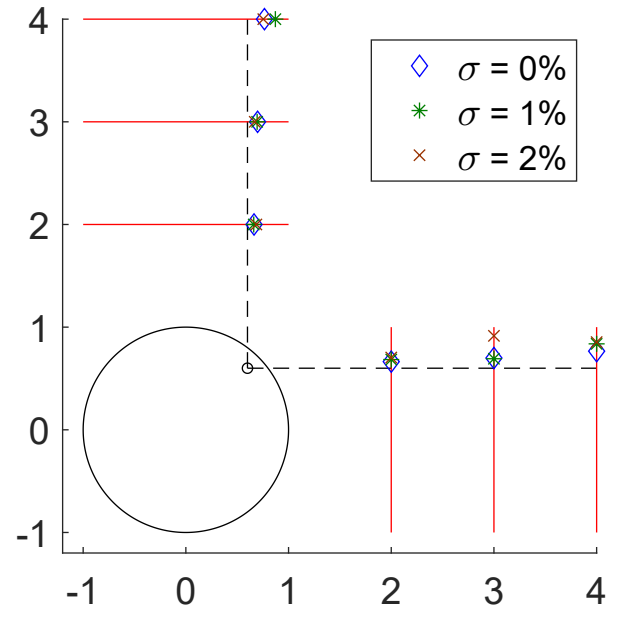}&
    \includegraphics[width=0.24\linewidth]{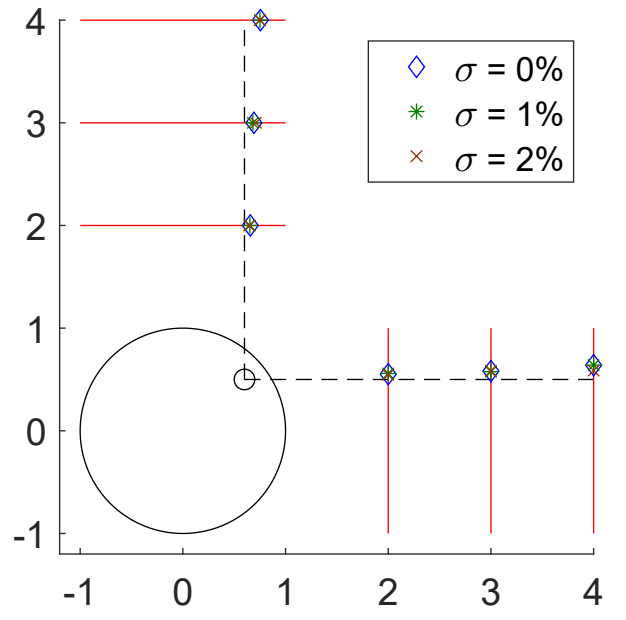}  \\
(a) ((0.2,0.3),0.05)& (b) ((0.2,0.3),0.1) & (c)((0.6,0.6),0.05) & (d) ((0.6,0.5),0.1)
\end{tabular}
\caption{\label{fig:noisy:reconstruction}
The reconstruction results for Example \ref{ex:1} with four configurations of one inclusion $(z_1,\varepsilon)$. Any pair of two non-parallel red line segments serves as the pair of domains $(\widetilde{\Sigma}_1,\widetilde{\Sigma}_2)$. The markers denote the solution $P_j$ to \eqref{eq:reconst:bytwolines}.}
\end{figure}

To recover one inclusion,
we employ the data $I_\Phi(U)$ generated with two harmonic background solutions $U_1(x,t) = (1,0) \cdot x$ and $U_2(x,t) = (0,1) \cdot x$, and the test function $\Phi(x,t)=\Psi_{(P_j,0),3}(x,T-t)$ in $I_{\Phi}(U_j)$ for $j = 1, 2$, where $P_j$ lies on the red line segment $\widetilde{\Sigma}_j$, cf. Fig. \ref{fig:noisy:reconstruction}. We test the performance of the algorithm in Section \ref{ssec:single}.

\begin{example} \label{ex:1}
Consider one circular inclusion $A = \{x\in\mathbb{R}^2\,:\,|x - z_1| < \varepsilon\}$, with  four configurations $(z_1,\varepsilon)\in\{ ((0.2,0.3),0.05),((0.2,0.3),0.1),((0.6,0.6),0.05),((0.6,0.5),0.1)\}$. The conductivity is $\gamma_1=50$ inside the inclusion.
\end{example}

\begin{figure}[hbt!]
\centering\setlength{\tabcolsep}{0pt}
\begin{tabular}{cccc}
\includegraphics[width=0.24\linewidth]{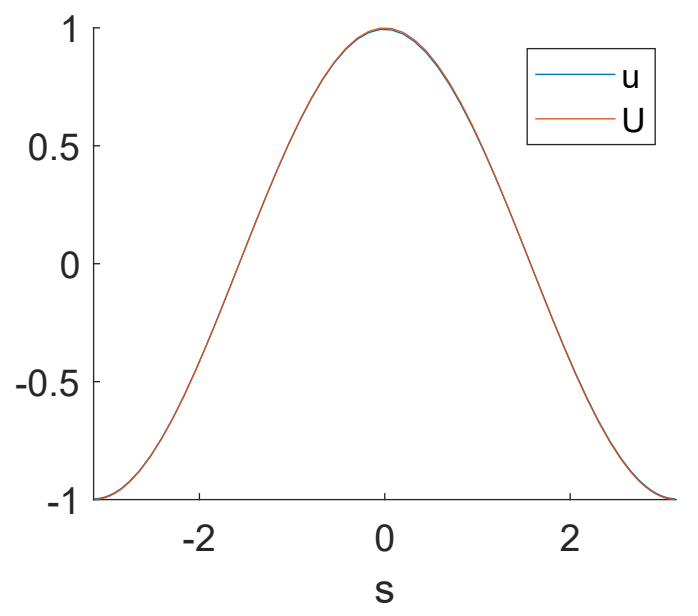}&
\includegraphics[width=0.24\linewidth]{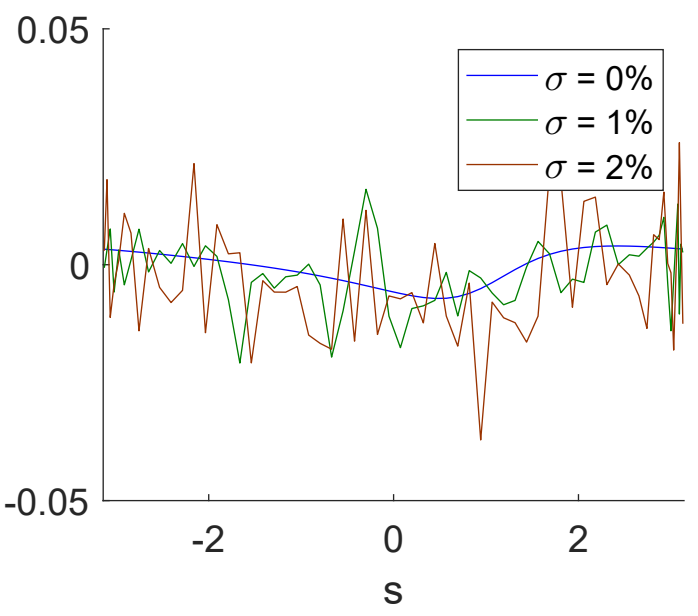}&
\includegraphics[width=0.24\linewidth]{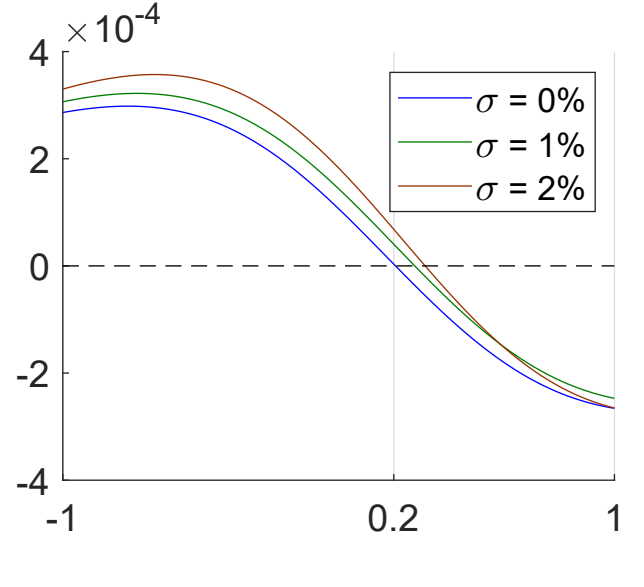} &
\includegraphics[width=0.24\linewidth]{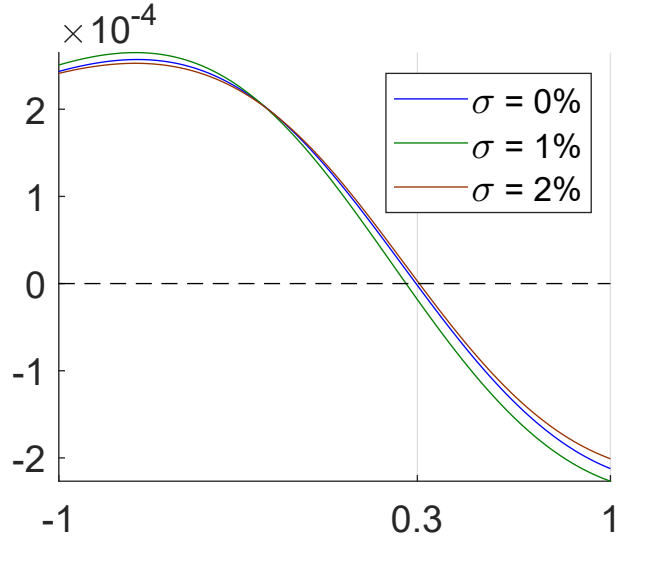}\\
\includegraphics[width=0.24\linewidth]{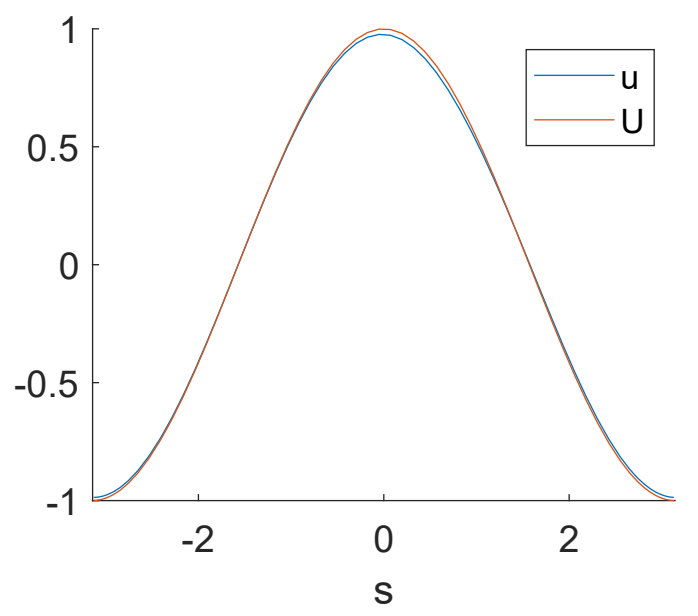}
&\includegraphics[width=0.24\linewidth]{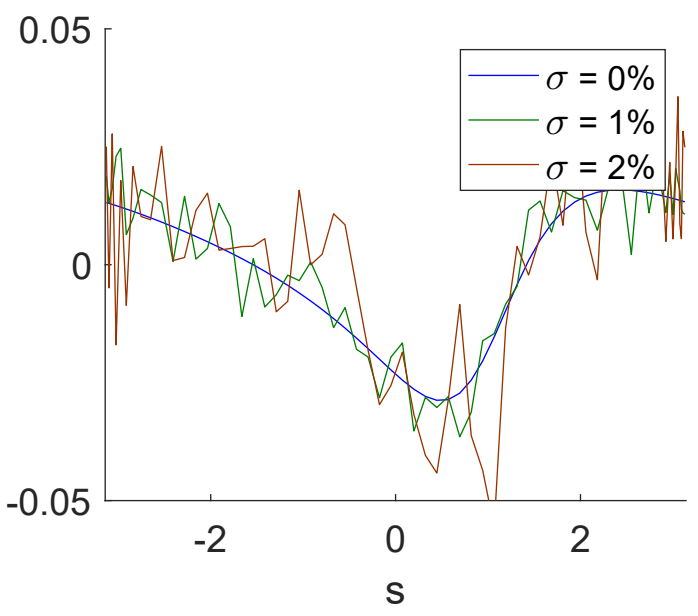}&
\includegraphics[width=0.24\linewidth]{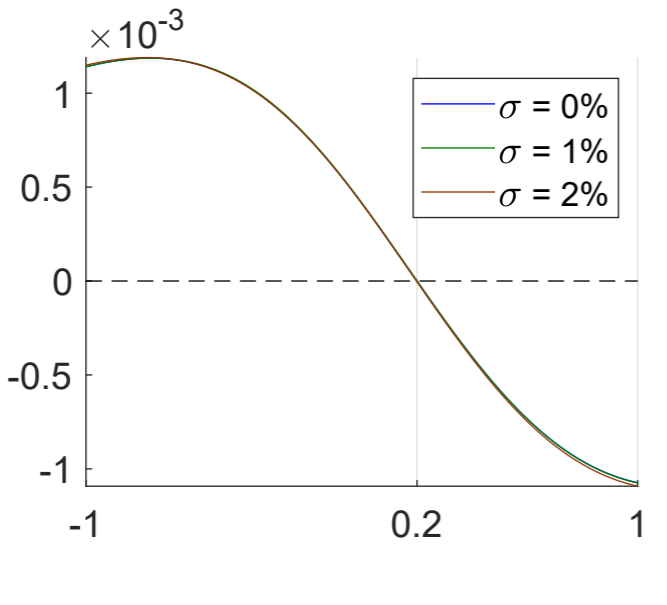}&
\includegraphics[width=0.24\linewidth]{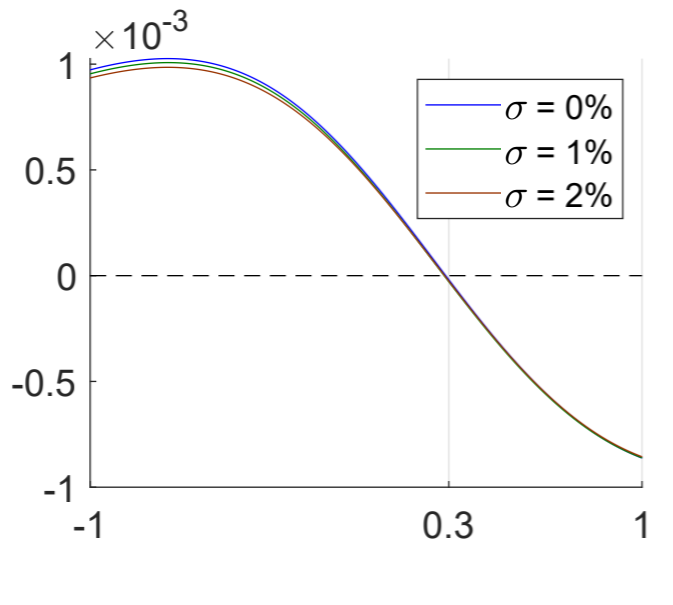}\\
(a) $u$ and $U$ & (b) $\widetilde{u} - U$ & (c) {$I_{\Phi}(U_1)$ versus $x_1$} & (d) {$I_{\Phi}(U_2)$ versus $x_2$}
\end{tabular}
\caption{\label{fig:noisy:rootfinding} The impact of different noise levels on measurement data for Example \ref{ex:1}. The top and bottom panels show the results for $(z_1,\varepsilon)=((0.2,0.3),0.05)$ and $(z_1,\varepsilon)=((0.2,0.3),0.1)$, respectively. Panels (a) and (b) show $\{u, U\}$ and $\widetilde{u} - U$ on $\p\Om$ at $T = 1$ with $U=U_1$, respectively, where $\p\Om$ is parameterized by $(\cos s,\sin s)$.
 Panels (c) and (d) show the value of $I_{\Phi}(U)$ for $U = U_1$ and $U = U_2$, respectively, with different noise levels, where for $U = U_j$, $\Phi = \Psi_{(P_j,0),3}(x,T-t)$, $j = 1, 2$, with $P_1 = (x_1, 2)$ and $P_2 = (2, x_2)$.}
\end{figure}

The reconstructions obtained by solving \eqref{eq:reconst:bytwolines} are shown in Fig. \ref{fig:noisy:reconstruction}. The inclusion locations are well resolved for both exact and noisy data. For the two cases $(z_1,\varepsilon) =((0.2,0.3),0.05)$ and $(z_1,\varepsilon) = ((0.2,0.3),0.1)$, we also show the impact of different noise levels on the data in Fig. \ref{fig:noisy:rootfinding}.
Note that the magnitude of $|u - U|$ increases with the size $\varepsilon$ of the inclusion. Thus, in the presence of data noise of a fixed level $\sigma$, the noise appears relatively large for small inclusions, in view of the $\mathcal{O}(\varepsilon^d)$ magnitude of the leading-order term in \eqref{eq:I:asymp:ti}, and it can cause larger errors when solving for $P$ from \eqref{eq:reconst:bytwolines}, i.e., a bigger reconstruction error for a smaller inclusion. However, thanks to the beneficial smoothing effect of the weighted integral of $u-U$ in $I_{\Phi}$, the reconstruction obtained by solving \eqref{eq:reconst:bytwolines} is fairly robust with respect to the additive noise in $u$. Indeed, Fig. \ref{fig:noisy:rootfinding} shows that the values of $I_\Phi$ are perturbed only very mildly even when the noise is relatively large compared with $|u - U|$.
	
\begin{example}\label{ex:1:variant}
Considers one elliptical inclusion $A = \{(a,b)\in\mathbb{R}^2\,:\,\varrho^{-1}(a-0.2)^2+ \varrho (b-0.3)^2\le0.01\}$, with an area $0.01\pi$ and various aspect ratios $\varrho\in\{2,3,5\}$. The conductivity $\gamma_1$ inside the inclusion is $=50$.
\end{example}

\begin{figure}[hbt!]
\centering\setlength{\tabcolsep}{0pt}
\begin{tabular}{ccc}
\includegraphics[width=0.25\linewidth]{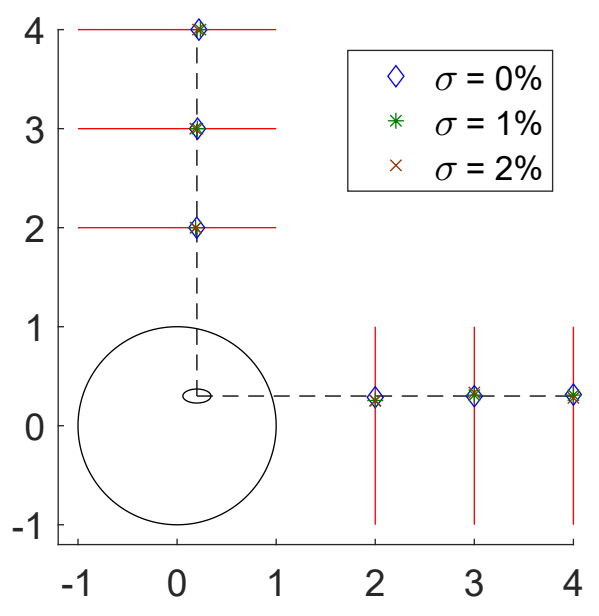}&\includegraphics[width=0.25\linewidth]{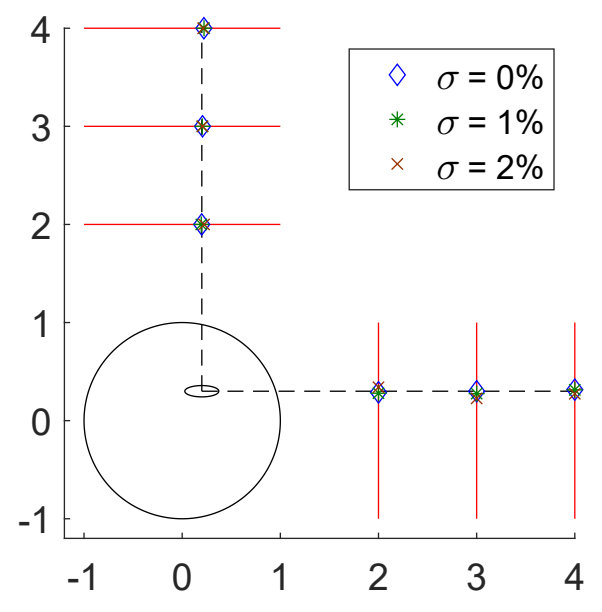}
&\includegraphics[width=0.25\linewidth]{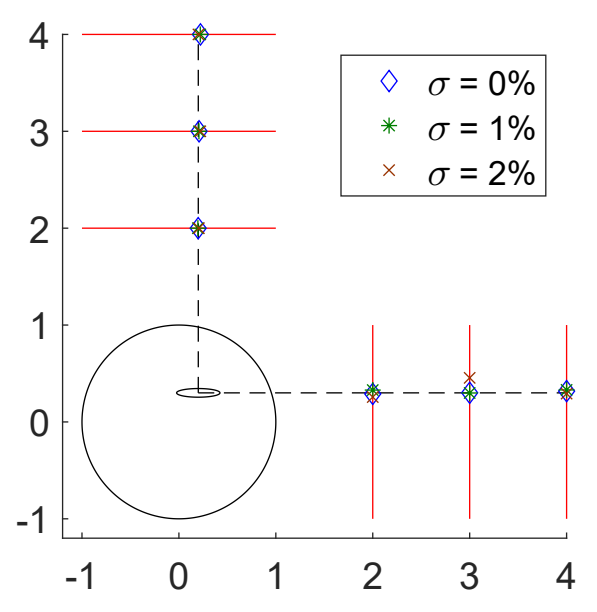}\\
(a) $\varrho=2$ & (b) $\varrho=3$ & (c) $\varrho=5$
\end{tabular}
\caption{\label{fig:one:ellipse} The reconstruction results for Example \ref{ex:1:variant} with three aspect ratios $\varrho=2$, $3$ and $5$. Any pair of two non-parallel red line segments serves as the pair of domains $(\widetilde{\Sigma}_1, \widetilde{\Sigma}_2)$. The markers denote the solution $P_j$ to \eqref{eq:reconst:bytwolines}.}
\end{figure}

\begin{figure}[hbt!]
\centering\setlength{\tabcolsep}{0pt}
\begin{tabular}{ccc}
\includegraphics[width=0.35\linewidth]{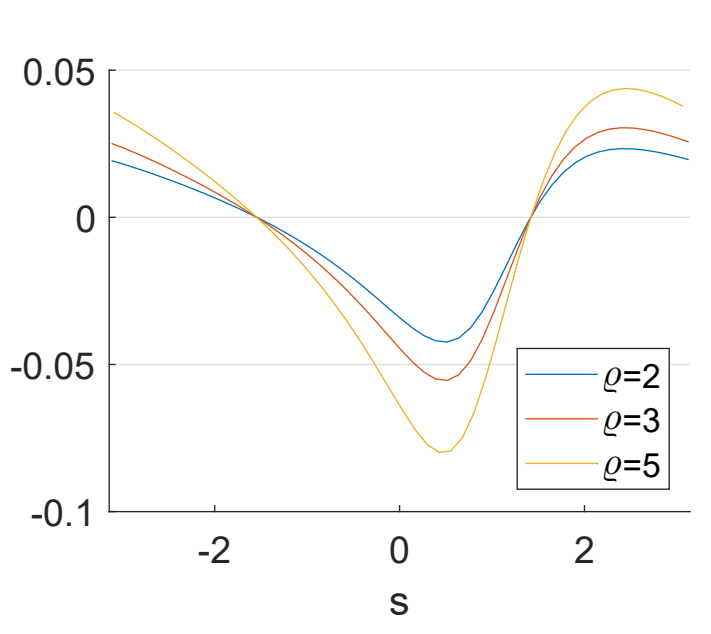}&\includegraphics[width=0.35\linewidth]{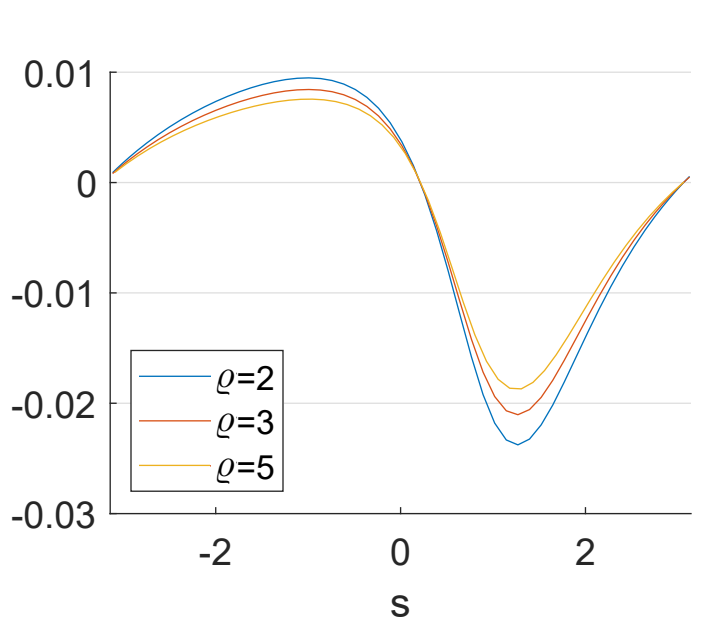}\\
(a) $U=U_1$ & (b) $U=U_2$
\end{tabular}
\caption{\label{fig:one:ellipse:u} $u - U$ for Example \ref{ex:1:variant} on $\p\Om$ at $T=1$, with $\p\Om$  parametrized by $(\cos s,\sin s)$. }
\end{figure}

Fig. \ref{fig:one:ellipse} shows the reconstruction results of one elliptical inclusion with various aspect ratios from noisy data. Interestingly, the reconstruction of the horizontal coordinate of the inclusion is more stable with respect to the noise than that of the vertical one. This can be explained by Fig. \ref{fig:one:ellipse:u}: the magnitude of $|u - U|$ is larger for $U=U_1$ than for $U=U_2$ for all aspect ratios $\varrho$. Then with a fixed level of noise, the noise is relatively larger for $I_{\Phi}(U)$ when $U = U_2$ than when $U = U_1$, and thus the solution of $P_1$ is more stable with respect to the noise than that of $P_2$.
Nonetheless, the results in all cases are fairly good, and moreover, the algorithm in Section \ref{ssec:single} is also applicable to non-circular inclusions.

\subsection{Multiple inclusions}
\label{subsec: multiple inclusions in numerical part}
\begin{figure}[hbt!]
\centering
\begin{tabular}{ccc}
\includegraphics[width=0.25\linewidth]{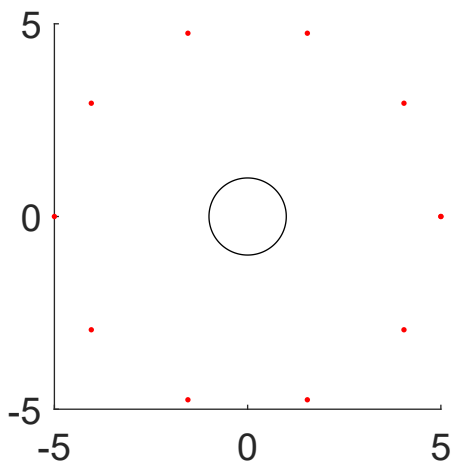}&
\includegraphics[width=0.25\linewidth]{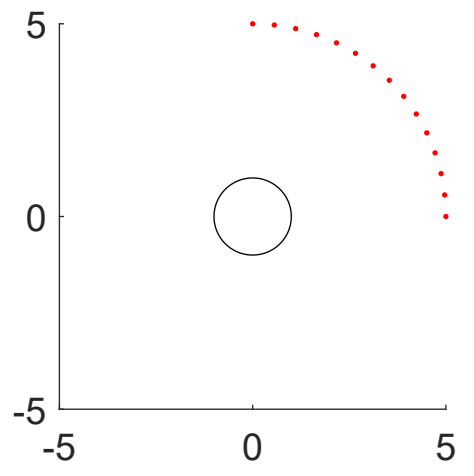}&		\includegraphics[width=0.25\linewidth]{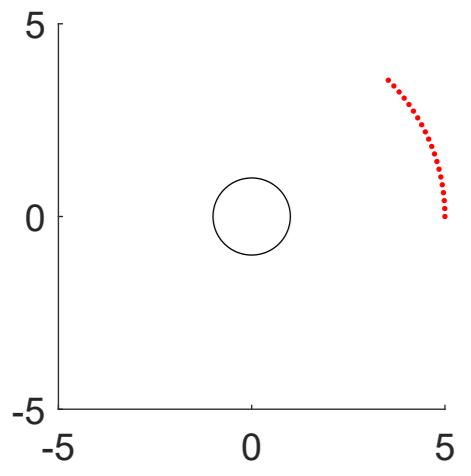}\\
(i) & (ii) & (iii)
\end{tabular}
\caption{\label{fig:sourcepts}
Schematic illustration of $\p\Om$ (black circle) and source points $\{ x_j \}_{j = 1}^n$ (red dots) in Section \ref{subsec: multiple inclusions in numerical part} for recovering multiple inclusions. }
\end{figure}

Now we consider the more challenging case of reconstructing multiple inclusions, using the algorithm in Section \ref{ssec:multiple}. The data $I_\Phi(U)$ are generated with the background solutions $U(x, t) = \Psi_{(x_{j'},0),3}(x, t)$ and the test function $\Phi(x, t) = \Psi_{(x_{j''},0),3}(x, T - t)$ for $j', j'' = 1, \ldots, n$. Below we fix $n$  at 10, and consider three configurations for the source points $\{ x_j \}_{j = 1}^n$, cf. Fig. \ref{fig:sourcepts}, including both full and limited aperture for measurement.
	
\begin{example}\label{ex:2}
Consider $m$ circular inclusions $A=\bigcup_{\ell=1}^m A_{\ell}$, with $A_{\ell} = \{x\in\mathbb{R}^2\,:\,|x-z_{\ell}|<0.05\}$ and $\gamma_{\ell}=3$ for $\ell=1,\dots,m$, in two cases: {\rm(i)} $m=2$ with $z_1=(0.3,0.2)$ and $z_2=(-0.4,0)$ and {\rm(ii)} $m=3$ with $z_1=(-0.2,0)$, $z_2=(0.2,0.3)$ and $z_3=(0.5,-0.1)$.
\end{example}
	
\begin{figure}[hbt!]
\centering
\setlength{\tabcolsep}{0pt}
\begin{tabular}{ccc}	
\includegraphics[width=0.30\linewidth]{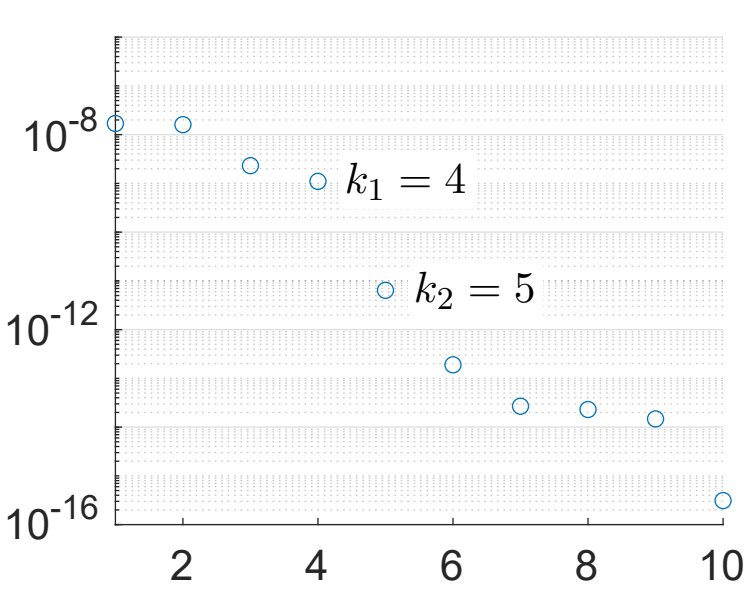}
&\includegraphics[width=0.25\linewidth]{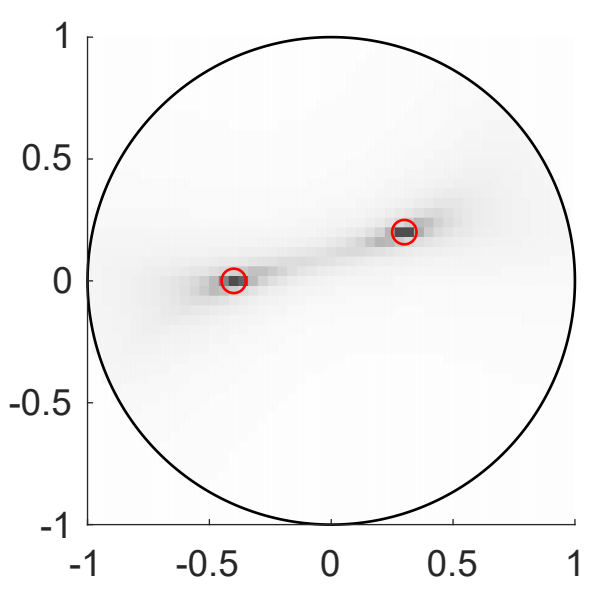}
&\includegraphics[width=0.25\linewidth]{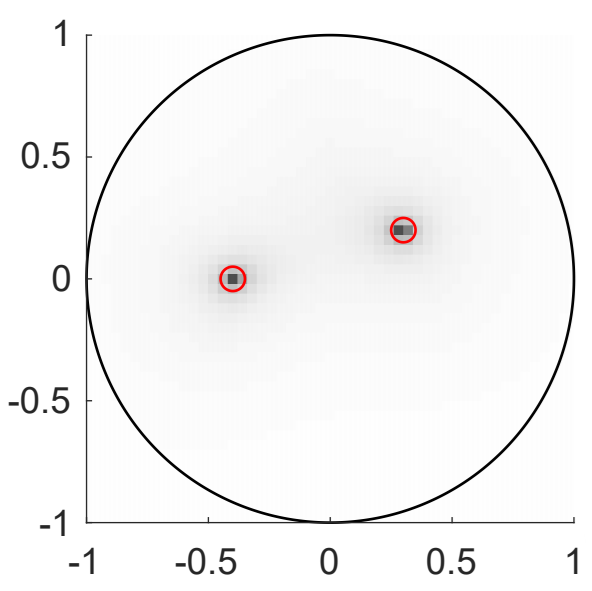}\\
\includegraphics[width=0.3\linewidth]{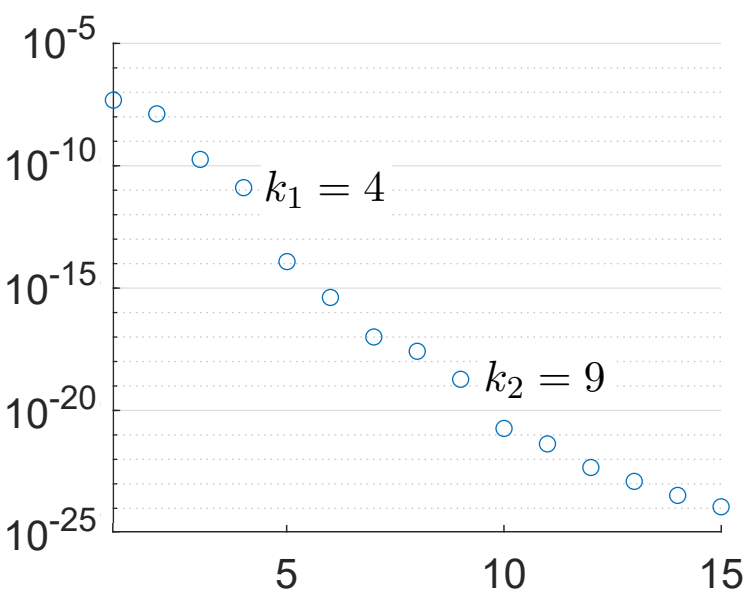}
&\includegraphics[width=0.25\linewidth]{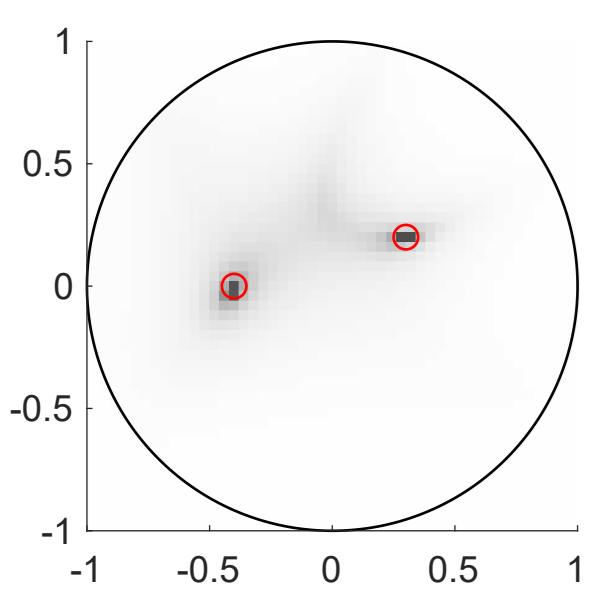}
&\includegraphics[width=0.25\linewidth]{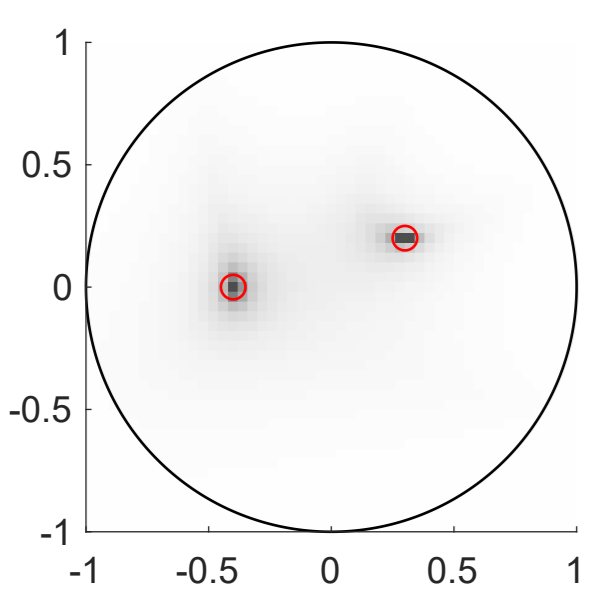}\\
\includegraphics[width=0.3\linewidth]{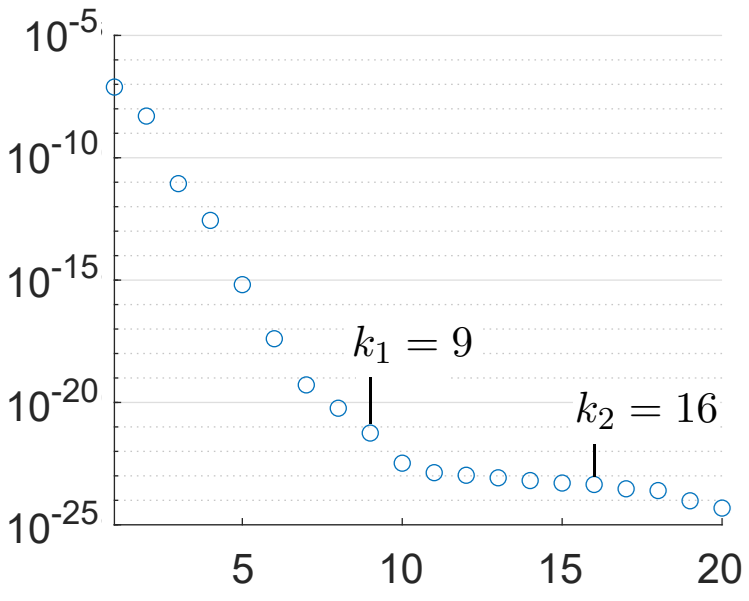}
&\includegraphics[width=0.25\linewidth]{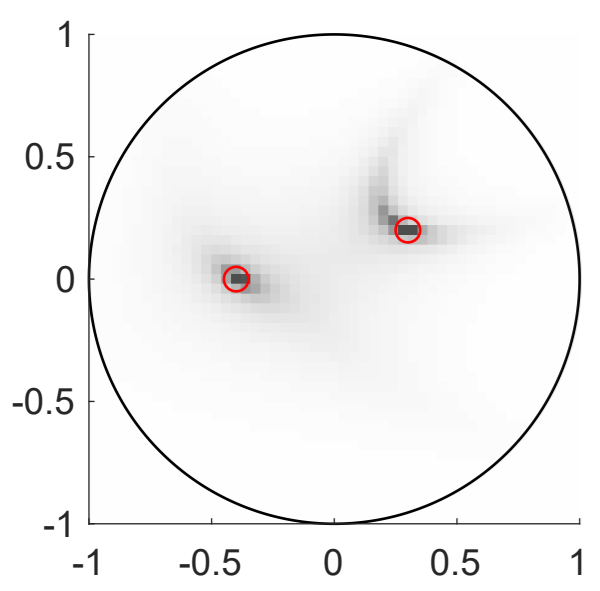}
&\includegraphics[width=0.25\linewidth]{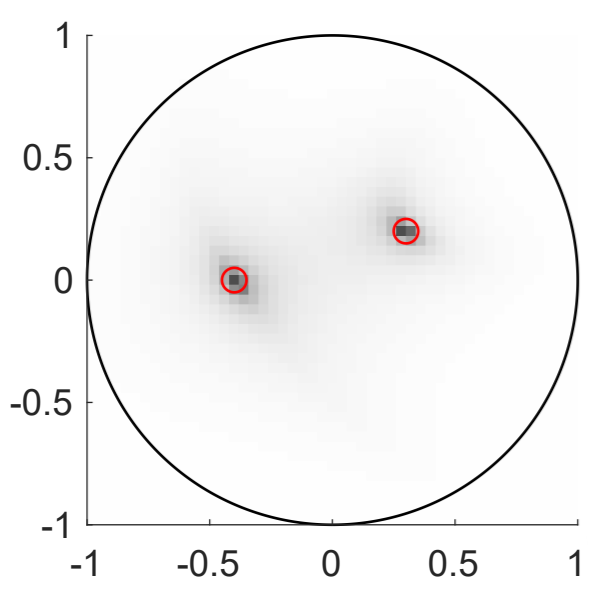}\\
(a) singular values of $B_I$ & (b) $k=k_1$ & (c) $k=k_2$
\end{tabular}
\caption{\label{fig:twoinclusion} The reconstruction results for case (i) in Example \ref{ex:2} with exact data. Panel (a) shows the singular values of $B_I$. Panels (b) and (c) show the contour plots of $W_{n,k}$ with $k=k_1$ and $k=k_2$, respectively. The red circles are the exact inclusion boundaries $\p A_1$ and $\p A_2$. The results from top to bottom correspond to the source configurations (i)--(iii) in Fig. \ref{fig:sourcepts}, respectively.}
\end{figure}
	
\begin{figure}[hbt!]
\centering
\setlength{\tabcolsep}{0pt}
\begin{tabular}{ccc}
\includegraphics[width=0.30\linewidth]{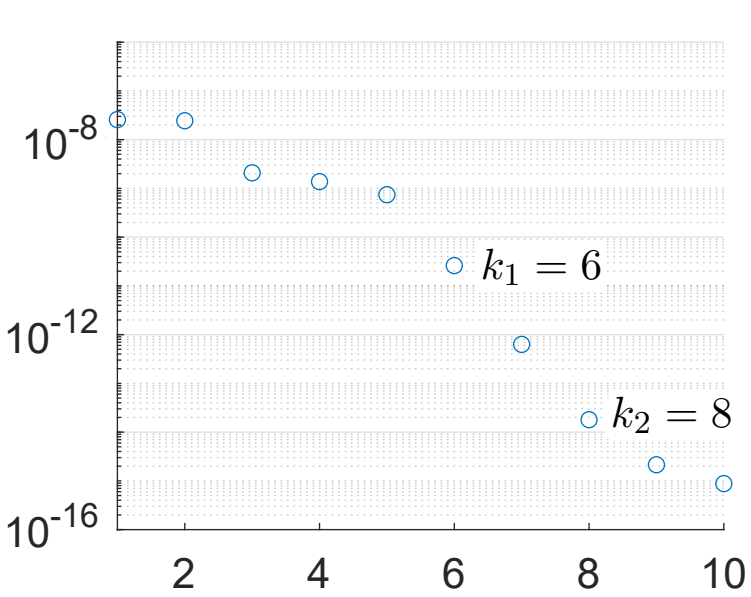}
&\includegraphics[width=0.25\linewidth]{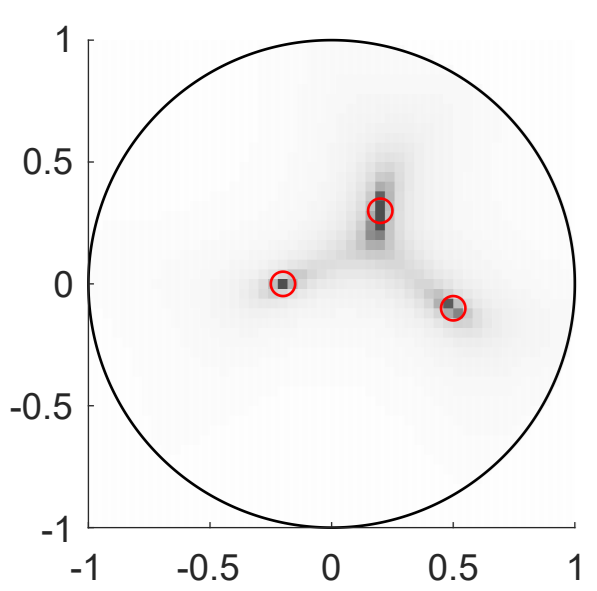}
&\includegraphics[width=0.25\linewidth]{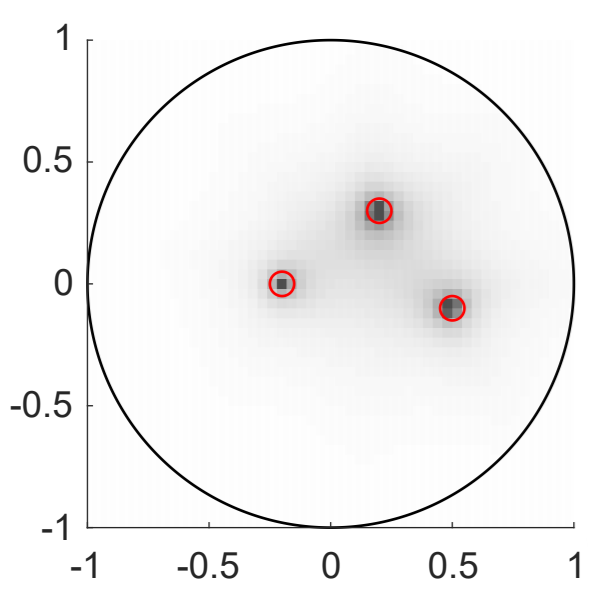}\\
\includegraphics[width=0.3\linewidth]{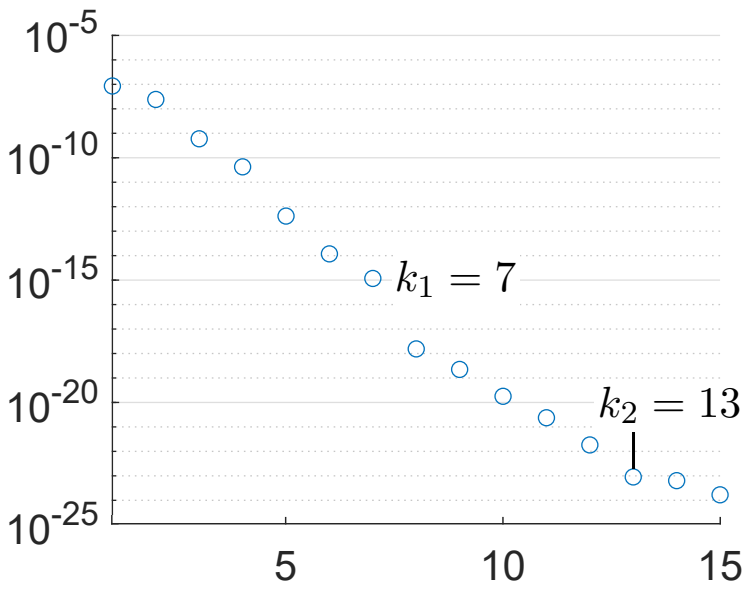}
&\includegraphics[width=0.25\linewidth]{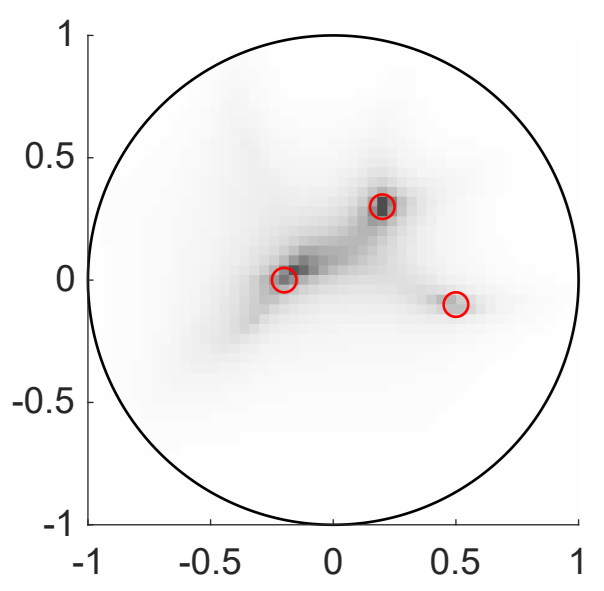}
&\includegraphics[width=0.25\linewidth]{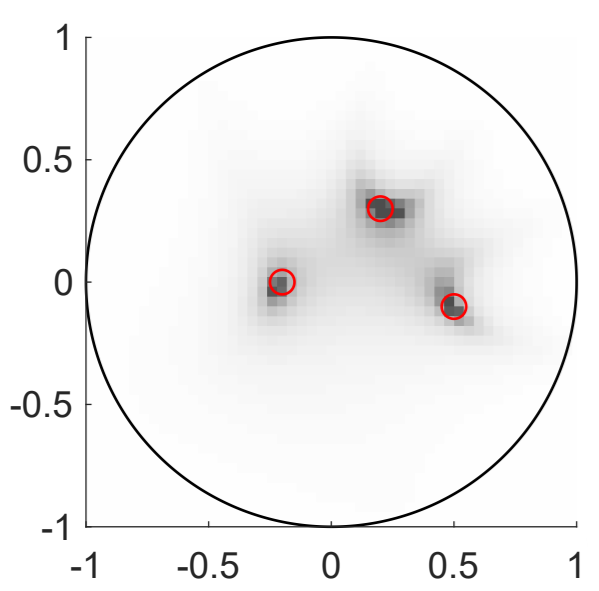}\\
\includegraphics[width=0.3\linewidth]{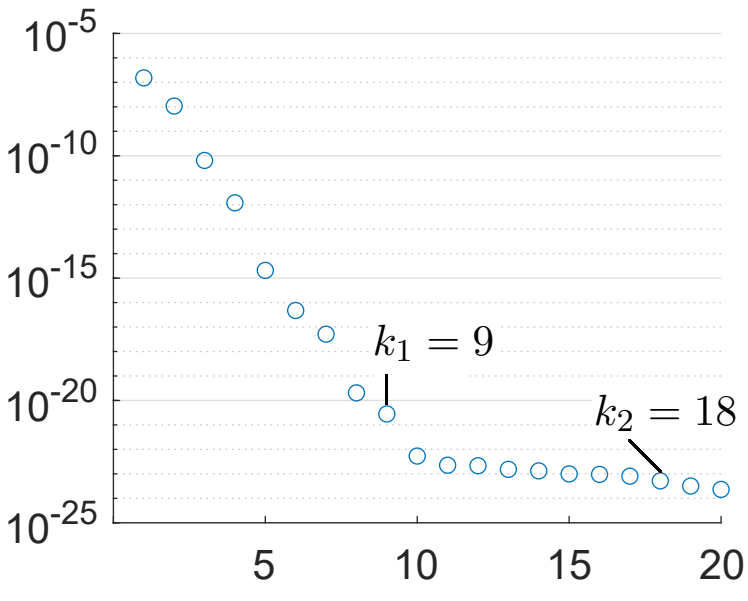}
&\includegraphics[width=0.25\linewidth]{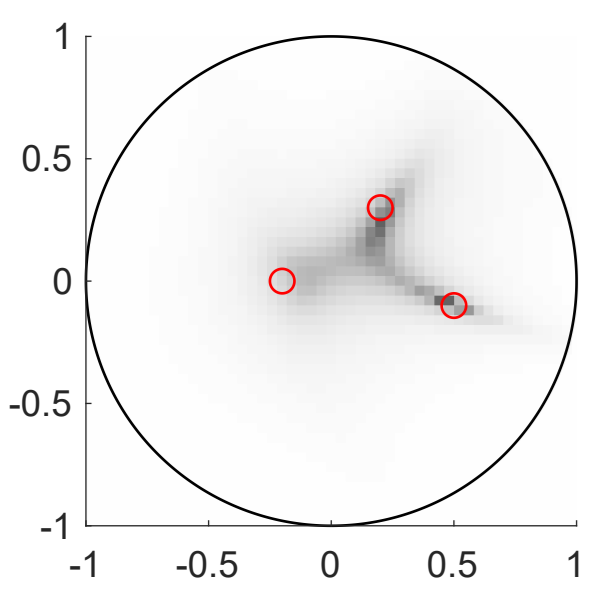}
&\includegraphics[width=0.25\linewidth]{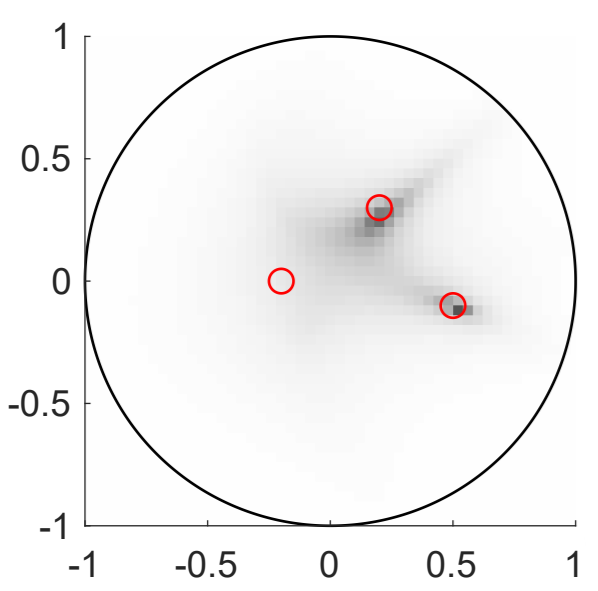}\\
(a) singular values of $B_I$ & (b) $k=k_1$ & (c) $k=k_2$
\end{tabular} 		\caption{\label{fig:threeinclusion} The reconstruction results for case (ii) in Example \ref{ex:2} with exact data. Panel (a) shows the singular values of $B_I$. Panels (b) and (c) show the contour plots of $W_{n,k}$ with $k=k_1$ and $k=k_2$, respectively. The red circles are the exact inclusion boundaries $\p A_1$, $\p A_2$ and $\p A_3$. The results from top to bottom correspond to the source configurations (i)--(iii) in Fig. \ref{fig:sourcepts}, respectively.}
\end{figure}

Figs. \ref{fig:twoinclusion}  and \ref{fig:threeinclusion} show the reconstruction results for Example \ref{ex:2} (i) and  (ii) with exact data, respectively. The singular values of $B_I$ (cf. \eqref{eq:datamatrix:realistic}) decay rapidly, indicating its low-rank structure. With a suitable truncation level $k$, the algorithm can accurately recover the locations, for all three configurations of observation apertures.
	
\begin{example} \label{ex:3}
Consider two circular inclusions $A=A_1\bigcup A_2$, where $A_{\ell} = \{ x \in \mathbb{R}^2: |x - z_{\ell}| < \varepsilon \}$ for $\ell = 1, 2$, with $\gamma_1 = \gamma_2 = 50$. We consider different inclusion sizes $\varepsilon$ and different locations $\{ z_{\ell} \}_{\ell = 1}^2$.
\end{example}

\begin{figure}[ht!]
\centering
\setlength{\tabcolsep}{0pt}
\begin{tabular}{ccccc}
\includegraphics[width=0.19\linewidth]{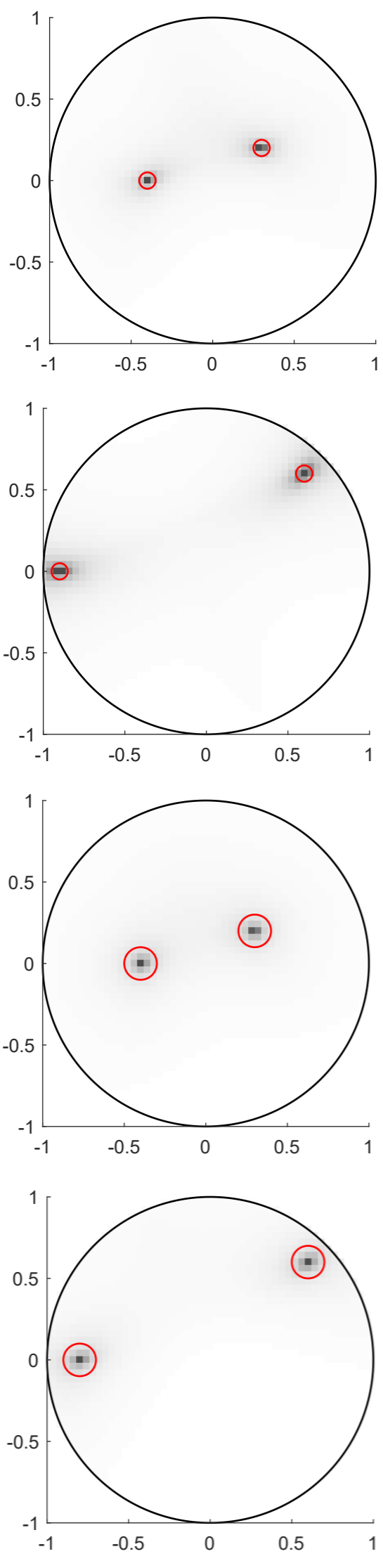}
&	\includegraphics[width=0.19\linewidth]{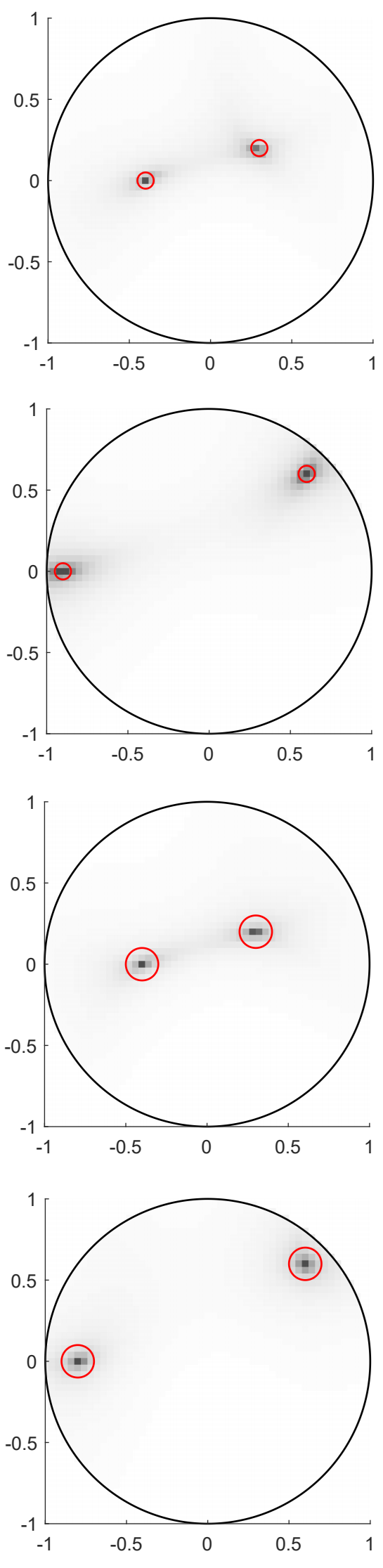}
&	\includegraphics[width=0.19\linewidth]{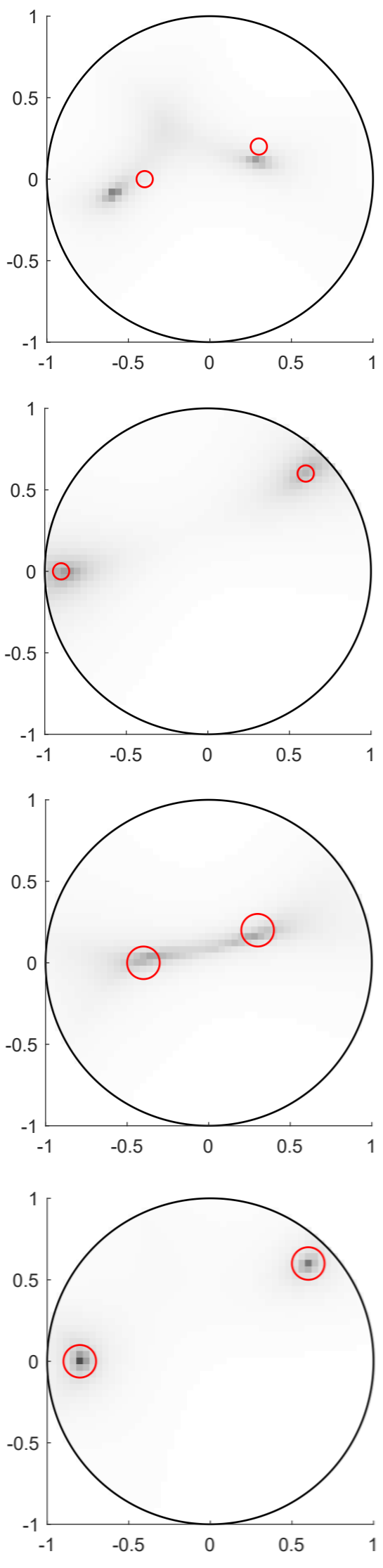}
&	\includegraphics[width=0.19\linewidth]{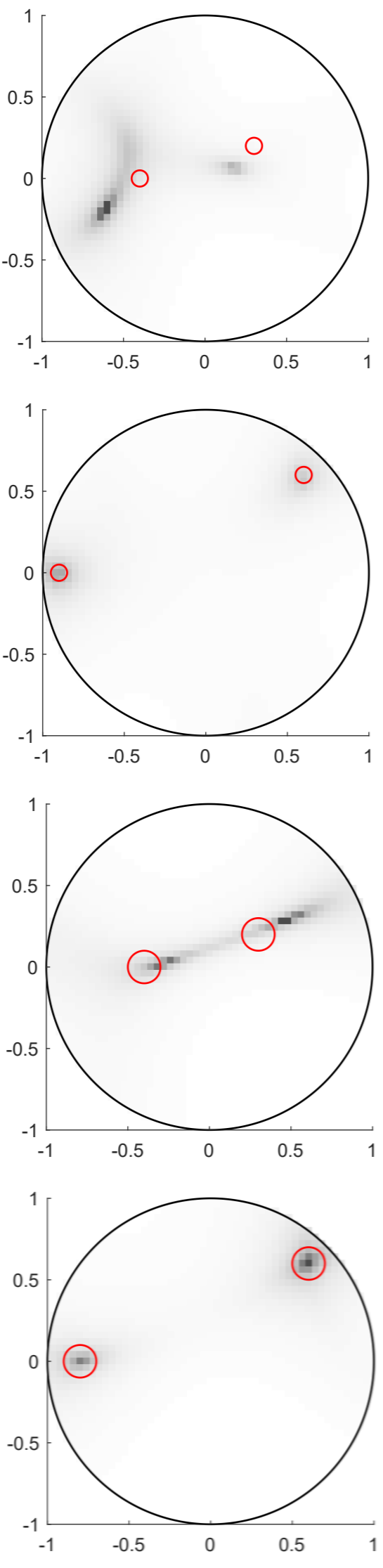}
&	\includegraphics[width=0.19\linewidth]{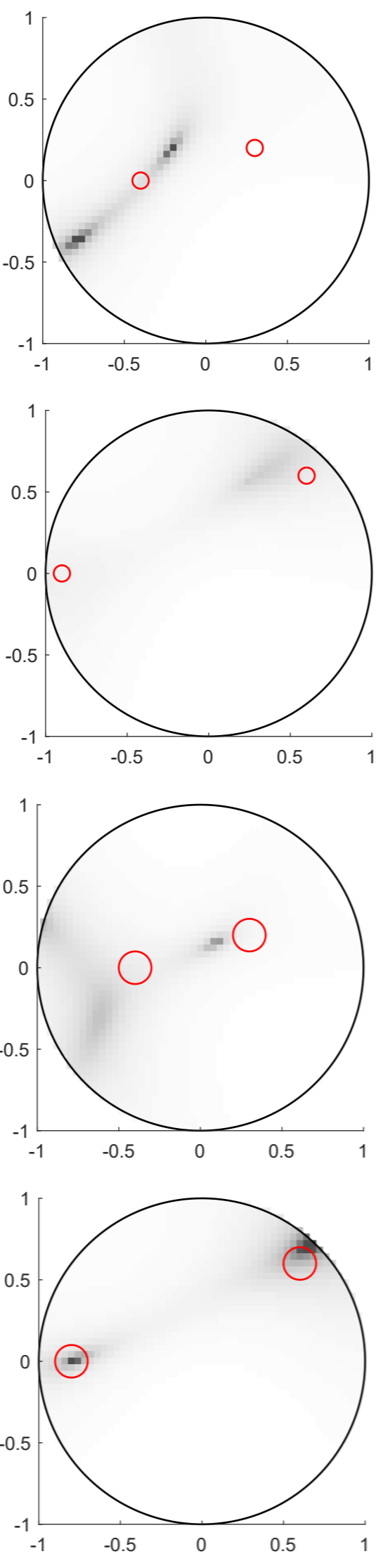}\\
$\sigma=0\%$&$\sigma=0.2\%$&$\sigma=1\%$&$\sigma=2\%$&$\sigma=5\%$
\end{tabular}
\caption{\label{fig:noisy:reconst} The contour plots of $W_{n,k}$ for Example \ref{ex:3} with various inclusion sizes $\varepsilon$, locations $\{ z_{\ell} \}_{\ell = 1}^2$ and noise levels $\sigma$. The truncation level $k$ in $W_{n,k}$ is $5$ for $\sigma = 0\%$, and varies within $[5, 9]$ for $\sigma > 0$ depending on the decay of singular values of $B_I$. The red circles are the exact inclusion boundaries $\{ \partial A_{\ell} \}_{\ell = 1}^2$.}
\end{figure}

In Fig. \ref{fig:noisy:reconst}, we present the reconstruction results for Example \ref{ex:3}, with the source points in Configuration (i) of Fig. \ref{fig:sourcepts}. With a fixed level of noise, the recovery of the inclusion locations is more accurate either as the size $\varepsilon$ of the inclusion increases or as the inclusions get closer to the boundary $\partial\Omega$ of the domain.
This concurs with intuition that the inclusions with larger volume fractions or closer to the boundary generate stronger effective signals on the boundary $\partial\Omega$ and thus are easier to locate.

\begin{example}\label{ex:many:ellipses}
Consider $m$ elliptical inclusions $A = \bigcup_{\ell=1}^m A_{\ell}$, with $A_{\ell} = \{(a,b)\in\mathbb{R}^2\,:\,\varrho^{-1}(a - z_{\ell, 1})^2+\varrho(b-z_{\ell, 2})^2<(0.05)^2\}$ and $\gamma_{\ell}=50$ for $\ell=1,\dots,m$, where the aspect ratio $\varrho\in\{2,3,5\}$, in two cases: {\rm(i)} $m=2$ with $(z_{1, 1},z_{1, 2})=(0.3,0.2)$ and $(z_{2, 1},z_{2, 2})=(-0.4,0)$ and {\rm(ii)} $m=3$ with $(z_{1, 1},z_{1, 2})=(-0.2,0)$, $(z_{2, 1},z_{2, 2})=(0.2,0.3)$ and $(z_{3, 1},z_{3, 2})=(0.5,-0.1)$.
\end{example}

\begin{figure}[hbt!]
\centering\setlength{\tabcolsep}{0pt}
\begin{tabular}{ccc}
\includegraphics[width=0.25\linewidth]{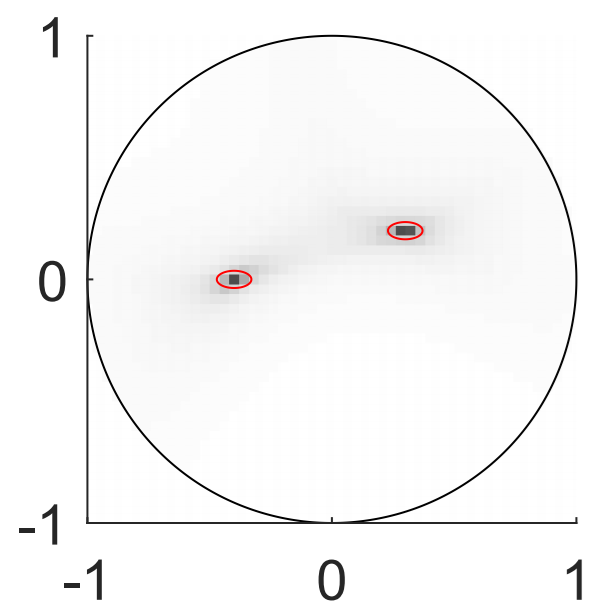}
&\includegraphics[width=0.25\linewidth]{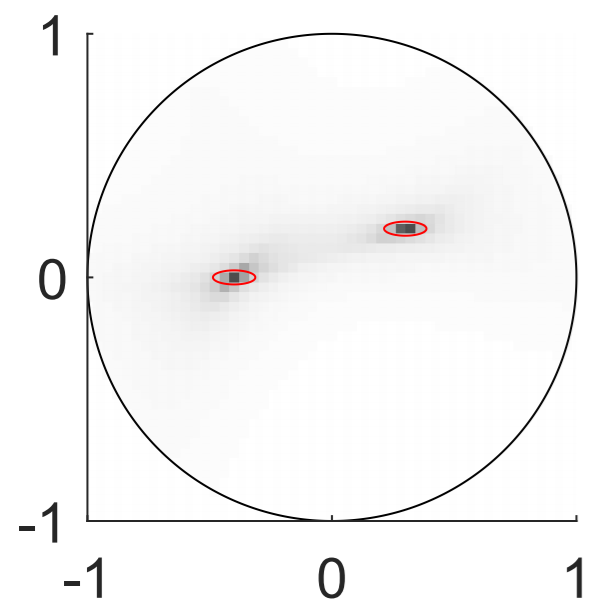}
&\includegraphics[width=0.25\linewidth]{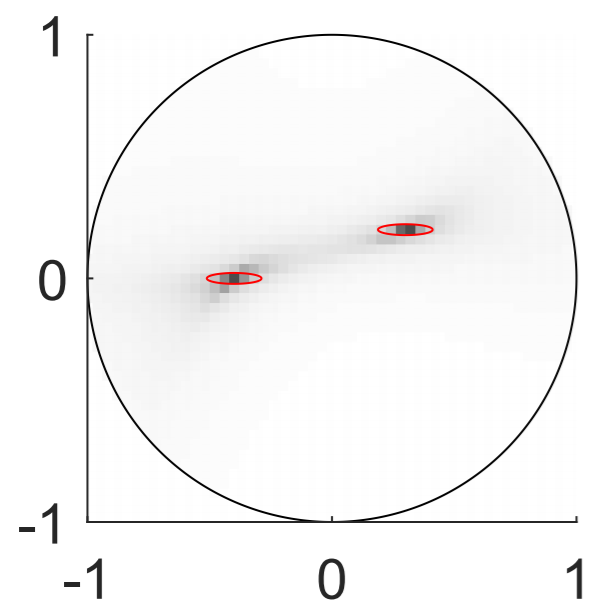}\\
\includegraphics[width=0.25\linewidth]{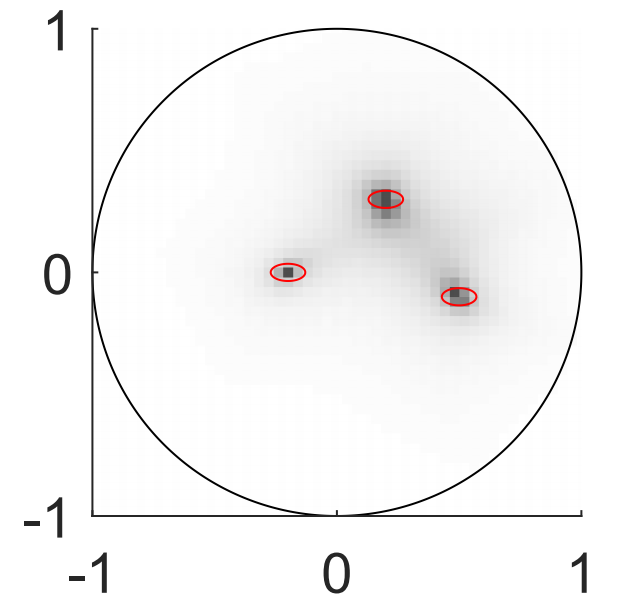}
&\includegraphics[width=0.25\linewidth]{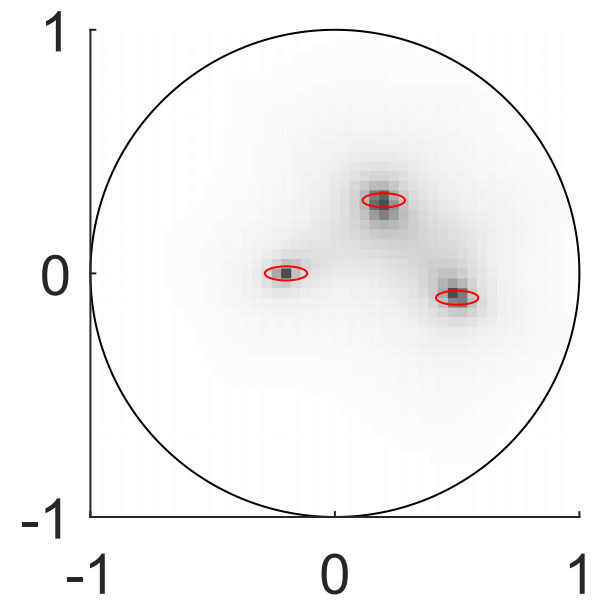}
&\includegraphics[width=0.25\linewidth]{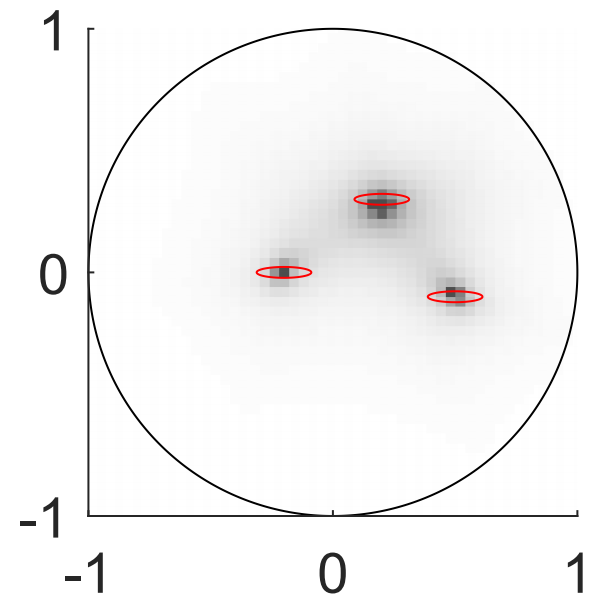}\\
(a) $\varrho=2$ & (b) $\varrho=3$ & (c) $\varrho=5$
\end{tabular} 		
\caption{\label{fig:many:ellipses}The contour plots of $W_{n,k}$ for Example \ref{ex:many:ellipses}. The top and bottom panels show the reconstruction results for cases (i) and (ii), respectively. The truncation level in $W_{n,k}$ is $k  =5$ for case (i) and $k = 7$ for case (ii). The red circles are exact inclusion boundaries $\{\p A_j\}_{j=1}^m$.}
\end{figure}

In Fig. \ref{fig:many:ellipses}, we show the reconstructions for Example \ref{ex:many:ellipses} with exact data.
The results indicate that the reconstructions remain fairly accurate for all the cases, even though the inclusions are not circular any more. Thus, the algorithm in Section \ref{ssec:multiple} is also applicable to general inclusions other than circular ones.

\section{Conclusion}
\label{sec: conclusion}

In this work we have developed two direct algorithms for recovering the locations of small conductivity inclusions in the subdiffusion model from boundary measurement. The algorithms are based on the asymptotic expansion of the boundary measurement with respect to the size of small inclusions and essentially make use of approximate fundamental solutions. We have provided theoretical underpinnings for the algorithms for one and multiple circular inclusions. Several numerical experiments indicate that the algorithms can indeed accurately predict the locations of small circular and elliptical inclusions, including the case where data noise is present. Future works may include the analysis of the algorithms for more general conductivity inclusions as well as other types of measurements, e.g. partial boundary measurements or posteriori boundary measurements (i.e., over the time interval $[T_0,T]$, with $T_0>0$).

\appendix

\section{Proof of Theorem \ref{thm:multi-char}}
\label{app:thm:multi-char}

We prove the two assertions by means of contradiction. Suppose that the identity \eqref{eq:error:oneside} holds for some $z = \widehat{z} \in\Om\backslash Z$. Fix any $x_0\in\Gamma$, and let $\Gamma_{x_0,\widehat{z}}=\{x\in\Gamma: (x-\widehat{z})\cdot(x_0-\widehat{z})\ge0\}$ and $f_0\in L^2(\Gamma)$ be the characteristic function  $\chi_{\Gamma_{x_0,\widehat{z}}}(x)$ of the set $\Gamma_{x_0,\widehat{z}}$.
Since $\Gamma$ is a nonempty open set containing $x_0$, we have $\|f_0\|_{L^2(\Gamma)}>0$.
Choose any $g_0\in L^2(\Gamma)$, and let $\varepsilon>0$ be small enough so that $\widehat{z}\notin\overline{A}$.
Then $\mathcal{M}[g_0]$ is analytic at $\widehat{z}$, and hence
\begin{align}\label{eq: second relation}
		\lim_{\rho\to 0}\int_{B_{\rho}(\widehat{z})}\gamma_0 \Delta_y \mathcal{M}[g_0](y)F_{\widehat{z},f_0}(y)\,{\rm d}y=0.
	\end{align}
Meanwhile, by the definition of the operator $\mathcal{L}_{\widehat{z}}[f_0]$, there holds
\begin{equation}\label{eq:LapLz:action}
	\int_{B_{\rho}(\widehat{z})}\!\!\!\!\!\!\!\!\gamma_0\Delta_y\mathcal{L}_{\widehat{z}}[f_0](y)F_{\widehat{z},f_0}(y)\,{\rm d}y=\int_{B_{\rho}(\widehat{z})}\!\!\int_0^T\!\!\!\!\nabla_y(\p_t^\alpha\Psi_{(y,0)}(\widehat{z},t)\!-\!\delta_{(\widehat{z},0)}(y,t))\cdot \widetilde{F}(y,t)\,{\rm d}t\,{\rm d}y,	
\end{equation}
with
$$\widetilde{F}(y,t):=\int_\Gamma\nabla\Psi_{(x,0)}(\widehat{z},T-t)f_0(x)\,{\rm d}\sigma_x\, F_{\widehat{z},f_0}(y).
$$
Note that the auxiliary function $\widetilde{F}(y,t)$ is smooth in $(y,t)\in \Om\times [0,T]$ since $\operatorname{dist}(\widehat{z},\Gamma)>\operatorname{dist}(\Om,\Gamma)>0$
and satisfies $\widetilde{F}(y,0)=\nabla F_{\widehat{z},f_0}(y)$ for $y\in\Om$. Thus, we have
\begin{equation}\label{eq:ddelta:action}
	\int_{B_{\rho}(\widehat{z})}\int_0^T\nabla_y(-\delta_{(\widehat{z},0)}(y,t))\cdot \widetilde{F}(y,t){\rm d}t\,{\rm d}y=\Delta F_{\widehat{z},f_0}(\widehat{z}).	
	\end{equation}
To bound the difference between \eqref{eq:LapLz:action} and \eqref{eq:ddelta:action}, i.e., the integral $\int_{B_{\rho}(\widehat{z})}\!\!\int_0^T\!\!\nabla_y\p_t^\alpha\Psi_{(y,0)}(\widehat{z},t)\cdot \widetilde{F}(y,t){\rm d}t\,{\rm d}y$, we split it into \begin{equation}\label{eq:something:to:degenerate}
		\int_{B_{\rho}(\widehat{z})}\int_0^T\nabla_y\p_t^\alpha\Psi_{(y,0)}(\widehat{z},t)\cdot \widetilde{F}(y,t){\rm d}t\,{\rm d}y=J_1(\rho )+J_2(\rho)+J_3(\rho),
	\end{equation}
where
$$J_k(\rho):=\int_{B_{\rho }(\widehat{z})}\int_0^T\nabla_y\p_t^\alpha\Psi_{(y,0)}(\widehat{z},t)\cdot \widetilde{J}_k(y,t){\rm d}t\,{\rm d}y,\quad k=1,2,3,$$
with
$$\widetilde{J}_1(y,t):=\widetilde{F}(y,t)-\widetilde{F}(y,0),\quad \widetilde{J}_2(y,t):=\widetilde{F}(y,0)-\widetilde{F}(\widehat{z},0),\quad \widetilde{J}_3(y,t):=\widetilde{F}(\widehat{z},0).$$
Since $\nabla\Psi_{(y,0)}(\widehat{z},t)=-\nabla\Psi_{(2\widehat{z}-y,0)}(\widehat{z},t)$ for $y\in B_\rho(\widehat{z})$, we have $J_3(\rho)=0$. Since $\widetilde{F}$ is smooth, there hold \begin{equation}\label{eq:J1J2:estimate}
		|\widetilde{J}_1(y,t)|\le C t\quad\mbox{and}\quad |\widetilde{J}_2(y,t)|\le C|y-\widehat{z}|,\quad\forall(y,t)\in B_\rho(\widehat{z})\times[0,T].
	\end{equation}
Also from \cite[Theorem 2.4 (i) and Theorem 2.5 (i-ii)]{Kim:2016:ABF}, we have, for some $C$ and $C'>0$,
	\begin{align}\label{eq:ddpsi:estimate}
		&|\nabla_y\p_t^\alpha\Psi_{(y,0)}(\widehat{z},t)|
  \le\begin{cases}
			C |y-\widehat{z}|^{-d+1}t^{-2\alpha},&\mbox{if }|y\!-\!\widehat{z}|^2 t^{-\alpha}\le1,\\
			C|y-\widehat{z}|^{-d-1}t^{-\alpha}e^{-C'|y-\widehat{z}|^{2/(2-\alpha)} t^{-\alpha/(2-\alpha)}},&\mbox{if }|y\!-\!\widehat{z}|^2 t^{-\alpha}\ge1.
		\end{cases}
	\end{align}
By \eqref{eq:J1J2:estimate}, \eqref{eq:ddpsi:estimate} and changing variables from $(y,t)$ to $(w,\tau)=(y-\widehat{z},|y-\widehat{z}|^2 t^{-\alpha})$, we derive
	\begin{align*}
	&|J_1(\rho)|\le C\int_{B_{\rho}(\widehat{z})}\int_{0}^T \left( |y-\widehat{z}|^{-d+1}t^{-2\alpha}\chi_{\tau\le1}+
	\frac{e^{-C'|y-\widehat{z}|^{2/(2-\alpha)} t^{-\alpha/(2-\alpha)}}}{|y-\widehat{z}|^{d+1}t^{\alpha}}\chi_{\tau\ge1}\right)t\,{\rm d}t\,{\rm d}y\\
	&=C\alpha^{-1}\int_{B_{\rho}(0)}|w|^{-d-3+\frac{4}{\alpha}}\left[\int_{|w|^2T^{-\alpha}}^1 \tau^{1-\frac{2}{\alpha}}\,{\rm d}\tau + \int_1^\infty \tau^{-\frac{2}{\alpha}}e^{-C'\tau^{1/(2-\alpha)}}\,{\rm d}\tau\right]{\rm d}w
   = o(\rho),
	\end{align*}
as $\rho\to0$.	Similarly, we derive
	\begin{align*}
		|J_2(\rho)|&\le C\alpha^{-1}\int_{B_{\rho}(0)}|w|^{-d-2+\frac{2}{\alpha}}\left[\int_{|w|^2T^{-\alpha}}^1 \tau^{1-\frac{1}{\alpha}}\,{\rm d}\tau + \int_1^\infty \tau^{-\frac{1}{\alpha}}e^{-C'\tau^{1/(2-\alpha)}}\,{\rm d}\tau\right]{\rm d}w
  =o(\rho),
	\end{align*}
as $\rho\to0$.	Therefore, the right-hand side of \eqref{eq:something:to:degenerate} converges to $0$ as $\rho\to0$.
This and the estimates \eqref{eq:LapLz:action}--\eqref{eq:ddelta:action} gives
	\begin{equation}\label{eq: first relation}
		\int_{B_{\rho}(\widehat{z})}\gamma_0\Delta_y\mathcal{L}_{\widehat{z}}[f_0](y)F_{\widehat{z},f_0}(y)\,{\rm d}y=\Delta F_{\widehat{z},f_0}(\widehat{z})+o(\rho).
	\end{equation}
Combining \eqref{eq: first relation} with \eqref{eq: second relation} gives
\begin{equation*} \lim_{\rho\to0}\int_{B_{\rho}(\widehat{z})}\!\!\!\!\gamma_0\Delta_y \left(\mathcal{L}_z[f_0](y)\!-\!\mathcal{M}[g_0](y)\right) F_{\widehat{z},f_0}(y){\rm d}y=\Delta F_{\widehat{z},f_0}(\widehat{z})=\left|\int_\Gamma \!\! \nabla\Psi_{(x,0)}(\widehat{z},T)f_0(x){\rm d}\sigma_x\right|^2.
\end{equation*}
Then by the definition of $\Psi$, we have
\begin{align*}
	&\left|\frac{x_0-\widehat{z}}{|x_0-\widehat{z}|}\cdot\int_\Gamma \nabla\Psi_{(x,0)}(\widehat{z},T)f_0(x)\,{\rm d}\sigma_x\right|
	=\int_{\Gamma_{x_0,\widehat{z}}}\!\!\!\! \frac{(x-\widehat{z})\cdot(x_0-\widehat{z})}{|x-\widehat{z}||x_0-\widehat{z}|}|\nabla\Psi_{(x,0)}(\widehat{z},T)|\,{\rm d}\sigma_x\ge C>0,
\end{align*}
where $C$ is independent of $g_0$ and $\varepsilon$. This contradicts the assumption \eqref{eq:error:oneside}. One can derive \eqref{eq:compare} in a similar way to \eqref{eq: first relation}.

\bibliographystyle{siam}
\bibliography{reference}
\end{document}